\documentclass{amsart}
\usepackage{amscd}
\usepackage[all]{xy}
\numberwithin{equation}{section}
\usepackage{amssymb}
\let\cal\mathcal
\def\Ascr{{\cal A}}

\def\Cscr{{\cal C}}
\def\Dscr{{\cal D}}
\def\Escr{{\cal E}}

\def\Gscr{{\cal G}}
\def\Hscr{{\cal H}}
\def\Iscr{{\cal I}}

\def\Lscr{{\cal L}}

\def\Oscr{{\cal O}}
\def\Pscr{{\cal P}}
\def\Qscr{{\cal Q}}

\def\Tscr{{\cal T}}
\def\Uscr{{\cal U}}
\def\Vscr{{\cal V}}
\def\Wscr{{\cal W}}
\def\Xscr{{\cal X}}

\let\blb\mathbb
\def\CC{{\blb C}}

\def \AA{{\blb A}}
\def \ZZ{{\blb Z}}

\def \NN{{\blb N}}

\def\cosk{\operatorname{cosk}}
\def\sk{\operatorname{sk}}
\def\Pre{\operatorname{Pre}}
\def\Sh{\operatorname{Sh}}
\def\Bilin{\operatorname{Bilin}}
\def\Fun{\operatorname{Fun}}
\def\Ho{\operatorname{Ho}}
\def\pro{\operatorname{pro}}

\def\Id{\operatorname{id}}

\def\Der{\operatorname{Der}}

\def\ctimes{\mathbin{\hat{\otimes}}}
\def\cbtimes{\mathbin{\hat{\boxtimes}}}

\def\Mod{\operatorname{Mod}}

\def\Ch{\mathop{\mathrm{Ch}}}

\def\Qch{\operatorname{Qch}}

\def\Spec{\operatorname {Spec}}

\def\Ext{\operatorname {Ext}}
\def\Hom{\operatorname {Hom}}

\def\RHom{\operatorname {RHom}}

\def\Open{\operatorname {Open}}

\def\im{\operatorname {im}}

\def\coker{\operatorname {coker}}
\def\ker{\operatorname {ker}}

\def\tot{\operatorname {tot}}

\def\r{\rightarrow}
\def\l{\leftarrow}

\DeclareMathOperator{\Proj}{Proj}
\DeclareMathOperator{\Ab}{Ab}
\DeclareMathOperator{\Pro}{Pro}

\DeclareMathOperator{\Alg}{Alg}

\DeclareMathOperator{\Aut}{Aut}

\DeclareMathOperator{\Ob}{Ob}

\let\dirlim\injlim
\let\invlim\projlim

\newtheorem{lemma}{Lemma}[section]
\newtheorem{proposition}[lemma]{Proposition}
\newtheorem{theorem}[lemma]{Theorem}

\newtheorem{lemmas}{Lemma}[subsection]
\newtheorem{propositions}[lemmas]{Proposition}
\newtheorem{theorems}[lemmas]{Theorem}
\newtheorem{corollarys}[lemmas]{Corollary}

\newtheorem{conventions}[lemmas]{Convention}
\newtheorem{conventionwarnings}[lemmas]{Convention-Warning}

\theoremstyle{definition}

\newtheorem{examples}[lemmas]{Example}
\newtheorem{definitions}[lemmas]{Definition}

{

\newtheorem{step}{Step}

}

\theoremstyle{remark}

\newtheorem{remark}[lemma]{Remark}
\newtheorem{remarks}[lemmas]{Remark}

\newdimen\uboxsep \uboxsep=1ex
\def\uboxn#1{\vtop to 0pt{\hrule height 0pt depth 0pt\vskip\uboxsep
\hbox to 0pt{\hss #1\hss}\vss}}

\def\uboxs#1{\vbox to 0pt{\vss\hbox to 0pt{\hss #1\hss}
\vskip\uboxsep\hrule height 0pt depth 0pt}}

\def\Spc{\operatorname{Spc}}
\def\Adic{\operatorname{Adic}}

\def\Set{\operatorname{Set}}
\def\Sch{\operatorname{Sch}}
\def\Jet{\operatorname{Jet}}
\def\poly{\operatorname{poly}}
\def\coord{\operatorname{coord}}
\def\cont{\operatorname{cont}}
\def\poly{\operatorname{poly}}
\def\aff{\operatorname{aff}}
\def\Diff{\operatorname{Diff}}
\def\HC{\mathop{\mathbf{C}}\nolimits}
\def\Gl{\operatorname{GL}}

\def\an{\operatorname{an}}
\keywords{Deformation quantization, Abelian categories}

\author{Michel Van den Bergh}
\address{Departement WNI\\Universiteit Hasselt\\ Universitaire Campus\\ Building D\\ 3590 
Diepenbeek\\ Belgium} 
\thanks{The author is a director of research at the FWO} 
 \email[M. Van den
Bergh]{michel.vandenbergh@uhasselt.be} 
\title{On global deformation quantization in the algebraic case}
\begin{document}
\begin{abstract}
  We give a proof of Yekutieli's global algebraic deformation
  quantization result which does not rely on the choice of local
  sections of the bundle of affine coordinate systems.  Instead we use
  an argument inspired by algebraic De Rham cohomology. 
\end{abstract}
\maketitle
\tableofcontents
\section{Introduction and motivation}
This paper has grown out of an attempt to understand \cite{Ko10,ye3}.
These papers deal with deformation quantization in an algebraic
setting.  After some consideration
we decided no harm would be done by writing down our own account. 
For simplicity we restrict ourselves to infinitesimal deformations.
The extension to formal deformations is routine. 

\medskip

Kontsevich's fundamental idea is that quantization of Poisson brackets
should take place in the setting of \emph{twisted presheaves}.  To
explain this let $X$ be a separated quasi-compact scheme over a field
$k$.  Choose an affine covering $\Uscr=\{U_1,\ldots,U_n\}$ of $X$. For
$J\subset \{1,\ldots,n\}$ define $U_J=\cap_{j\in J} U_j$.  A twisted
presheaf of $k$-algebras on $\Uscr$ is a collection of $k$-algebras
$\Ascr(U_J)$ together with restriction maps $\rho_{J,J'}:\Ascr(U_J)\r
\Ascr(U_{J'})$ for $J\subset J'$ which are compatible with
compositions up to an explicit inner automorphism. In addition the
units defining these inner automorphisms should satisfy a natural
cocycle conditions for triple compositions.  The motivation is that
under suitable flatness conditions one may define a category of
quasi-coherent sheaves over a twisted presheaf.

\medskip

Assume now that $X$ is a smooth and separated over $k=\CC$. For
$(l,m)$ an Artinian local $k$-algebra with residue field $k$,
Kontsevich constructs in \cite{Ko10} a ``quantization'' arrow
\begin{equation}
\label{ref-1.1-0}
\begin{gathered}
  \{\text{Poisson brackets on $X$ with coefficients in $m$}\}/\cong\\
\r \\
 \{\text{Flat $l$-deformations of $\Oscr_X$ in the category of twisted presheaves}\}/\cong
\end{gathered}
\end{equation}
To make the connection with the deformation theory of abelian
categories \cite{lowenvdb1,lowenvdb2} note that in
\cite{lowen5} Tor Lowen constructs  a natural bijection
\begin{gather*}
 \{\text{Flat $l$-deformations of $\Oscr_X$ in the category of twisted presheaves}\}/\cong\\
\leftrightarrow\\
 \{\text{Flat $l$-deformations of $\Qch(\Oscr_X)$}\}/\cong
\end{gather*}
Our aim is to explain \eqref{ref-1.1-0}. This explanation will make
it clear what the obstruction is against reversing the arrow.

\medskip

Let us recall some basic constructions. The complex of sheaves of
poly-differential operators $\Dscr^{\poly,\cdot}_X$ is defined as
follows. For an open $U$ of $X$, $\Dscr^{\poly,n}_X(U)$ is given by
the multilinear maps $\Oscr_U^{\otimes n} \r \Oscr_U$ which are
differential operators in each argument.\footnote{Note that we do not include
the usual shift by one (e.g. \cite{Ko3}) in the  definition of the complex of poly-differential operators.}

$\Dscr^{\poly,\cdot}_X$ is 
equipped with the standard Hochschild differential and
Gerstenhaber bracket. In this way $\Dscr^{\poly}_X[1]$ becomes a sheaf of
DG-Lie algebras.

Likewise $\Tscr^{\poly,\cdot}_X$ is defined as the graded sheaf on $X$
whose sections of degree $n$ on an open $U$ are given by the
multilinear maps $\Oscr_U^{\otimes n} \r \Oscr_U$ which are fully
anti-symmetric and derivations in each argument.  When equipped with the
Schouten-Nijenhuis bracket and trivial  differential
$\Tscr^{\poly}_X[1]$ also becomes a DG-Lie algebra.

\medskip

The key result in algebraic deformation quantization is the following
\begin{theorem} \cite[Thm 0.2]{ye3}  \label{ref-1.1-1}
There is an isomorphism 
\begin{equation}
\label{ref-1.2-2}
\Tscr^{\text{poly},\cdot}_X [1]\cong  \Dscr^{\text{poly},\cdot}_X[1]
\end{equation}
in the homotopy category of sheaves of DG-Lie algebras. 
Furthermore
if $X$ has a system of parameters $(x_i)_i$ then the resulting map 
on homology
\[
\Tscr^{\text{poly},\cdot}_X\r H^\cdot(\Dscr^{\text{poly},\cdot}_X)
\]
is given by the HKR-formula.
\[
\partial_{i_1}\wedge \cdots \wedge \partial_{i_n}\mapsto
\frac{1}{n!}\sum_{\sigma \in S_n}
(-1)^\sigma\partial_{i_{\sigma(1)}}\otimes \cdots \otimes \partial_{i_{\sigma(n)}}
\]
where $\partial_i=\partial_{x_i}$.
\end{theorem}
We will give a self contained proof of this result in this paper.  We
mimic to some extent Yekutieli's arguments and we use many of his
technical contributions.  However there is a substantial
simplification that we do not have to choose (local) sections of the
bundle of affine coordinate systems and thus we avoid the machinery of
simplicial sections.  Instead we use an argument inspired by algebraic
De  Rham cohomology (see \S\ref{ref-6.6-73} below). We believe
this idea is new even in the classical case of Fedosov quantization.
It follows in particular that \eqref{ref-1.2-2} is compatible with
automorphisms of $X$.

\medskip

An analogue of Theorem \ref{ref-1.1-1} has been proved in the complex
analytic case in \cite{cdh}. The proof in loc.\ cit.\ does not
extend immediately to the algebraic case as it depends on the choice
of a global connection. In \cite{DTT} the authors prove a version of
Theorem \ref{ref-1.1-1} using operadic methods which is uniformly
valid for the $C^\infty$, algebraic and complex analytic cases.
In \cite{CVdB} we will give yet another approach to these
results.

\medskip

Let us now explain how Proposition \ref{ref-1.1-1} is relevant to
deformation theory. For a sheaf of DG-Lie algebras $\Gscr$ on a
topological space $X$ one may take its derived global sections
$R\Gamma(X,\Gscr)^{\text{tot}}$ which is also a DG-Lie algebra and
which is canonically quasi-isomorphic to $R\Gamma(X,\Gscr)$ as complexes.
In particular the formation of $R\Gamma(X,\Gscr)^{\text{tot}}$ is
compatible with quasi-isomorphisms.  One possible construction using
pro-hypercoverings is outlined in Appendix \ref{ref-B-141}. A
different construction follows from Hinich's model structure on the
category of presheaves of algebras over an operad \cite{hinich2}. See
\S\ref{ref-B.10-198}. If $\Gscr$ is quasi-coherent and $X$ is separated
then $R\Gamma(X,\Gscr)^{\tot}$ is given by applying the Thom-Sullivan
normalization (see Appendix \ref{ref-A-137}) to the cosimplicial
DG-Lie algebra associated to an affine covering of $X$ (see
\cite{hsI,hsII}).  We note that it will be clear below that only the
properties of the functor $R\Gamma(X,-)^{\text{tot}}$ matter, not its
actual construction.

\medskip

Applying $R\Gamma(X,-)^{\tot}$ to \eqref{ref-1.2-2} we obtain in
particular an isomorphism in the homotopy category of DG-Lie algebras
\begin{equation}
\label{ref-1.3-3}
R\Gamma(X,\Tscr^{\text{poly},\cdot}_X [1])^{\tot}\cong R\Gamma(X,\Dscr^{\text{poly},\cdot}_X[1])^{\tot}
\end{equation}

Let $\frak{g}$ be a DG-Lie algebra.  The
\emph{Maurer-Cartan equation} in $m\otimes_k\frak{g}_1$ is given by
\begin{equation}
\label{ref-1.4-4}
d\pi+\frac{1}{2}[\pi,\pi]=0
\end{equation}
There is a a natural action on the solutions of this equation by the
``gauge'' group $\exp(m\otimes \frak{g}_0)$. It is well known that the
set of equivalence classes of solutions to the Maurer-Cartan equation 
is invariant under quasi-isomorphisms.

\medskip

If $A$ is a $k$-algebra then it is well known that the flat
$l$-deformations of $A$ correspond to solutions of the Maurer-Cartan
equation in $m\otimes_k\bold{C}(A)[1]$
where $\bold{C}(A)$ is the Hochschild complex
 of $A$ (equipped with the Gerstenhaber bracket).
Similar results for abelian and linear categories were proved in
\cite{lowenvdb1,lowenvdb2}.

\medskip

Let $\Uscr$ be as above and let $\frak{u}$ be the linear category with
objects  $\emptyset \subsetneq J\subset \{1,\ldots,n\}$ and 
\[
\frak{u}(J,J') = \begin{cases} \Hom_{\Oscr_X}(j_\ast \Oscr_{U_J},j_\ast \Oscr_{U_{J'}})=\Oscr_X(U_{J'}) &\text{if }\,J \subset J' \\ 0 & 
\text{otherwise}
\end{cases}
\]
where $j:U_J\r X$ denotes the inclusion map. We prove (Theorem
\ref{ref-3.1-13} below)
\begin{equation}
\label{ref-1.5-5}
R\Gamma(X,\Dscr^{\text{poly},\cdot}_X[1])^{\tot}\cong \bold{C}(\frak{u})[1]
\end{equation}
Since deformations of $\frak{u}$ are readily seen to correspond to deformations of $\Oscr_X$ as twisted
presheaf, and vice versa (see \cite{lowen5}), it follows from \eqref{ref-1.3-3} and 
\eqref{ref-1.5-5}
 that we have bijections.
\begin{equation}
\label{ref-1.6-6}
\begin{gathered}
  \{\text{Solutions to the MC equation in $R\Gamma(X,\Tscr_{X}^{\poly,\cdot}[1])^{\tot}\otimes_k m$}\}/\cong\\
  \leftrightarrow\\
  \{\text{Solutions to the MC equation in $R\Gamma(X,\Dscr_{X}^{\poly,\cdot}[1])^{\tot}\otimes_k m$}\}/\cong\\
  \leftrightarrow\\
 \{\text{Flat $l$-deformations of $\Oscr_X$ in the category of twisted presheaves}\}/\cong
\end{gathered}
\end{equation}
By Proposition \ref{ref-B.8.1-194} below there is a 
canonical map
\begin{equation}
\label{ref-1.7-7}
\Gamma(X,\Tscr_{X}^{\poly,\cdot})[1]\r R\Gamma(X,\Tscr_{X}^{\poly,\cdot}[1])^{\tot}
\end{equation}
and this is an isomorphism provided
\begin{equation}
\label{ref-1.8-8}
H^i(X,\wedge^j \Tscr_X)=0\qquad \text{for $i>0$}
\end{equation}
The solutions to the Maurer-Cartan equation in
$\Gamma(X,\Tscr_{X}^{\poly,\cdot})$ are the global Poisson brackets on
$X$. Thus combining \eqref{ref-1.6-6} with \eqref{ref-1.7-7} we
now  obtain the arrow \eqref{ref-1.1-0} and we see that
it is a bijection if \eqref{ref-1.8-8} holds. 

\medskip

Using similar ideas (see Proposition \ref{ref-4.1-24}) one proves that if
$X$ is proper there is a bijection 
\begin{gather*}
  \{\text{Solutions to the MC equation in $\oplus_{i,j} \Gamma(X^{\an},\Tscr^{i,0}_{X^{\an}}\otimes_{\Oscr_{X^{\an}}^\infty} \Omega^{0,j}_{X^{\an}})[1]\otimes_k m$}\}/\cong\\
  \leftrightarrow \\
  \{\text{Flat $l$-deformations of $\Oscr_X$ in the category of
    twisted presheaves}\}/\cong
\end{gather*}
thereby making the connection with the work of Barannikov and Kontsevich
\cite{BaKo}. The last bijection allows one, through the work of
Gualtieri \cite[\S 5.3]{gualtieri}, to associate a category of
coherent sheaves to an infinitesimal deformation of $X^{\an}$ as
\emph{generalized complex manifold}. Generalized complex manifolds
form a common generalization of complex and symplectic manifolds and
as such are important for mirror symmetry.  It is not known in general
how to define a (derived?)  category of coherent sheaves over a
generalized complex manifold. In the case of a symplectic manifold
this should be some variant of the Fukaya category.

\medskip

Other papers relevant for algebraic deformation quantization are
\cite{BezKal,BezKal1,schapira}. \cite{BezKal1} is especially
interesting as it discusses Fedosov quantization in positive
characteristic. This falls totally outside the reach of methods based
on DG-Lie algebras and the Maurer-Cartan equation.

\medskip

We now give a quick outline of the current paper. In \S\ref{ref-3-9} we
explain the connection between poly-differential operators and the Hochschild
complex of schemes. In \S\ref{ref-4-22} we discuss the application
to the analytic case mentioned above. 

The proof of  Theorem \ref{ref-1.1-1} uses crucially infinite dimensional
formal schemes. We discuss the relevant topological notions in
\S\ref{ref-5-27}. 

In \S\ref{ref-6-35} we use formal schemes to give 
an account of formal geometry in the algebraic case. See also \cite[\S5]{ye3}. 
Theorem \ref{ref-6.6.1-74} is our crucial acyclicity result for the bundle
of affine coordinate systems. 

In \S\ref{ref-7-82} we present a reminder on DG-Lie and
$L_\infty$-algebras.  An important notion is the twist of an
$L_\infty$-morphism by a solution of the Maurer-Cartan equation (which
I learnt from Yekutieli). We also discuss descent for
$L_\infty$-morphisms under an algebraic group action and its
compatibility with twisting. This is used to descend constructions on
the bundle of local coordinate systems to the bundle of affine local
coordinate systems.

In \S\ref{ref-8-103} we have a new look at poly-differential
operators and poly-vector fields. We introduce Kontsevich's local
$L_\infty$-quasi-isomorphism and remind the reader of its properties.
An interesting remark is that the linearity property (P3) thought to
be essential for globalization actually follows from (P5).

Finally in \S\ref{ref-9-111} we prove Theorem \ref{ref-1.1-1}.

\medskip

\noindent
\textbf{Acknowledgment\ } The author is very grateful to Damien
Calaque, Gilles Halbout, Vladimir Hinich, Bernhard Keller, Tor Lowen
and Amnon Yekutieli for many helpful discussions.
\section{Notations and conventions}
For simplicity of exposition we assume throughout that our base field
is algebraically closed of characteristic zero (and usually $\CC$). It
is clear that with sufficient care one can get by with weaker
hypotheses.

Many of the objects we use are equipped with some
kind of topology, but if an object is introduced without a specified
topology we assume that it is equipped with the discrete topology.

If an object carries a natural grading then all constructions associated
to it are implicitly performed in the graded context.  This implies
in particular to completions. 
\section{Going from poly-differential operators to the
  Hochschild complex}
\label{ref-3-9}
The main result of this section (Theorem \ref{ref-3.1-13}) was used in the introduction.

Let $k$ be a field.  If $\frak{u}$ is a $k$-linear category then the
Hochschild complex of $\frak{u}$ is defined as
\begin{equation}
\label{ref-3.1-10}
\bold{C}^n(\frak{u})=\prod_{U_0,\dots,U_p\in \Ob(\frak{u})}\Hom_k(\frak{u}(U_{p-1},U_p)
\otimes_k \dots \otimes_k \frak{u}(U_0,U_1), \frak{u}(U_0,U_p))
\end{equation}
with the usual differential.

It is well-known that the Hochschild complex of a linear category has
a lot of ``higher structure''. In particular it is a DG-Lie algebra
when equipped with the Gerstenhaber bracket. This is the structure
we will use below. 

The Hochschild complex of a linear category is contravariantly functorial
for fully faithful functors $\frak{v}\r \frak{u}$. The resulting map
$\bold{C}(\frak{u})\r \bold{C}(\frak{v})$ will be called the restriction map.

\medskip

Assume now that $k$ has characteristic zero and let $X$ be a smooth
separated scheme over $k$. It will be convenient to use the notations
$D^{\poly,\cdot}(U)=\Dscr^{\poly,\cdot}_X(U)$ for $U\subset X$ open and
$D^{\poly,\cdot}(R)=D^{\poly,\cdot}(X)$ for $X=\Spec R$, and similarly for
poly-vector fields.

\medskip

Assume first
that $X=\Spec R$ is affine. We obtain an inclusion of complexes
\begin{equation}
\label{ref-3.2-11}
D^{\poly,\cdot}(R)[1]\r \bold{C}(R)[1]
\end{equation}
which is compatible with the DG-Lie algebra structures on both sides. In
\cite{ye2} it is shown that \eqref{ref-3.2-11} is a quasi-isomorphism. 

We have isomorphisms
 \begin{equation}
\label{ref-3.3-12}
 D^{\poly,\cdot}(R)[1]\cong D^{\poly,\cdot}(X)[1] \cong R\Gamma(X,\Dscr_X^{\poly,\cdot})^{\tot}[1]
 \end{equation}
The first isomorphism is a tautology and the second one follows from   Proposition \ref{ref-B.8.1-194} below.

Now we drop the restriction that $X$ is affine.   Select an affine open
covering $\Uscr=\{U_1,\ldots,U_n\}$ of $X$ and let the associated notations
be as in the introduction.
\begin{theorem} \label{ref-3.1-13} 
\begin{enumerate} \item There is an
isomorphism in the homotopy category of DG-Lie algebras
\begin{equation}
\label{ref-3.4-14}
R\Gamma(X,\Dscr^{\poly}_X)^{\tot}[1]\cong \bold{C}(\frak{u})[1]
\end{equation}
\item If $X=\Spec R$ is affine and $\Uscr=\{X\}$ then \eqref{ref-3.4-14}
  coincides with the composition of \eqref{ref-3.2-11} and the inverse of
  \eqref{ref-3.3-12}.
\end{enumerate}
\end{theorem}
\begin{proof}
 If $p:\frak{h}\r
\frak{g}$ is a map of cosimplicial DG-Lie algebras then as in
\S\ref{ref-B.7-178} below we say that $p$ is a weak equivalence if $p$
induces a quasi-isomorphism between the cochain complexes
$C^\ast(\frak{h})$ and $C^\ast(\frak{g})$.

We first prove (1).  
Let $\Ch\nolimits^o(\Uscr,\Dscr^{\text{poly}}_X)[1]$ be the ordered {\v{C}ech}
cosimplicial DG-Lie algebra associated to $\Dscr^{\text{poly}}_X[1]$ and
the cover $\Uscr$ (see \S\ref{ref-B.9-195}).  We will construct a cosimplicial DG-Lie algebra
$\Cscr[1]$ together with weak equivalences
\begin{equation}
\label{ref-3.5-15}
\Ch\nolimits^o(\Uscr,\Dscr^{\text{poly}}_X)[1]\r \Cscr[1] \l \bold{C}({\frak{u}})[1]
\end{equation}
where we view $\bold{C}({\frak{u}})[1]$ as a constant cosimplicial
DG-Lie algebra.

Applying the Thom-Sullivan normalization functor $N(-)^{\text{TS}}$ 
(see Appendix \ref{ref-A-137}) we obtain isomorphisms in the homotopy category of DG-Lie algebras
\begin{multline}
\label{ref-3.6-16}
R\Gamma(X,\Dscr^{\poly}_X)^{\tot}[1]\xrightarrow{(1)}N(\Ch\nolimits^o(\Uscr,
\Dscr^{\poly}_X)[1])^{\operatorname{TS}}\r \\
N(\Cscr[1])^{\operatorname{TS}} \l
N(\bold{C}({\frak{u}})[1])^{\operatorname{TS}} \xleftarrow{(2)}
\bold{C}({\frak{u}})[1]
\end{multline}
Here arrow (1) is obtained from \S\ref{ref-B.9-195} below and (2)
is obtained from \eqref{ref-A.1-140}. The composed isomorphism in
\eqref{ref-3.6-16}
yields part (1) of the theorem. Part (2) will follow from the construction of $\Cscr$. 

\medskip 

So now we concentrate on \eqref{ref-3.5-15}. 
 For clarity we will sometimes omit
the shift $[1]$ in the formulas. 
For $\emptyset \neq J\subset I$ let $\frak{u}_J$ be the full subcategory
of $\frak{u}$ spanned by $J'$, $J'\supset J$.  Since $\frak{u}_{J'}\r
\frak{u}_J$ is fully faithful
it follows that there are restriction maps
\begin{equation}
\label{ref-3.7-17}
 \bold{C}(\frak{u}_{J})\r \bold{C}(\frak{u}_{J'})
\end{equation}
for $J'\supset J$. Furthermore it follows from \cite[Lemma
7.5.2]{lowenvdb2} that the 
restriction map
\begin{equation}
\label{ref-3.8-18}
\bold{C}(\frak{u}_{J})\r \bold{C}(\Oscr_{U_J})
\end{equation}
is a quasi-isomorphism. 

Define 
\[
\Cscr^m=\prod_{j_0\le\cdots\le j_m} \bold{C}(\frak{u}_{\{j_0,\ldots,j_m\}})
\]
We make $\Cscr=(\Cscr^m)_m$ into a cosimplicial DG-Lie algebra using
the restriction maps \eqref{ref-3.7-17}.

Consider $\bold{C}(\frak{u})$ as a constant cosimplicial DG-Lie algebra. We
claim that  the restriction map
\[
\bold{C}(\frak{u})\r  \Cscr
\]
is a weak equivalence. To this end it is sufficient to check that we
obtain a quasi-isomorphism on the corresponding (totalized) normalized cochain
complexes. 

The normalized cochain complex of $\Cscr$ it the total complex of a double
complex with columns
\[
\Cscr^{(m)}=\prod_{j_0<\cdots< j_m} \bold{C}(\frak{u}_{\{j_0,\ldots,j_m\}})
\]
Write $\frak{u}_{\emptyset}=\frak{u}$ and $\Cscr^{(-1)}=\bold{C}(\frak{u})$. We
have to show that the total complex associated to the double complex
\begin{equation}
\label{ref-3.9-19}
0\r \Cscr^{(-1)}\r \Cscr^{(0)}\r \cdots \r \Cscr^{(n)}\r 0
\end{equation}
is acyclic. We do this by showing that it is a long exact sequence of
complexes.

$\bold{C}(\frak{u}_J)$ is a direct product of abelian groups of the form
\begin{equation}
\label{ref-3.10-20}
\Hom_k(\Oscr_X(U_{J_1})\times\cdots\times \Oscr_X(U_{J_t}),\Oscr_X(U_{J_t}))
\end{equation}
for $J\subset J_0\subset\cdots\subset J_t$. It follows that the
summands in \eqref{ref-3.9-19} corresponding to a given sequence $J_0,\cdots,J_t$
are parametrized by $J\subset J_0$. Since the signs are the usual alternating ones it
follows easily that for the horizontal differential \eqref{ref-3.9-19} is a sum
of acyclic complexes. This proves what we want.

Our next aim is to construct a cosimplicial map 
\[
c:\Ch\nolimits^o(\Uscr,\Dscr^{\text{poly}}_X)\r \Cscr
\]
We do this by combining maps
\begin{equation}
\label{ref-3.11-21}
c_J:D^{\text{poly}}(U_J)
\r \bold{C}(\frak{u}_J)
\end{equation}
If $d\in D^{\text{poly,t}}(U_J)$ then $d$ is a differential
operator in
\[
\Hom(\Oscr_X(U_{J})^t,\Oscr_X(U_{J}))
\]
It follows that $d$ extends uniquely to a differential operator with
$t$ arguments in \eqref{ref-3.10-20}. We define $c_J(d)$
as this extension.

Now we claim that $c$ is a weak equivalence. To do this it is sufficient to show that
the maps $c_J$ are quasi-isomorphisms. Then it is sufficient that the
composition 
\[
D^{\text{poly}}(U_J)\r \bold{C}(\frak{u}_J)\r 
\bold{C}(\Oscr_{U_J})
\]
of $c_J$ with the quasi-isomorphism \eqref{ref-3.8-18} is a quasi-isomorphism. 

This composition in nothing but \eqref{ref-3.2-11} for $R=\Oscr(U_J)$
and hence we are done.
\end{proof}
\begin{remark} It seems quite likely that the fact that
  $R\Gamma(X,\Dscr_X)^{\tot}$ controls the deformation theory of
  $\Oscr_X$ in the category of twisted presheaves follows also from
  Vladimir Hinich's descent theorem \cite{hinich} given the fact that
  this is true if $X$ is affine. One has to check that the global
  deformation functor is given by gluing the local deformation
  fuctors. Note that for this to work one should view these
  deformation functors as taking values in $2$-groupoids. One should
  also use the fact, first observed by Deligne, that the solutions to
  the Maurer-Cartan equation in a DG-Lie algebra concentrated in
  degrees $\ge -1$ form naturally a $2$-groupoid as well. See
  \cite{Getzler} and the reference to a downloadable copy of Deligne's
  letter contained therein.
\end{remark}
\section{The analytic case}
The main result in this section (Proposition \ref{ref-4.1-24}) was used in the introduction. 
\label{ref-4-22}
Assume now that $k=\CC$ and $X$ is a separated smooth proper scheme over $k$. 
Let $\Oscr^\infty_{X^{\an}}$ be the sheaf of $C^\infty$-functions on
$X^{\an}$ and let  $\Tscr^{i,0}_{X^{\an}}$ and $\Omega^{0,j}_{X^{\an}}$ be the
$C^\infty$-vector bundles respectively generated by the holomorphic vector fields
and by the anti-holomorphic
differential forms. Below we consider the sheaf of DG-Lie algebras
\begin{equation}
\label{ref-4.1-23}
\bigoplus_{i,j} \Tscr^{i,0}_{X^{\an}}\otimes_{\Oscr_{X^{\an}}^\infty} \Omega^{0,j}_{X^{\an}}[1]
\end{equation}
where the differential is obtained by linearly extending
the differential on $\Omega^{0,\cdot}_{X^{\an}}$ and the Lie bracket
is obtained by linearly extending the Lie bracket on
$\Tscr^{\cdot,0}_{X^{\an}}$.
\begin{proposition}
\label{ref-4.1-24}
There is an isomorphism
\[
R\Gamma(X,\Tscr^{\poly,\cdot}_X)^{\tot}[1]\cong \bigoplus_{i,j} \Gamma(X^{\an},\Tscr^{i,0}_{X^{\an}}\otimes_{\Oscr_{X^{\an}}^\infty} \Omega^{0,j}_{X^{\an}})[1]
\]
in the homotopy category of DG-Lie algebras. 
\end{proposition}
\begin{proof}
We first prove that there is an isomorphism 
\begin{equation}
\label{ref-4.2-25}
R\Gamma(X,\Tscr^{\poly,\cdot}_X)^{\tot}[1]\r 
R\Gamma(X^{\an},\Tscr^{\poly,\cdot}_{X^{\an}})^{\tot}[1]
\end{equation}
Choose an affine covering $\Uscr=\{U_1,\ldots,U_n\}$ for $X$. By 
 Lemma \ref{ref-B.9.1-197} and the fact that affine varieties are Stein
we have
\begin{align*}
  R\Gamma(X,\Tscr^{\poly,\cdot}_X)^{\tot}&=N(\Ch\nolimits^o(\Uscr,\Tscr^{\poly,\cdot}_X))^{TS}\\
  R\Gamma(X^{\an},\Tscr^{\poly,\cdot}_{X^{\an}})^{\tot}&=
N(\Ch\nolimits^o(\Uscr,\Tscr^{\poly,\cdot}_{X^{\an}}))^{TS}
\end{align*}
The obvious map
\[
\Ch\nolimits^o(\Uscr,\Tscr^{\poly,\cdot}_X)\r \Ch\nolimits^o(\Uscr,\Tscr^{\poly,\cdot}_{X^{\an}})
\]
yields \eqref{ref-4.2-25}. To prove that \eqref{ref-4.2-25} is an isomorphism it is 
sufficient to prove that it induces an isomorphism on cohomology. The cohomology
on the left and on the right are respectively given by 
\[
H^i(X,\wedge^j \Tscr_X)\qquad\text{and}\qquad H^i(X^{\an},\wedge^j \Tscr_{X^{\an}})
\]
These are are equal because of GAGA. 

We now note that the sheaf of DG-Lie algebras 
\eqref{ref-4.1-23} is an acyclic resolution for $\Tscr^{\poly,\cdot}_{X^{\an}}$ 
(the
$\bar{\partial}$-resolution). 
Using Lemma \ref{ref-B.9.1-197} we have
\begin{equation}
\label{ref-4.3-26}
R\Gamma(X,\Tscr_{X^{\an}}^{\poly,\cdot})^{\tot}[1]\cong 
\bigoplus_{i,j} \Gamma(X^{\an},\Tscr^{i,0}_{X^{\an}}\otimes_{\Oscr_{X^{\an}}^\infty} \Omega^{0,j}_{X^{\an}})[1]
\end{equation}
This concludes the proof.
\end{proof}
\begin{remark} In \cite{BaKo} Barannikov and Kontsevich show that if
  $X$ is Calabi-Yau then 
$\bigoplus_{i,j} \Gamma(X^{\an},\Tscr^{i,0}_{X^{\an}}\otimes_{\Oscr_{X^{\an}}^\infty} \Omega^{0,j}_{X^{\an}})[1]
$
is isomorphic in the homotopy
  category of DG-Lie algebra to the vector space
\[
\bigoplus_{i,j} H^i(X,\wedge^j \Tscr_{X})[1]
\]
with zero differential and Lie bracket. It follows that for $(l,m)$
an Artinian local ring with residue field $k$ the $l$-deformations of
$\Qch(\Oscr_X)$ correspond to elements of
\[
m\otimes_k\biggl(\bigoplus_{i+j=2}
  H^i(X,\wedge^j \Tscr_{X})\biggr)
\]
\end{remark}
\section{Topological notions}
\label{ref-5-27}
Below we will naturally encounter topological rings and modules.
Rather than using Yekutieli's category of
$\operatorname{Dir}\operatorname{Inv}$-abelian groups \cite{ye3} we
 work in the classical setting of filtered topological
abelian groups.
Below we list the few facts we need. 
\subsection{Topological abelian groups}
Below we will encounter exclusively linear topological abelian groups.
 I.e.\ topological abelian groups equipped with a topology such that $0$ has a 
system of  open neighborhoods consisting of subgroups.

The following trivial lemma is often used implicitly.
\begin{lemmas} Assume that $V$ is a linear topological abelian group
and $W$ is an open subgroup. Then $W$ is also closed and the quotient topology
on $V/W$ is discrete.
\end{lemmas}
\begin{proof} $W$ is the complement of the open set $\bigcup_{w\not\in W} (w+W)$
and hence $W$ is closed.
\end{proof}
From now on we assume in addition that there is a countable basis of
neighborhoods of $0$. It is convenient to take for such a basis a
descending chain of subgroups $V=F_0V\supset F_1V\supset\cdots$,
constructed by taking successive intersections in an arbitrary linear
basis indexed by the natural numbers.

We define the completion of a
linear topological abelian group $V$ as
\[
\hat{V}=\invlim_p V/F_pV
\]
It is easy to see that this is the same as the usual completion using Cauchy
sequences.
$\hat{V}$ is a linear topological abelian group with neighborhood basis
of $0$ given by the $(F_pV)\,\hat{}$. $V$ is complete if the map $V\r \hat{V}$
is an homeomorphism.

If $V$, $W$ is are linear topological abelian groups then we make
$V\otimes_{\ZZ} W$ into a linear topological abelian group by
selecting as system of neighborhoods for $0$ the images of the abelian
groups $F_pV\otimes_{\ZZ} W+V\otimes_{\ZZ} F_pW$. The completed tensor
product $V\ctimes_{\ZZ} W$ is the completion of $V\otimes_{\ZZ} W$ for
this topology.

Since now the categories of linear topological abelian groups and complete 
linear topological abelian groups are tensor categories we can define
vector spaces, rings, modules etc\dots\ in them. 

\subsection{Filtered linear topological abelian groups}
\label{ref-5.2-28}
It will be necessary to use filtered linear topological abelian groups. A filtered
linear topological abelian group
is by definition an abelian group $V$, equipped with a filtration $F^\cdot$ such
that each $F^m V$ is equipped with a linear topology with the property
that the inclusion maps $F^m V\r F^{m+1} V$ are continuous.
The filtration is $F^\cdot$ is considered part of the structure of a filtered
linear topological abelian group. 

 We say that
$V$ is complete if each $F^m V$ is complete. The completion of $V$ is
defined as $\bigcup (F^m V )\hat{}$. The category of filtered and
complete filtered linear topological abelian groups have the obvious structure
of a monoidal category, so we can define rings, modules etc\dots in them.
\begin{examples} The ring of differential operators of $k[[t]]$, which is
equal to $k[[t]][\partial_t]$ has
an obvious structure of complete filtered linear topological 
ring.
\end{examples}
\subsection{Finite adic rings}
In this section all rings are commutative.
We will need to complete non-noetherian rings. This
is a somewhat dangerous operation for the following reason. Let $T$
be a ring with an ideal $I$ and put $\hat{T}=\invlim T/I^n$.  If 
we  equip all $T/I^n$ with the discrete topology then $\hat{T}$ 
becomes a topological ring. However in general we will not have
$T/I\cong \hat{T}/I\hat{T}$. The reason for this is that $I\hat{T}$ is
not closed, i.e. $I\hat{T}\neq (I\hat{T})\,\hat{}=\hat{I}$. So we should
replace $I\hat{T}$ by $\hat{I}$. Unfortunately this does not
resolve all our problems since in general  $\hat{I^n}\neq (\hat{I})^n$. 
So we still don't have an isomorphism $R/I^n\cong R/(\hat{I})^n$. The
following example clarifies this.
\begin{examples} Let $T=k[x_1,x_2,\ldots]$ be the polynomial ring
  in infinitely many variables over $k$ and let $I=(x_1,x_2,\ldots)$. Then
  $\hat{T}=T[[x_1,x_2,\ldots]]$.

The ideals $(I^n)\hat{}$ are topologically generated by the monomials
$x_{i_1}\cdots x_{i_n}$ with $i_1\le i_2\le\cdots \le i_n$. 
The following element  shows
\[
x_1+x_2x_3+x_4x_5x_6+\cdots
\]
that $\hat{I}$   is not generated by $x_1,x_2,\ldots$ as an ordinary ideal.

We now show that $\hat{I}^2\neq (I^2)\hat{}$.  If $f\in \hat{I}^2$
then the partial derivatives $f_i=\partial f/\partial x_i$ are in a 
finitely generated ideal in $\hat{I}$ (exercise). Consider the
element 
\[
f=x_1^2+x_2^3+x_3^4+\cdots
\]
By working modulo $x_m$ for $m>n$ we see by looking at heights that an
ideal containing the $f_i$ needs at least $n$
generators. Since $n$ is arbitrary this means that the $f_i$ cannot be 
contained in a finitely generated ideal.
\end{examples}
Luckily all problems go away if we consider completions at finitely
generated ideals.
\begin{definitions} An \emph{adic} ring \cite[0.7.1.9]{EGA0} is a
  linear topological ring such that
\begin{equation}
\label{ref-5.1-29}
T=\invlim_n T/I^n
\end{equation}
where of course $T/I^n$ is equipped with the discrete topology. 
\end{definitions}
An alternative way of stating \eqref{ref-5.1-29} is by saying that $I^n$ is
a fundamental system of \emph{open} neighborhoods of $0$ and that the
topology on $T$ is separated and complete. An ideal $I$ with this
property is called an \emph{ideal of definition} of $T$. Note that if $T$
has the discrete topology then $T$ is adic and the zero ideal is an
ideal of definition.

\medskip

The following definition is non-standard but convenient.
\begin{definitions} A \emph{finite adic ring} is an adic ring with a finitely generated ideal of definition.
\end{definitions}
If $T$ is an adic ring with
ideal of definition $I$ and $M$ is a topological $T$ module then we say
that $M$ is adic if 
\[
M=\invlim_n M/I^nM
\]
where as usual  $M/I^nM$ is equipped with the discrete topology. It
is clear that this definition does not depend on the choice of $I$.
\begin{theorems}  
\label{ref-5.3.4-30}
Assume that $T$ is a
ring with an ideal $I$ such that $I/I^2$ is finitely generated and let $M$ be a $T$-module. 
Let $\hat{T}$ and $\hat{M}$ be respectively the completions of
$T$ and $M$ for the $I$-adic topology. Then we have the following:
\begin{enumerate}
\item $\hat{T}$ is an adic ring with ideal of definition $I\hat{T}$.
\item If $f_1,\ldots,f_d$ are lifts in $I$ of generators of $I/I^2$
then the images of $f_1,\ldots,f_d$ in $\hat{T}$ are generators of $I\hat{T}$.
In particular $\hat{T}$ is finite adic.
\item $\hat{M}$ is adic.
\item We have $(I^nM)\hat{} =I^n \hat{M}$
\item The canonical map $M/I^nM\r \hat{M}/I^n\hat{M}$ is an isomorphism.
\end{enumerate}
\end{theorems}
\begin{proof} This is a slight refinement of \cite[Prop. 0.7.2.7]{EGA0}.
\end{proof}
The category of  adic topological rings has tensor products.
More concretely let $C\r A$, $C\r B$ be continuous maps of linear
topological rings such  that $A,B$ have finitely generated ideals
of definition $I,J$. Then it is easy to see that the topology on
$A\otimes_C B$ induced from $A\otimes B$ has a finitely generated
defining ideal which is given as the image of $K=I\otimes B+A\otimes J$.
We define $A\ctimes_C B$ as the completion of $A\otimes_C B$ for this 
topology. Note that for this definition the topology on $C$ does not have
to be adic. Also note that if $A,B$ are finite adic then so is $A\ctimes_C B$.
\begin{conventions} \label{ref-5.3.5-31} Let $A,B$ be $k$-algebras. It will often happen below
that there is given some finitely finitely generated ideal $I\subset A\otimes B$ and that we are interested in the $I$-adic completion of $A\otimes B$
at $I$. To avoid confusion we will write this completion as $A\cbtimes B$. A we use a similar convention for modules. If $M$, $N$ are
respectively $A$, $B$-modules then $M\cbtimes N$ denotes the $I$-adic
completion of $M\otimes N$ at~$I$.  
\end{conventions}

\medskip

If $T$ is an adic ring then we will write $\Adic(T)$ for the additive
category of adic $T$-modules. Note that any $T$-linear map between
objects in $\Adic(T)$ is automatically continuous. If $T$ has a
finitely generated defining ideal then completion defines a left
adjoint to the forgetful functor
\[
\Adic(T)\r \Mod(T)
\]

\subsection{Differentials}
\label{ref-5.4-32}
If $M$ is a $T$-module then the symmetric group acts on $M^{\otimes
  n}$ (the $n$-fold tensor product of $M$ over $T$) by permuting the factors.

We put
\[
\bigwedge\nolimits^i _TM=\coker \left(\bigoplus_{\sigma\in S_n}M^{\otimes
  n}\xrightarrow{\sum_\sigma\phi_\sigma} M^{\otimes n}\right)
\]
where $\phi_\sigma$ acts by $1-(-1)^{\operatorname{sign}(\sigma)}\sigma$ on the 
summand indexed by $\sigma$. It follows immediately that 
$\bigwedge^i _TM$ is compatible with base change in $T$.

If $T$ is an $R$-algebra then as usual we write 
\[
\Omega^i_{T/R}=\bigwedge\nolimits_T^i \Omega_{T/R}
\]
and we call the collection of $T$-modules $\Omega^i_{T/R}$, together
with its natural differential the \emph{relative De Rham complex} $\Omega^\cdot_{T/R}$
of $T/R$.

If $T$ is an  adic $R$-algebra then we write.
\[
\Omega_{T/R}^{i,\text{cont}}=(\Omega^i_{T/R})\,\hat{}
\]
As derivations are automatically continuous with respect to the $I$-adic
topology (for $I$ an ideal of definition of $T$) the differential on $\Omega^\cdot_{T/R}$ lifts to $\Omega_{T/R}^{\cdot,\text{cont}}$.
We call the resulting complex
the \emph{continuous relative De Rham complex} of $T/R$.

Assume now that $T$ is finite adic. By Theorem \ref{ref-5.3.4-30} the $\Omega_{T/R}^{i,\text{cont}}$ are adic $T$-modules. 

The following formula is convenient
\begin{propositions} 
\label{ref-5.4.1-33} We have
\[
\Omega_{T/R}^{i,\mathrm{cont}}=\invlim_n\Omega^i_{(T/I^n)/R}
\]
\end{propositions}
\begin{proof} We have a standard exact sequence
\[
I^n/I^{2n}\r \Omega_{T/R}/I^n \Omega_{T/R}\r \Omega_{(T/I^n)/R}\r 0
\]
Hence modulo essentially zero systems we have an isomorphism between
the inverse systems $(\Omega_{T/R}/I^n \Omega_{T/R})_n$ and
$(\Omega_{(T/I^n)/R})_n$.  It is easy to see that we obtain from this
an isomorphism between the inverse systems $(\Omega^i_{T/R}/I^n
\Omega^i_{T/R})_n$ and $(\Omega^i_{(T/I^n)/R})_n$, modulo essentially
zero systems. Taking the inverse limit proves what we want.
\end{proof}
\begin{examples} (A.\ Yekutieli) Consider $T=k[[t]]$. As 
$\Omega_{T/k}$ is compatible with localization we have $(\Omega_{T/k})_t=
\Omega_{T_t/k}$. Since $T_t=k((t))$ is a field of infinite transcendence
degree over $k$, it follows that $\Omega_{T/k}$ is a very large object. On the
other hand $\Omega^{\text{cont}}_{T/k}$ is equal to $Tdt$.
\end{examples}
For $M\in \Mod(T)$ let us denote by $\Der^i_R(T,M)$ the set of anti-symmetric
multilinear maps $T\otimes_R \cdots \otimes_R T\r M$ with $i$
arguments, which are derivations in each of their arguments.

Clearly $\Omega_{T/R}^{i}$ represents the functor 
\[
\Der^i_R(T,-):\Mod(T)\r \Mod(R):
\]
Similarly $\Omega_{T/R}^{i,\text{cont}}$ also represents $\Der^i_R(T,-)$ but now
considered as a functor
$\Adic(T)\r \Mod(R)$.

\begin{propositions}
\label{ref-5.4.3-34}  Assume that $T$ is a
  ring with a finitely generated ideal $I$. Let $\widehat{(-)}$ stand for
  $I$-adic completion. Then  the canonical map
\[
\Omega^i_{T/R}/I^n\Omega^i_{T/R}\r
\Omega_{\hat{T}/R}^{i,\text{cont}}/I^n\Omega_{\hat{T}/R}^{i,\text{cont}}
\]
is an isomorphism. In particular
\[
(\Omega^i_{T/R})\hat{}=\Omega_{\hat{T}/R}^{i,\text{cont}}
\]
\end{propositions}
\begin{proof}
$\Omega^i_{T/R}/I^n\Omega^i_{T/R}$ represents the functor
\[
\Der^i_R(T,-):\Mod(T/I^n)\r \Mod(R)
\]
and likewise $\Omega_{\hat{T}/R}^{i,\text{cont}}/
I^n\Omega_{\hat{T}/R}^{i,\text{cont}}$ represents the functor
\[
\Der_R(\hat{T},-):\Mod(\hat{T}/I^n\hat{T})\r \Mod(R)
\]
It is easy to see that these functors are naturally isomorphic
if we make the identification $\Mod(T/I^n)\cong \Mod(\hat{T}/I^n\hat{T})$
(using Theorem \ref{ref-5.3.4-30}).
\end{proof}
\subsection{Formal schemes}
A standard reference for the basic material on formal schemes is
\cite{EGA0}. Since the formal schemes we use are not noetherian, 
we recall the basics.  Let $I$ be an ideal in a ring $T$. For $f\in
T/I$ let $\tilde{f}\in T$ stand for an arbitrary lift of $f$. Then it
is easy to see that the completion of $T_{\tilde{f}}$ with respect to
the ideal $I_{\tilde{f}}$ is (canonically) independent of the choice
of $\tilde{f}$. We will write $\hat{T}_f$ for this completed
localization. If $M$ is a $T$-module then $\hat{M}_f$ can be define
likewise.

If $T$ is finite adic then follows from Theorem
\ref{ref-5.3.4-30} that all $\hat{T}_f$ are finite adic as well and furthermore
all $\hat{M}_f$ are adic $\hat{T}_f$-modules.

An alternative
definition for $\hat{M}_f$ is the following.  $f$ defines an open subset $D(f)$
of $\Spec T/I^n=\Spec T/I$. Then $\hat{M}_f$ is given by the global sections of
\[
\invlim_n (M/I^n M)\,\tilde{} \mid D(f)
\]
(where as usual $(-)\,\tilde{}$ denotes the quasi-coherent sheaf associated to
a module).

Now we sheaffify these constructions. Recall the technically useful fact that 
an inverse limit of sheaves can be computed as presheaves. 

If $T$ is a finite adic ring with finitely generated ideal of definition we 
define 
\[
\Spc T=(\Spec T/I, \invlim_n (T/I^n)\,\tilde{}\,)
\]
for an ideal of definition $I$. So $\Spc T$ is a topologically ringed $(\Tscr,\Oscr_\Tscr)$  with  $\Tscr=\Spec T/I$.  It is clear that
the definition of $\Spc T$ is independent of $I$ and
\[
(D(f),\Oscr_{\Tscr}{\mid} D(f) )=\Spc \hat{T}_f
\]
We use a special type of formal scheme. For the full definition see
\cite{EGA0}.
\begin{definitions} An \emph{finite adic affine formal scheme} is a
  topologically ringed space which is isomorphic to $\Spc T$ for a
  finite adic ring $T$.
\end{definitions}
Let $T$ now be a finite adic ring with finitely generated ideal of
definition $I$.  Put $(\Tscr,\Oscr_\Tscr)=\Spc T$.  If $M$ is an adic
$T$-module then we define
\[
M^{\triangle}=\invlim_n (M/I^n M)\tilde{}
\]
Thus $M^{\triangle}$ is a sheaf of topological $\Oscr_\Tscr$-modules such that 
\[
M^{\triangle}\mid D(f)=(\hat{M}_f)^{\triangle}
\]
$\Oscr_\Tscr$ itself contains a sheaf of ideals $\Iscr=I^{\triangle}$ such that
\[
\invlim_n\Oscr_\Tscr/\Iscr^n=\Oscr_\Tscr
\]
We call $\Iscr$ an ideal of definition of $\Oscr_\Tscr$. Being an ideal of
definition is a local property (see \cite[Prop 10.3.5]{EGA0}).

\medskip 

Chaining the various definitions and using Proposition \ref{ref-5.4.3-34} we also find
\[
(\Omega_{T/R}^{i,\text{cont}})^{\triangle}\mid D(f)=(\Omega_{\hat{T}_f/R}^{i,\text{cont}})^\triangle
\]
\begin{definitions} Let $(\Xscr,\Oscr_\Xscr)$ be a topologically ringed space.
  We say that $\Xscr$ is a finite adic formal scheme
  if $\Xscr$ is locally a finite adic affine topologically ringed space.
\end{definitions}
As usual we call $\Oscr_\Xscr$ the structure sheaf of $\Xscr$.  We say
that a topological $\Oscr_\Xscr$-module is adic if it is locally of
the form $M^{\Delta}$. If $X\r Y$ is a map of a finite adic formal
scheme to a scheme then we define $\Omega_{Y/X}^{i,\text{cont}}$ as
the $\Oscr_\Xscr$-module is which is locally, on open affine formal
subschemes $\Spc(T)$, of the form $\Omega_{T/R}^{i,\text{cont}}$.

Morphisms between formal schemes are by definition morphisms of locally
topologically ringed spaces \cite[Def 10.4.5]{EGA0}. If $S,T$ are finite adic rings
then 
\[
\Hom_{\text{formal schemes}}(\Spc S,\Spc T)=\Hom_{\text{topological
    rings}}(T,S)
\]
\cite[Prop 10.2.2]{EGA0}.

Clearly the category of (affine) schemes is a full subcategory of the
category of finite adic (affine) formal schemes.

The category of finite adic formal schemes has fibered products \cite[Prop
10.7.3]{EGA0}. As usual it is sufficient to construct these for affine
formal schemes. If  $C\r A$, $C\r B$  are continuous maps of
adic rings then we put
\[
\Spc A\times_{\Spc C}\Spc B=\Spc (A\hat{\otimes}_C B)
\]
Suppose $X$ is a scheme and $Y$ is a closed subscheme defined by a
quasi-coherent ideal $\Iscr$ which is (locally) of finite type. Then
$\hat{X}_Y$ (or simply $\hat{X}$) is the finite adic formal scheme whose
underlying space is $Y$ and whose structure sheaf is
\[
\Oscr_{\hat{X}_Y}=\hat{\Oscr}_{X,Y}=\invlim_n \Oscr_X/\Iscr^{n}
\]
\section{Formal geometry}
\label{ref-6-35}
Many ideas in this section are taken from \cite[\S 5]{ye3}. However we put
more emphasis on the language of  formal schemes. 
\subsection{Basic definitions}
Everything will be over an algebraically closed base field $k$ of
characteristic zero.  Fix an integer $d$. For a finite adic scheme $Y$
with locally finitely generated ideal of definition $\Iscr$ we let
$Y[[t_1,\ldots,t_n]]$ be the finite adic formal scheme which is the
completion of $Y\times \AA^d_k$ at the ideal $\Iscr+(t_1,\ldots,t_d)$.
The inclusion/projection $Y=Y\times \{0\}\hookrightarrow Y\times
A^d_k\r Y$ yields a canonical projection map
\begin{equation}
\label{ref-6.1-36}
p_Y:Y[[t_1,\ldots,t_n]]\r Y
\end{equation}
with section
\begin{equation}
\label{ref-6.2-37}
i_Y:Y\r Y[[t_1,\ldots,t_n]]
\end{equation}
If $Y=\Spc S$ is affine then 
\[
Y[[t_1,\ldots,t_d]]\cong \Spc S[[t_1,\ldots,t_d]]
\]
We like to view \eqref{ref-6.1-36} as an infinite dimensional vector bundle
over $Y$ with zero-section given by \eqref{ref-6.2-37}.

Let $\operatorname{FSch}/k$ be the category of finite adic formal schemes over
$k$.
\begin{propositions} \label{ref-6.1.1-38} The functor 
\[
\operatorname{FSch}/k\r \operatorname{FSch}/k: 
Y\mapsto Y[[t_1,\ldots,t_d]]
\]
has a right adjoint.
\end{propositions}
\begin{proof}  Let $X\in \operatorname{FSch}/k$. We need to show
that the contravariant functor 
\begin{equation}
\label{ref-6.3-39}
\Phi:\operatorname{FSch}/k\r \Set: Y\mapsto \Hom_{\operatorname{FSch}/k}(Y[[t_1,\ldots,t_d]],X)
\end{equation}
is representable.

Since maps of finite adic formal schemes are compatible with gluing
we reduce to $Y=\Spc S$, 
$X=\Spc R$. So we may work with the category of finite adic affine formal schemes,
or equivalently, the category of finite adic rings. Thus 
\[
\Phi(S)=\Hom(R,S[[t_1,\ldots,t_d]])
\]
Let $R^{\mathbf d,\circ}$ be the $k$-algebra generated by symbols $f_{\underline i}$ for $f\in R$ and
$\underline i=(i_1,\ldots,i_d)$ with relations
\begin{align*}
(f+g)\,\tilde{}
&=\tilde{f}+\tilde{g}\\
(fg)\,\tilde{}&=\tilde{f}\tilde{g}\\
\tilde{\lambda}&=\lambda\qquad \text{(for $\lambda\in k$)}
\end{align*}
where $\tilde{f}$ is the generating function $ \sum_{\underline i}
f_{\underline i} t^{\underline i} $ with $t^{\underline
  i}=t_1^{i_1}\cdots t_d^{i_d}$.  In particular there is a ring homomorphism
\[
R\r R^{\bold{d},\circ}:f\mapsto f_{(0,\ldots,0)}
\]
 Let $I$ be a finitely generated
defining ideal for $R$ and consider the ideal $J\subset R^{\mathbf
  d,\circ}$ the ideal generated by $f_{(0,\cdots,0)}$ for $f\in I$. 
Then clearly $J$ is also finitely generated. Let $R^{\mathbf d}$ be the completion
of $R^{\mathbf
  d,\circ}$ at $J$.

It is easy to see that we have a functorial isomorphism for any
finite adic $k$-algebra $S$:
\begin{equation}
\label{ref-6.4-40}
\mu:\Hom(R,S[[t_1,\ldots,t_d]])\r \Hom(R^{\textbf{d}},S)
\end{equation}
where $\mu$ is defined by 
\[
\sum_{\underline{i}}\mu(\phi)(f_{\underline{i}})\,t^{\underline{i}}=\phi(f)
\]
Hence
$R^{\mathbf d}$ represents $\Phi$.
\end{proof}
Below we denote the right adjoint to $Y\mapsto Y[[t_1,\ldots,t_d]]$ by 
$X\mapsto X^{\mathbf d}$. The proof of the previous proposition shows
that if $X$ is affine then so is  $X^{\mathbf d}$.

\medskip

Now assume that $\phi:Y\r X$ is a map between  $k$-schemes
such that $X$ is  separated and of finite type.  Then the graph
$\Gamma_\phi=(\phi,\Id_Y):Y\r X\times Y$ of $\phi$ is closed and its
defining ideal is of finite type. We let $\Jet_{\phi,\infty}$ be the
completion of $Y\times X$ along $\Gamma_\phi$. Thus $\Jet_{\phi,\infty}$
comes equipped with a map of formal schemes
\[
\Jet_{\phi,\infty}\r X\times Y
\]
An interesting special case is when $\phi$ is the inclusion of a closed point
$x$ in $X$. In that case
\begin{equation}
\label{ref-6.5-41}
\Jet_{\phi,\infty}=\Spc \hat{\Oscr}_{X,x}
\end{equation}

It is easy to see that
$\Jet_{\phi,\infty}$ is compatible with base extension in the sense that
if there is a commutative diagram
\[
\begin{CD}
Z @>>> Y\\
@V\theta VV @VV\phi V\\
X @= X
\end{CD}
\]
then 
\begin{equation}
\label{ref-6.6-42}
\Jet_{\theta,\infty}=Z\times_Y \Jet_{\phi,\infty}
\end{equation}
We will write 
\[
\Jet_{X,\infty}=\Jet_{\Id_X,\infty}
\]
We obtain in particular
\[
Y\times_X \Jet_{X,\infty}=\Jet_{\phi,\infty}
\]

Now let $Y$ be a an arbitrary $k$-scheme. Fix a map
\[
\phi:Y[[t_1,\ldots,t_d]]\r X
\]
We then get a commutative diagram
\begin{equation}
\label{ref-6.7-43}
\begin{CD}
Y @= Y \\
@Vi_Y VV @VV (\phi\circ i_Y,\Id_Y) V\\
Y[[t_1,\ldots,t_d]] @>>(\phi,p_Y)  > X\times Y\\
@V p_Y VV @VV\operatorname{pr}_2V\\
Y @= Y 
\end{CD}
\end{equation}
where the composition of the vertical arrows is the identity. Put
$\phi_0=\phi\circ i_Y$.

If we complete the middle arrow at $Y$ we get a map.
\begin{equation}
\label{ref-6.8-44}
\hat{\phi}_0:Y[[t_1,\ldots,t_d]]\r \Jet_{\phi_{0,\infty}}
\end{equation}
\begin{definitions}
\label{ref-6.1.2-45} 
Let the notations be as above. We say that $\phi$ is a local coordinate
system parametrized by $Y$ if $\hat{\phi}_0$ is an isomorphism of formal 
schemes.
\end{definitions}
Assume that $\phi$ is a local system of parameters and fix a $k$-point
$y:\Spec k\r Y$. Pulling back \eqref{ref-6.8-44} in the category of formal
schemes and using \eqref{ref-6.5-41}\eqref{ref-6.6-42} we get an isomorphism of $k[[t_1,\ldots,t_n]]$ and $\hat{\Oscr}_{X,x}$
where $x=(\phi\circ i_Y)(y)$. In particular $X$ is smooth at the image of $x$.

It follows that it is meaningless to talk about local coordinate
systems for non-smooth schemes. So we now assume that $X$ is smooth of
dimension $d$ over $k$.

Assume in addition that $X$ has a system of parameters
$x_1,\ldots,x_d$.  Let $\Iscr\subset \Oscr_{X\times X} $ be the
defining ideal of the diagonal $\Delta$.  Then the sections
$x'_i=x_i\otimes 1-1\otimes x_i$ of $\Iscr$ form a generating regular
sequence of the ideal of definition $\hat{\Iscr}$ of $\Jet_{X,\infty}$
(we are in the context of noetherian schemes and noetherian formal
schemes so there are no subtleties).  Invoking Proposition \ref{ref-6.1.3-47}
below (for affine open subsets of $\Jet_{X,\infty}$) we find:
\[
\Jet_{X,\infty}\cong \Delta[[x'_1,\ldots, x'_d]]
\]
and hence by base extension
\begin{equation}
\label{ref-6.9-46}
\Jet_{\phi,\infty}\cong \Gamma_\phi[[x'_1,\ldots, x'_d]]
\end{equation}
\begin{propositions}
\label{ref-6.1.3-47}
Assume that we have a  map $R\r T$ where $R$ is a ring and $T$ is finite adic. Assume that $T$ has a defining ideal $I$
generated by a regular sequence $x_1,\ldots,x_n$ such that $T/I\cong R$ and
such that the composition $R\r T\r R$ is an isomorphism. Then
$T\cong R[[x_1,\ldots,x_n]]$.
\end{propositions}
\begin{proof} Putting $\phi(x_i)=x_i$ defines a continuous map 
\[
\phi:R[[x_1,\ldots,x_n]]\r T
\]
Since in both rings the $(x_i)_i$ form a regular sequence this map becomes
an isomorphism after taking associated graded rings. And since the topologies
involved are separated and complete, this easily implies that $\phi$
is an isomorphism.
\end{proof}
\begin{theorems} \label{ref-6.1.4-48} The subfunctor
\[
\Phi^0:\Sch/k \r \Set:Y\mapsto \{\text{local coordinate systems on $X$}\}
\]
of $\Phi$ (as in \eqref{ref-6.3-39}) is representable by an open subscheme $X^{\coord}$ of $X^{\mathbf d}$ which 
is still affine over $X$.
\end{theorems}
\begin{proof} 
As usual this is a local statement on $X$. Hence we may assume that $X$
has a system of parameters $x_1,\ldots,x_d$.

Assume given a local coordinate system on $X$, indexed by $Y$:
\[
\phi:Y[[t_1,\ldots,t_d]]\r X
\]
Thus by the above discussion we obtain a map
\begin{equation}
\label{ref-6.10-49}
\hat{\phi}_0:Y[[t_1,\ldots,t_d]]\r \Gamma_\phi[[x'_1,\ldots,x'_d]]
\end{equation}
We may write the pullbacks of the $x'_i$ as $\sum_j a_{ij} t_j+\cdots$ for functions
$a_{ij}$  on $Y$. Put $\det \phi_0=\det a_{ij}$. Then $\hat{\phi_0}$ is
an isomorphism if and only if $\det \phi_0$ is a unit. Note that $\det \phi_0$
is well defined up to a unit.  

Since it is easy to see that the formation of \eqref{ref-6.10-49} is
compatible with pullbacks this implies that $\Phi^0$ is an open
subfunctor of $\Phi$. So it is representable by an open subscheme
$X^{\coord}$ of $X^{\mathbf d}$. It is in fact represented by the open
subset defined by $\det \theta_0$ for $\theta:X^{\textbf
  d}[[t_1,\ldots,t_d]]\r X$ the universal map.  From the fact that
$X^{\textbf d}$ is affine over $X$, we deduce that $X^{\coord}$ is
affine as well.
\end{proof}
Let $\theta:X^{{\coord}}[[t_1,\ldots,t_d]]\r X$ be the universal local
coordinate system on $X$. From Definition \ref{ref-6.1.2-45} we
obtain a canonical isomorphism
\begin{equation}
\label{ref-6.11-50}
\hat{\theta}_0:X^{{\coord}}[[t_1,\ldots,t_d]]\r \Jet_{\theta_0,\infty} 
\end{equation}
of finite adic $X^{{\coord}}$-schemes.

If $X=\Spec R$ for a $d$-dimensional smooth $k$-algebra  we write $R^{\coord}$
for the coordinate ring of $X^{\coord}$.
We obtain an
isomorphism
\begin{equation}
\label{ref-6.12-51}
R^{\coord}\cbtimes
R \r R^{\coord}[[t_1,\ldots,t_d]]:r\otimes f
\mapsto r\tilde{f}
\end{equation}
where following Convention \ref{ref-5.3.5-31} we let $R^{\coord}\cbtimes
R$ be the completion of $R^{\coord}\otimes R$ at the kernel of the
multiplication map
\begin{equation}
\label{ref-6.13-52}
R^{\coord} \otimes R\r R^{\coord}:r\otimes f\mapsto r \tilde{f}
\end{equation}
\begin{examples} 
  \label{ref-6.1.5-53} It is instructive to understand the isomorphism
  \eqref{ref-6.12-51} in the simplest possible case, namely when
  $R=k[x]$. In that case
\[
R^{\coord}=k[x_0,x_1,\dots]_{x_1}
\]
and \eqref{ref-6.12-51} is given by
\begin{equation}
\label{ref-6.14-54}
(k[x,x_0,x_1,\dots]_{x_1})\,\hat{}\r
(k[t,x_0,x_1,\dots]_{x_1})\,\hat{}: x_i\mapsto x_i,x\mapsto \sum_{i\ge 0} x_i t^i
\end{equation}
where the first completion is at the ideal $(x-x_0)$ and the second completion
is at the ideal $(t)$. To see directly that \eqref{ref-6.14-54} is
an isomorphism we look at associated graded rings. Putting $\delta=(x-x_0)$
the associated graded map to \eqref{ref-6.14-54} is given by
\[
(k[x_0,x_1,\dots]_{x_1})[\delta]\r (k[x_0,x_1,\dots]_{x_1})[t]:
x_i\mapsto x_i, \delta\mapsto x_1 t
\]  
This is clearly an isomorphism.
\end{examples}
\subsection{Groups and  actions}
\label{ref-6.2-55}
 Recall 
that by definition for $R$ a smooth $d$-dimensional $k$-algebra we have
\[
\Hom(R^{\mathbf d},S)\cong \Hom(R,S[[t_1,\ldots,t_d]])
\]
where $S$ is an arbitrary finite adic $k$-algebra. According to
Definition \ref{ref-6.1.2-45} a map $\phi:R^{\mathbf d}\r S$ represents a local
coordinate system on $\Spec R$, parametrized by $\Spc S$ (i.e.\ an
element of $\Hom(R^{{\coord}},S)$) if the
corresponding map $\phi:R\r S[[t_1,\ldots,t_d]]$ induces an
isomorphism
\begin{equation}
\label{ref-6.15-56}
S\cbtimes R\r S[[t_1,\ldots,t_d]]
\end{equation}
where $S\cbtimes R$ is the completion of $S\otimes R$ at the kernel of
the ideal $S\otimes R\r S:s\otimes r\mapsto s\phi_0(r)$ where $\phi_0$
is the kernel of the composition 
\[
R\xrightarrow{\phi} S[[t_1,\ldots,t_d]]\xrightarrow{t_i\mapsto 0} S
\]
\begin{lemmas}
The functor  which sends a finite adic $k$-algebra $S$
to the group
\begin{equation}
\label{ref-6.16-57}
\Aut_S(S[[t_1,\cdots,t_d]])
\end{equation}
is representable by an affine finite adic formal $k$-scheme. 
\end{lemmas}
\begin{proof}
  We sketch the proof which is similar to the proof of Proposition
  \ref{ref-6.1.1-38}. Let $A$ be the $k$-algebra generated by
  variables $z_{i,j_1,\ldots,j_d}$ for $i=1,\ldots, d$, $j_l\ge 0$,
  localized at the determinant of the matrix $z_{i,e_j}$ where
  $e_{j}=(0,\ldots,1,\ldots,0)$ has its $1$ in the $j$'th position.

Let $\hat{A}$ be the completion of $A$ at the ideal
generated by $(z_{i,0,\ldots,0})_i$. There is a bijection
\[
\mu:\Hom(\hat{A},S)\r \Aut_S(S[[t_1,\ldots,t_d]])
\]
defined by
\[
\mu(\phi)(t_i)=\sum_{j_1,\ldots,j_d} \phi(z_{i,j_1,\ldots,j_d}) t_1^{j_1}\cdots t_d^{j_d}
\]
 Thus $\hat{A}$ represents
\eqref{ref-6.16-57} (as finite adic ring). Hence $\Spc \hat{A}$ represents
\eqref{ref-6.16-57} in the category of finite adic formal schemes.
\end{proof}
Below we denote the representing object of \eqref{ref-6.16-57} by $G$. It is 
a group object in the category of formal schemes. 

\medskip

The canonical action of $G(S)$ on $\Hom(R,S[[t_1,\ldots,t_d]])$ now
defines an action on $\Hom(R^{\textbf{d}},S)$.  Since this action is
functorial in $S$ we obtain an action of $G$ on $R^{\textbf{d}}$.
Since $G(S)$ also acts on local coordinate systems it is clear that we
obtain in addition an action of $G$ on $R^{{\coord}}$.  Finally since
everything is compatible with base change these actions globalize to
the case of not necessarily affine $d$-dimensional smooth $k$-schemes.
\begin{propositions} 
\label{ref-6.2.2-58}
Let $X$ be a separated $d$-dimensional smooth $k$-scheme.
  Then the action of $G$ on $X^{\coord}$ is free in the sense that
 \begin{equation}
\label{ref-6.17-59}
G\times X^{\coord}\r X^{\coord}\times X^{\coord}:(g,x)\r (x,gx)
\end{equation}
is a monomorphism.
\end{propositions}
\begin{proof} For $S$ finite adic we have to show that \eqref{ref-6.17-59} induces
an injection
\[
G(S)\times X^{\coord}(S)\r X^{\coord}(S)\times X^{\coord}(S)
\]
As usual we may reduce to the case that $X=\Spec R$ is affine. Then the
statement amounts to proving that $G(S)$ acts freely on morphisms
\[
R\r S[[t_1,\ldots,t_d]]
\]
defining a local coordinate systems parametrized by $\Spc S$. This follows from the existence of 
the isomorphism \eqref{ref-6.15-56}.
\end{proof}
The Lie algebra of $G$ is defined as the kernel
$G(k[\epsilon])\r G(k)$. It can be naturally identified with the Lie
algebra $\frak{g}$ of derivations of $k[[t_1,\cdots,t_d]]$.

The following remarkable result is the main result of ``formal
geometry'' \cite{GeKa}. It says that in a suitable sense $X^{\coord}$ is
a principal homogeneous space over $G$.
We will not explicitly use it however.
\begin{propositions} As before let $X$
  be a separated smooth $k$-scheme of dimension $d$. For $x\in X^{\coord}$ let
  $T_x(X^{\coord})$ be the tangent space at $x$, i.e. the set of maps
  $\Spec k[\epsilon]/(\epsilon^2)\r X^{\coord}$ such that the
  composition $\Spec k\r \Spec k[\epsilon]/(\epsilon^2)\r X^{\coord}$
  is $x$. Then the map $\frak{g}\r T_x(X^{\coord})$  induced by the
  $G$-action on $X^{\coord}$ is an isomorphism of vector spaces.
\end{propositions} 
\begin{proof} Since \eqref{ref-6.17-59} is a monomorphism we obtain an injection
\[
\frak{g}\times T_x(X^{\coord})\r T_x(X^{\coord})\times T_x(X^{\coord})
\]
we have to prove that this is a bijection. That is, if $x_1$, $x_2$ are
$k[\epsilon]/(\epsilon^2)$-points of $X^{\coord}$ mapping $x$ then there
is an $g\in G(k[\epsilon]/(\epsilon^2))$,  mapping to the identity in $G(k)$ such
that $gx_1=x_2$. 

We may assume that $X=\Spec R$ is affine.  Let $x^\circ$ be the image
of $x$ in $X$.  Then $x$ is given by a map
\[
x:R\r k[[t_1,\ldots,t_d]]
\]
inducing an isomorphism
\[
\hat{R}\r k[[t_1,\ldots,t_d]]
\]
where $\hat{R}$ is the completion of $R$ at the maximal ideal of $R$
defining $x^\circ$.

Then $x_1,x_2$ are maps maps making the following diagram commutative
\begin{equation}
\label{ref-6.18-60}
\xymatrix{ 
R \ar[r]^{x_1} \ar[dr]^x\ar[d]_{x_2} &
k[\epsilon]/(\epsilon^2)[[t_1,\ldots,t_d]]\ar[d]\\
k[\epsilon]/(\epsilon^2)[[t_1,\ldots,t_d]]\ar[r]& k[[t_1,\ldots,t_d]] 
}
\end{equation}
Both $x_1$, $x_2$ induce  $k[\epsilon]/(\epsilon^2)$-linear
isomorphisms
\[
\hat{R}\otimes k[\epsilon]/(\epsilon^2)\r
k[\epsilon]/(\epsilon^2)[[t_1,\ldots,t_d]]
\]
From this it follows we can complete the diagram
\eqref{ref-6.18-60} with a $k[\epsilon]/(\epsilon^2)$-linear
diagonal arrow $k[\epsilon]/(\epsilon^2)[[t_1,\ldots,t_d]]\r
k[\epsilon]/(\epsilon^2)[[t_1,\ldots,t_d]]$. This is the required
element of $G(k[\epsilon]/(\epsilon^2))$.
\end{proof}
The action of $G$ on $R^{\bold{d}}$ may be differentiated to an action of
$\frak{g}$ on $R^{\bold{d}}$.
The following proposition describes the nature of this action.
\begin{propositions} 
\label{ref-6.2.4-61} Let $R$ be a smooth affine $k$-algebra of
  dimension $d$. For $f\in R$ let $\tilde{f}\in R^{\mathbf
    d}[[t_1,\ldots,t_d]]$ be as in the proof of Proposition
  \ref{ref-6.1.1-38}. For $v\in \frak{g}$ let $L_v$ be the action $v$ on
  $R^{\mathbf d}[[t_1,\ldots,t_d]]$ obtained by linearly extending the action of $v$ on
  $k[[t_1,\ldots,t_d]]$ (recall that $\frak{g}=\Der_k(k[[t_1,\ldots,t_d]])$). Let $L_{\bar{v}}$ be the action of $v$ on
  $R^{\mathbf d}[[t_1,\ldots,t_d]]$ by linearly extending the action of $v$ on
  $R^{\mathbf d}$.  Then we have
\begin{equation}
\label{ref-6.19-62}
L_{\bar{v}} (\tilde{f})=-L_v(\tilde{f})
\end{equation}
\end{propositions}
\begin{proof} 
 We have the universal map
\[
R\r R^{\textbf{d}}[[t_1,\ldots,t_d]]:f\mapsto \tilde{f}
\]
which which is easily seen to be $G$-invariant (for $G$-acting
trivially on $R$).  Formula
\eqref{ref-6.19-62} expresses the fact that the differentiated $G$-action,
given by $L_{\bar{v}}+L_v$ acts trivially on $\tilde{f}$ for $f\in R$.
\end{proof}
\begin{examples} It is again interesting to consider the simple case $R=k[x]$.
We have 
\[
R^{\bold{d}}=k[x_0,x_1,\ldots]
\]
and 
\[
\frak{g}=k[[t]]\partial_t
\]
To compute the action of $\frak{g}$ we note that $\frak{g}$ has 
a $k$-linear topological basis given by $\delta_i=t^i\partial_t$. We have
\[
[\delta_i,\delta_j]=(j-i)\delta_{i+j-1}
\]
To compute the action of $\delta_i$ we use the method of proof of
Proposition \ref{ref-6.2.4-61}. Thus we use
\begin{align*}
  0&=\delta_i(\sum_j x_j t^j)\\
  &=\sum_j \delta_i(x_j)t^j+\sum_j x_j j t^{i+j-1}\\
  &=\sum_j \delta_i(x_j)t^j+\sum_{j\ge i-1}(j-i+1) x_{j-i+1} t^j\\
  &=\sum_{j<i-1} \delta_i(x_j)t^j+ \sum_{j\ge i-1} (\delta_i(x_j)+
  (j-i+1) x_{j-i+1}) t^j
\end{align*}
Thus it follows
\[
\delta_i(x_j)=
\begin{cases}
0&\text{if $j<i-1$}\\
-(j-i+1) x_{j-i+1}&\text{if $j\ge i-1$}
\end{cases}
\]
or  simply
\begin{equation}
\label{ref-6.20-63}
\delta_i(x_j)=-(j-i+1) x_{j-i+1}
\end{equation}
using the convention $x_j=0$ for $j<0$. 

We obtain for $\alpha_i\in k$
\[
\biggl(\sum_{i \ge 0} \alpha_i \delta_i\biggr) \cdot x_j=
-\sum_{i\ge 0} \alpha_i (j-i+1) x_{j-i+1}
\]
This is a finite sum so the action of $\frak{g}$ on 
$R^{\bold{d}}$ is indeed well defined. 

\medskip

It is clear from \eqref{ref-6.20-63} that the action of $\delta_0$ of
$R^{\bold{d}}$ cannot be exponentiated (i.e.\ the action of
$e^{\delta_0}$ on $x_i$ does not yield a finite sum). So the
$\frak{g}$-action cannot be exponentiated to a group action on
$R^{\bold{d}}$. \emph{However} the results in this section show that
 $\frak{g}$ can be exponentiated to a group $G$ \emph{in the category
  of finite adic schemes}.
\end{examples}
\subsection{Affine coordinate systems}
By restricting ourselves to  linear coordinate changes we may view $\Gl_d$ as
a subgroup of $G$. The action \eqref{ref-6.17-59} now restricts to 
a free action
\[
\Gl_d\times X^{{\coord}}\r X^{{\coord}}
\]
Since $X^{\coord}/X$ is affine we may define the scheme
$X^{{\aff}}=X^{\coord}/\Gl_d$.  Following our usual practice we write
$R^{\aff}$ for the coordinate ring of $(\Spec R)^{\aff}$.

\medskip

The advantage of $X^{{\aff}}$ over $X^{\coord}$ is the following property.
\begin{propositions} 
\label{ref-6.3.1-64}
$X^{{\aff}}$ is a bundle of ($\infty$-dimensional)
affine spaces over $X$.
\end{propositions}
\begin{proof} Assume that $X=\Spec R$ is affine and that $R$ has a
  system of parameters $x_1,\ldots,x_d$. Consider the closed subscheme
  $Y$ of $X^{\coord}$ whose $S$-points are given by maps
\[
\phi:R\r S[[t_1,\ldots,t_d]]
\]
such that $\phi(x_i)=a_i+t_i+\cdots$ for certain $a_i\in S$. It is clear
that the obvious map $\Gl_d\times Y\r X$ defines a bijection on $S$-points
and hence is an isomorphism. Thus $Y\cong X/\Gl_d$.

Using the fact that $R/k[x_1,\ldots,x_d]$ is etale (and hence formally
etale) we see that any diagram
\[
\begin{CD}
R @>>> S\\
@AAA @AA t_i\r 0 A\\
k[x_1,\ldots ,x_d] @>>> S[[t_1,\ldots,t_d]]
\end{CD}
\]
may be completed uniquely with a diagonal arrow $R\r S[[t_1,\ldots,t_d]]$.

It is now clear that sending $\phi$ to $\phi_0$ together with the coefficients
of $\phi(x_i)$ of the terms of degree $\ge 2$ defines a bijection between
$Y$ and the $S$-points of the product of $\Spec R$ with an infinite dimensional affine space. 
This proves what we want.
\end{proof}
\subsection{The abstract formalism of Maurer Cartan forms}
\label{ref-6.4-65}
This is an abstract section whose results will be employed in the next
section. 
We consider
\[
\frak{g}=\Der_{k[[t_1,\ldots,x_d]]}(k[[t_1,\ldots,t_d]])
\]
together with its natural  topological Lie algebra structure. Let $T$ be a finite adic $k$-algebra. Then the DG-Lie
algebra 
\[
\Omega^\cdot_{T/k}\ctimes_k \frak{g}
\]
may be written as 
\[
\sum_i \Omega^\cdot_{T/k}[[t_1,\cdots,t_d]]\left[\frac{\partial\ }{\partial t_i}\right]
\]
and hence its action on
\[
\Omega^\cdot_{T/k}[[t_1,\ldots,t_d]]
\]
is clearly faithful.

  We want to classify
the derivations of degree one on
$\Omega^\cdot_{T/k}[[t_1,\ldots,t_d]]$ such that the natural map of
algebras
\[
\Omega^\cdot_{T/k}\r \Omega^\cdot_{T/k}[[t_1,\ldots,t_d]]
\]
becomes a map of DG-algebras. Below $d$ is such a differential.

Clearly $d$ is determined by the values
 $\omega_i=dt_i\in \Omega^1_{T/k}[[t_1,\ldots,t_d]]$ or
equivalently by the restriction of $d$ to $k[[t_1,\dots,t_d]]$. This
yields us a derivation
\[
\delta:k[[t_1,\dots,t_d]]\r \Omega^1_{T/k}[[t_1,\ldots,t_d]]
\]
such that $\delta(t_i)=\omega_i$. 

Put
\[
\omega=\sum_i\omega_i \frac{\partial }{\partial t_i}\in \Omega^1_{T/R}\ctimes_k \frak{g}
\]
Then with a slight abuse of notation $d$ may be written as
\[
d=d_0+\omega
\]
where $d_0$ is the extension of the differential on
$\Omega^\cdot_{T/R}$. The fact that $d^2=0$ translates into the
identity
\begin{equation}
\label{ref-6.21-66}
d_0\circ \omega+\omega\circ  d_0+\omega\circ \omega=0
\end{equation}
as operations on $\Omega^{\cdot}_{T/R}[[t_1,\ldots,t_d]]$.
The left hand side of this identity is the image of
\[
d_0\omega+\frac{1}{2}[\omega,\omega]
\]
in $\Omega^\cdot_{T/k}\ctimes \frak{g}$.  Hence (using faithfulness) \eqref{ref-6.21-66} is nothing but the 
Maurer-Cartan equation
\[
d_0\omega+\frac{1}{2}[\omega,\omega]=0
\]
in the DG-Lie algebra $\Omega^\cdot_{T/k}\ctimes \frak{g}$ (compare with \eqref{ref-1.4-4}).
\subsection{The Maurer-Cartan form on coordinate spaces}
\label{ref-6.5-67}
Tensoring \eqref{ref-6.12-51} on the left by the graded $R^{\coord}$-module
$\Omega^\cdot_{R^{\coord}}$ and completing we obtain an isomorphism of
graded commutative algebras.
\begin{equation}
\label{ref-6.22-68}
\Omega^{\cdot}_{R^{\coord}}\cbtimes R\cong \Omega^\cdot_{R^{\coord}}[[t_1,\ldots,t_d]]:
\eta\otimes f\mapsto \eta f
\end{equation}
The DG-algebra structure on
$\Omega^{\cdot}_{R^{{\coord}}}\cbtimes R$ now induces a
DG-algebra structure on $\Omega^\cdot_{R^{\coord}}[[t_1,\ldots,t_d]]$
and thus according to the abstract discussion in \S\ref{ref-6.4-65} 
there is an associated Maurer-Cartan form
\[
\omega_{\text{MC}}\in \Omega^1_{R^{\coord}}\ctimes \frak{g}
\]
such that  for $\eta\in \Omega_{R^{\coord}}^i$, $f\in R$
\[
(d\eta) \tilde{f}=(d+\omega_{\text{MC}})(\eta \tilde{f})
\]
which is equivalent to
\begin{equation}
\label{ref-6.23-69}
(d+\omega_{\text{MC}})(\tilde{f})=0
\end{equation}
The following lemma will be used below.
\begin{lemmas} 
\label{ref-6.5.1-70}
For $v\in \frak{g}$ let $i_{\bar{v}}$ be the contraction
on $\Omega^\cdot_{R^{{\coord}}}$ with  the
  derivation on $R^{{\coord}}$ induced by $v$ (cfr \eqref{ref-6.19-62}). Extend
$i_{\bar{v}}$ to a map of degree $-1$ from
 $\Omega^\cdot_{R^{{\coord}}}\ctimes \frak{g}$  to itself.
Then we have
\begin{equation}
\label{ref-6.24-71}
i_{\bar{v}}\omega_{MC}=1\otimes v
\end{equation}
where both sides are elements of $R^{{\coord}}\ctimes \frak{g}$.
\end{lemmas}
\begin{proof}
  It is easy to see that for any $\omega\in
  \Omega_{R^{\coord}}^1\ctimes \frak{g}$ we have $(i_{\bar{v}}
  \omega)(\tilde{f})=i_{\bar{v}}(\omega(\tilde{f}))$.
  Applying $i_{\bar{v}}$ to \eqref{ref-6.23-69} and using this fact we obtain
\[
(L_{\bar{v}}+i_{\bar{v}}
\omega_{MC})(\tilde{f})=0
\]
Or using \eqref{ref-6.19-62}
\[
(i_{\bar{v}}
\omega_{MC})(\tilde{f})=L_v(\tilde{f})
\]
The operators on both sides are $R^{{\coord}}$-linear. Since $R^{{\coord}}[[t_1,\ldots,t_d]]$ is topologically generated by the $\tilde{f}$, $f\in R$ and
$R^{{\coord}}$ (by the isomorphism \eqref{ref-6.12-51})
we obtain as operators on $R^{\coord}[[t_1,\ldots,t_d]]$
\begin{equation}
\label{ref-6.25-72}
i_{\bar{v}}\omega_{MC}=L_v
\end{equation}
Then \eqref{ref-6.24-71} is the same equation as \eqref{ref-6.25-72} but interpreted in 
$\Omega^\cdot_{R^{{\coord}}} \ctimes \frak{g}$ using faithfulness (see \S\ref{ref-6.4-65}).
\end{proof}
\subsection{An acyclicity result} 
\label{ref-6.6-73}
Assume that $X$ is a  separated smooth $k$-scheme of dimension $d$. Let $\theta_0:X^{\text{aff}}\r X$ be
the canonical map. If $X^{\aff}$ were finite dimensional then the
following result would follow trivially from the theory of algebraic
De Rham cohomology \cite{Hartshorne} together with Proposition \ref{ref-6.3.1-64}.
\begin{theorems}
\label{ref-6.6.1-74}
Put $J=\Jet_{\theta_0,\infty}$. Then the canonical map
\[
\Oscr_{X}\r \pi_{\ast}\Omega^{\cdot,\text{cont}}_{J/X}
\]
is a quasi-isomorphism where $\pi:J\r X$ is the composition of the
map $J\r X\times X^{{\aff}}$ with the projection on the first factor. 
\end{theorems}
\begin{proof}
Since this result is local on $X$ we may assume $X=\Spec R$ and $R$ has
a system of parameters $x_1,\ldots,x_d$. 
Put $x'_i=\theta_0(x_i)\otimes 1-1\otimes x_i$ and

Let $I\subset R^{\text{aff}}\otimes R$ be the kernel of the
multiplication map $R^{\text{aff}}\otimes R\r R^{\aff}$.  Then $J=\Spc
R^{{\aff}}\cbtimes R$ where $R^{{\aff}}\cbtimes R$ is the completion
of $ R^{\text{aff}}\otimes R$ at the ideal $I$ and
\[
\Gamma(X,\pi_{\ast}\Omega^{\cdot,\text{cont}}_{J/X})=\Omega^{\cdot,\text{cont}}_{
R^{\aff}\cbtimes R /R}
\]
We have an $R$-linear isomorphism $R^{\aff}\cong (R^{\text{aff}}\otimes R)/I$
(where on the right hand side $R$ acts on the nose and on the left hand
side it acts via the map $\theta_0$). By Proposition \ref{ref-6.3.1-64} $R^{\aff}/R$ is formally smooth.  Using formal smoothness we obtain a (non-canonical)
$R$-linear splitting of the map
$(R^{\text{aff}}\otimes R)/I^n\r (R^{\text{aff}}\otimes R)/I\cong R^{\aff}$.

Taking the inverse limit over $n$ we obtain a commutative diagram
\begin{equation}
\label{ref-6.26-75}
\begin{CD}
R^{\aff}\cbtimes R @>r\otimes f\mapsto r\theta_0(f)>> R^{{\aff}}\\
@A1\otimes \Id_R AA @AA\theta_0 A\\
R @= R
\end{CD}
\end{equation}
where the top map in \eqref{ref-6.26-75} is (non-canonically) split
as $R$-algebras.

Using Proposition \ref{ref-6.1.3-47} we obtain an isomorphism
\begin{equation}
\label{ref-6.27-76}
R^{\aff}\cbtimes R\cong R^{{\aff}}[[x'_1,\ldots,x'_d]]
\end{equation}
of $R$-algebras.

According to Theorem \ref{ref-6.7.1-79} below the canonical
map
\[
\Omega_{R^{{\aff}}/R}^\cdot\r \Omega^{\text{cont},\cdot}_{R^{{\aff}}[[x'_1,\ldots,x'_d]]/R}
\]
is a quasi-isomorphism. Hence the left inverse of this  map
(coming from the map $R^{\aff}[[x'_1,\ldots,x'_d]]\r R^{\aff}$ given by sending $x'_i\r 0$)
\[
\Omega^{\text{cont},
\cdot}_{R^{{\aff}}[[x'_1,\ldots,x'_d]]/R}\r \Omega_{R^{\aff}/R}^{\cdot}
\]
is also a quasi-isomorphism.

Combining this with \eqref{ref-6.27-76} we see that
the top map in \eqref{ref-6.26-75} induces a quasi-isomorphism on relative
De Rham complexes.

Using the proof of Proposition \ref{ref-6.3.1-64} we see that $R^{\aff}$ is 
a direct limit of finitely generated polynomial rings $R_i$ over $R$.
Thus we have $\Omega^{\cdot}_{R^{{{\aff}}}/R}=\dirlim_i \Omega^{\cdot}_{R_i/R}$ and since
it is well-known that $R\r\Omega_{R_i/R}$ is a quasi-isomorphism
(the Poincare lemma) we obtain that $\theta_0$ induces
a quasi-isomorphism
\[
R=\Omega^{\cdot}_{R/R}\r \Omega_{R^{{{\aff}}}/R}
\]
Thus the right most map in \eqref{ref-6.26-75} also induces a quasi-isomorphism
on relative De Rham complexes. Therefore the left most one does as well.
\end{proof}
\begin{examples} As usual it is instructive to consider the case $R=k[x]$. 
As in Example \ref{ref-6.1.5-53} we have
\[
R^{\coord}=k[x_0,x_1,\dots]_{x_1}
\]
The one dimensional torus $\Gl_1$ acts with weight $-i$ on $x_i$ (this follows form the
fact that $\tilde{x}=\sum_i x_i t^i$ must be invariant). Hence
\[
R^{\aff}=(k[x_0,x_1,\dots]_{x_1})^{\Gl_1}=k[y_0,y_2,\ldots]
\]
where $y_i=(x_1)^{-i} x_i$. The map $R\r R^{\aff}$ is still given by
$x\mapsto x_0=y_0$ and the ideal $I=\ker(R^{\aff}\otimes R\r R^{\aff})$ is generated
by $y_0-x$. 

We have
\[
\Omega_{R^{\aff}}^{\cdot}=k[y_0,dy_0,y_2,dy_2,\ldots]
\]
where $\deg dy_i=1$. Put
\[
\Delta\overset{\text{def}}{=}\Omega_{R^{\aff}}^{\cdot}\otimes R=k[y_0,dy_0,y_2,dy_2,\ldots,x]
\]
and thus
\begin{equation}
\label{ref-6.28-77}
\Omega_{R^{\aff}}^{\cdot}\cbtimes R=k[y_0,dy_0,y_2,dy_2,\ldots,x]\,\hat{}
\end{equation}
where the completion is graded completion with respect to the ideal
$y_0-x$.  To prove directly that the homology of \eqref{ref-6.28-77} is
$R$ it is sufficient to construct a \emph{continuous} homotopy between
the maps of DG-$k[x]$algebras
\[
\phi_0:\Delta \r \Delta:y_i\mapsto \begin{cases} x &\text{if $i=0$}\\
0&\text{otherwise}
\end{cases}, \qquad dy_i\mapsto 0
\]
and 
\[
\phi_1:\Delta\r \Delta:y_i\mapsto y_i, dy_i\mapsto dy_i
\]
(viewed as maps of complexes). 

Introducing an auxiliary variable $z$, a functional homotopy between
these two maps is given by the map of DG-$k[x]$-algebras
\[
H:\Delta\mapsto \Delta\otimes k[z,dz]
\]
\[
H(y_i)=
\begin{cases}
 z(y_0-x)+x,  &\text{if $i=0$}\\
zy_i&\text{otherwise}
\end{cases}
\]
and $H(dy_i)=d(H(y_i))$. By this we mean that
$\phi_i=H\left|_{z=i,dz=0}\right.$.

The following formula then yields a continuous homotopy between $\phi_0$
and $\phi_1$ 
\begin{equation}
\label{ref-6.29-78}
h(\omega)=\int_{z=0}^{z=1}  H(\omega)
\end{equation}
The meaning of the right hand side of \eqref{ref-6.29-78} is as follows. Write
an element $\eta$ of $\Delta\otimes k[z,dz]$ as $\eta_0(x,y,z)+
\eta_1(x,y,z)dz$ where $\eta_0$ does not contain $dz$. Then
\[
\int_{z=0}^{z=1}  \eta\overset{\text{def}}{=}\int_{z=0}^{z=1}\eta_1(x,y,z)dz
\]
\end{examples}
\subsection{De Rham complexes of formal power series rings}
The following abstract result was used in the previous section.
\begin{theorems}
  \label{ref-6.7.1-79} Assume that $T_0$ is an $R$-algebra and that $R$
  is a $k$-algebra. Put
  $\hat{T}=T_0[[x_1,\ldots,x_n]]$.  Then the canonical map
\[
\Omega^{\cdot}_{T_0/R}\r \Omega^{\cdot,\text{cont}}_{\hat{T}/R}
\]
is a quasi-isomorphism of complexes of $R$-modules.
\end{theorems}
\begin{proof} Put $T=T_0[t_1,\ldots,t_n]$. By Proposition \ref{ref-5.4.3-34}
  $\Omega^{\cdot,\text{cont}}_{\hat{T}/R}$ is the (graded) completion of 
$\Omega^\cdot_{{T}/R}$.

 The latter is equal to
$
\Omega^{\cdot}_{T_0/R}\otimes_R \Omega^{\cdot}_{R[t_1,\ldots,t_n]/R}
$
as graded commutative differential graded $R$-algebras.
 We view $\Omega^{\cdot}_{T/R}$ as a first quadrant
double complex with the horizontal direction being given by 
$ \Omega^{\cdot}_{R[t_1,\ldots,t_n]/R}$.

 Hence it is sufficient to prove that for any
$T_0$-module $M$ the completion of 
\begin{equation}
\label{ref-6.30-80}
M\otimes_R \Omega^\cdot_{R[t_1,\ldots,t_n]/R}
\end{equation}
has homology in degree zero and is acyclic elsewhere. By the Poincare
lemma for polynomial rings this is true before completion.
\end{proof}
 Now if we put
$\deg t_i=\deg dt_i=1$ then \eqref{ref-6.30-80} is a graded complex. Hence for every
$n$ the part of degree $n$ 
\begin{equation}
\label{ref-6.31-81}
(M\r M\otimes_R \Omega^\cdot_{R[t_1,\ldots,t_n]/R})_n
\end{equation}
in \eqref{ref-6.30-80} is exact (with $M$ in degree zero). Now since \eqref{ref-6.30-80} is a complex with positively graded
components, its completion (augmented with $M$) is simply the product of the complexes
\eqref{ref-6.31-81}. Hence it is exact also.
\section{Reminder on DG-Lie and $L_\infty$-algebras}
\label{ref-7-82}
\subsection{Coderivations}
\label{ref-7.1-83}
Let $V$ be a graded vector space and set $SV=\oplus_{n=0}^\infty S^n V$ considered
as an augmented coalgebra  such that
\begin{align*}
\Delta(v)&=v\otimes 1+1\otimes v\\
\epsilon(v)&=0
\end{align*}
for $v\in V$.  Fix $w=w_1\cdots w_n$ where the $w_i$ are homogeneous
elements of $V$.  Let $N=\{1,\ldots,n\}$.  For $I\subset
N$ put $w_I=\prod_{i\in I}w_i$. For a disjoint decomposition
$N=I_1\cup\cdots\cup I_p$ we define $\epsilon(I_1,\ldots,I_n)$ as the
sign which makes the following formula formally correct
\[
w_{I_1}\cdots w_{I_p}=
\epsilon(I_1,\ldots,I_p) w
\]
A coderivation $Q$ of degree one on
$SV$ is determined by its ``Taylor coefficients'' $(\partial^n Q)_{n\ge 0}$ which
are the compositions
\[
S^nV\xrightarrow{\text{inclusion}} SV \xrightarrow{Q} SV
\xrightarrow{\text{projection}} V
\]
$Q$ can be computed from its Taylor coefficients by 
a kind of Leibniz rule. One has
\begin{equation}
\label{ref-7.1-84}
Q(w)=\sum_{I\subset N} \epsilon(I,N-I) (\partial^{|I|} Q) (w_I) w_{N-I}
\end{equation}
The Taylor coefficients of $Q^2$ are thus given by
\begin{equation}
\label{ref-7.2-85}
(\partial^n Q^2)(w)=\sum_{I\subset N} \epsilon(I,N-I) 
(\partial^{n-|I|+1}Q)((\partial^{|I|} Q) (w_I) w_{N-I})
\end{equation}

We assume throughout that $Q$ is compatible with the augmented
structure.  I.e.\ $Q(1)=0$, or equivalently $\partial^0 Q=0$. If
$\partial^nQ=0$ for $n>1$ then \eqref{ref-7.1-84} implies that $Q$ is a
\emph{derivation} for the canonical algebra structure on $SV$.
\subsection{Coalgebra maps}
If $V,W$ are graded vector spaces then an augmented coalgebra map of
degree zero $\psi:SV\r SW$ is determined its ``Taylor coefficients''
$(\partial^n \psi)_{\ge 1}$ which are the compositions
\[
S^nV\xrightarrow{\text{inclusion}} SV \xrightarrow{\psi} SW
\xrightarrow{\text{projection}} W
\]
$\psi$ can be computed from its Taylor coefficients as follows. 
\begin{equation}
\label{ref-7.3-86}
\psi(w)=\sum_{N=I_1\cup \cdots \cup I_p} \frac{1}{p!}
\epsilon(I_1,\ldots,I_n) (\partial^{|I_1|} \psi) (w_{I_1})
\cdots
(\partial^{|I_p|} \psi) (w_{I_p})
\end{equation}
Here $N=I_1\cup \cdots \cup I_p$ is an ordered partition of $N$
into $p$ disjoint subsets (with $p$ variable).

It follows from \eqref{ref-7.3-86} that if
$\partial^n\psi=0$ for $n>1$ then $\psi$ is an \emph{algebra
  homomorphism} $SV\r SW$.

Assume that $SV$ and $SW$ are equipped with a coderivation of degree one,
denoted by $Q$. One may show that the condition 
\[
\psi\circ Q=Q\circ \psi
\]
is equivalent to the corresponding ``first order condition''
\[
\partial^n(\psi\circ Q)=\partial^n(Q\circ \psi)
\]
The latter condition maybe expanded as
\begin{multline}
\label{ref-7.4-87}
\sum_{I\subset N} \epsilon(I,N-I) 
(\partial^{n-|I|+1}\psi)((\partial^{|I|} Q) (w_I) w_{N-I})=\\
\sum_{N=I_1\cup \cdots \cup I_p}\frac{1}{p!}
\epsilon(I_1,\ldots,I_n) (\partial^{p} Q)((\partial^{|I_1|} \psi )(w_{I_1})
\cdots
(\partial^{|I_p|} \psi) (w_{I_p})
\end{multline}
For further reference we note that in case $\partial^i Q=0$ for $i\neq 1,2$
this formula specializes to
\begin{multline}
\label{ref-7.5-88}
  \sum_{1\le i\le n} \epsilon(i,N-\{i\}) (\partial^n\psi)((\partial^1
  Q)(w_i)w_{N-\{i\}}) + \\
\sum_{1\le i< j\le n}
  \epsilon(i,j,N-\{i,j\})(\partial^{n-1}\psi)((\partial^2 Q)(w_iw_j)
  w_{N-\{i,j\}})= \\
(\partial^1 Q)((\partial^n \psi)(w))+ \frac{1}{2}\sum_{N=I_1\cup I_2}\epsilon(I_1,I_2)
(\partial^2 Q)((\partial^{|I_1|} \psi)(w_{I_1})(\partial^{|I_2|} \psi)(w_{I_2}))
\end{multline}

\subsection{$L_\infty$-algebras and morphisms}
\begin{definitions}
  An $L_\infty$-structure on a vector space $\frak{g}$ is a coderivation $Q$ of
  degree one on $S(\frak{g}[1])$ which has square zero. 
\end{definitions}
One puts for $a\in \frak{g}$
\begin{equation}
\label{ref-7.6-89}
\begin{split}
da&=-\partial^1 Q(a)\\
[a,b]&=(-1)^{|a|} \partial^2Q(a,b)
\end{split}
\end{equation}
(where $|a|$ is the degree of $a\in \frak{g}$). It then follows from 
\eqref{ref-7.2-85} that $d^2=0$ and that
$d$ is a derivation of degree one of $\frak{g}$ with respect to the
binary operation of degree zero $[-,-]$.  If $\partial^i Q=0$ for
$i>2$ then $\frak{g}$ is a DG-Lie algebra.  Conversely any DG-Lie
algebra can be made into an $L_\infty$-algebra by defining $\partial^1
Q$, $\partial^2 Q$ according \eqref{ref-7.6-89} and by putting
$\partial^iQ=0$ for $i>2$.

\medskip

A morphism of $L_\infty$-algebras $\frak{g}\r \frak{h}$ is by definition
a coalgebra map of degree zero $S(\frak{g}[1])\r S(\frak{h}[1])$ 
commuting with $Q$.  It is customary to write $\psi_i=\partial^i \psi$
where $\psi_i$ is considered as a map $\wedge^i\frak{g}\r \frak{h}$ 
of degree $1-n$. It follows from \eqref{ref-7.4-87} that $d\psi_1=\psi_1 d$.
Hence $\psi_1$ defines a morphism of complexes. 
\subsection{The topological case} The above notions make sense in any 
symmetric monoidal category. We will use them in the case of filtered complete linear
topological vector spaces. 
\subsection{Twisting}
\label{ref-7.5-90}
Assume that $\psi:\frak{g}\r \frak{h}$ is a $L_\infty$-morphism
between $L_\infty$-algebras equipped with some type of topology and
let $\omega\in \frak{g}_1$ be a solution of the $L_\infty$-Maurer-Cartan equation
\[
\sum_{i\ge 1} \frac{1}{i!} (\partial^i Q)(\omega^i)=0
\]
in $\frak{g}$. Here and below we assume that we are in a situation
where all occurring series are convergent and standard series
manipulations are allowed.  This will be the case in the application
in \S\ref{ref-9.1-112} where the series will in fact be finite.

 Define $Q_\omega$, $\psi_\omega$ and $\omega'$ by  \cite{ye3}
\begin{align}
\label{ref-7.7-91}
(\partial^i Q_\omega)(\gamma)&=\sum_{j\ge 0} \frac{1}{j!} (\partial^{i+j} Q)
(\omega^j \gamma)\qquad \text{(for $i>0$)}\\
\label{ref-7.8-92}
(\partial^i \psi_\omega)(\gamma)&=\sum_{j\ge 0} \frac{1}{j!} (\partial^{i+j} \psi)
(\omega^j \gamma)\qquad \text{(for $i>0$)}\\
\label{ref-7.9-93}
\omega'&=\sum_{j\ge 1} \frac{1}{j!} (\partial^{j} \psi)
(\omega^j )
\end{align}
for $\gamma\in S^i(\frak{g}[1])$. Yekutieli shows in \cite{ye3} that
$\omega'$ is a solution of the Maurer-Cartan equation on $\frak{h}$
and that furthermore $\frak{g}$, $\frak{h}$, when equipped with
$Q_\omega$, $Q_{\omega'}$ are again $L_\infty$-algebras.  Let us
denote these by $\frak{g}_\omega$ and $\frak{h}_{\omega'}$.  Yekutieli
also shows that $\psi_\omega$ is an $L_\infty$ map $\frak{g}_\omega\r
\frak{h}_{\omega'}$. Variants of this principle occur at other places
in the literature. See e.g.\ \cite[Corollary
4.0.3]{Tsygan}\cite[\S2.4]{Dolgushev}.
Let us see what the definition of $Q_\omega$ means in case $\frak{g}$ is 
a DG-Lie algebra. In this case the $L_\infty$-Maurer-Cartan equation translates into the usual Maurer-Cartan equation
\[
d\omega+\frac{1}{2}[\omega,\omega]=0
\]
Then
\begin{align*}
(\partial^1 Q_\omega)(\gamma)&=(\partial^1 Q)(\gamma)+ (\partial^2 Q)(\omega \gamma)\\
(\partial^2 Q_\omega)(\gamma)&=(\partial^2 Q)(\gamma)\\
(\partial^i Q_\omega)(\gamma)&=0\qquad (\text{for $i\ge 3$})
\end{align*}
Or translated into differentials and Lie brackets
\begin{equation}
\label{ref-7.10-94}
\begin{split}
d_\omega&=d+ [\omega ,-]\\
[-,-]_\omega&=[-,-]
\end{split}
\end{equation}
\subsection{Descent for $L_\infty$-morphisms}
\label{ref-7.6-95}
This is a somewhat abstract section. It is an explicitation of 
\cite[\S7.3.3]{Ko3} (in particular the last paragraph).
The result will be used to descend
a $L_\infty$-morphism under a rational group action.

We now assume that $\frak{g}$ is a DG-Lie algebra and $\frak{s}$ is a
set. We assume there is an ``action'' 
of $\frak{s}$ on $\frak{g}$
such that $v\in \frak{s}$ acts by a derivation of degree $-1$ on $\frak{g}$, denoted by $i_v$. Put $L_v=di_v+i_vd$. This is a derivation of $\frak{g}$ of
degree zero. 

By the discussion in \S\ref{ref-7.1-83} there exist unique
coderivations $\tilde{\imath}_v$ and $\tilde{L}_v$ on $S(\frak{g}[1])$
such that $\partial^1 \tilde{\imath}_v=j_v\overset{\text{def}}{=}-i_v$
, $\partial^1 \tilde{L}_v=L_v$ and $\partial^i
\tilde{\imath}_v=\partial^i \tilde{L}_v=0$ for $i\neq 1$ (the sign
change on $\tilde{\imath}_v$ occurs because of the fact that $i_v$ is
an odd map $\frak{g}\r \frak{g}$ and $\partial^1 \tilde{\imath}_v$ is
the corresponding  map $\frak{g}[1]\r \frak{g}[1]$).
\begin{lemmas} \label{ref-7.6.1-96} One has
\[
\tilde{L}_v=[Q,\tilde{\imath}_v]
\]
\end{lemmas}
\begin{proof} We know that $[Q,\tilde{\imath}_v]$ is a coderivation.
  Hence we have to compute $\partial^i[Q,\tilde{\imath}_v]$. Since
  $\tilde{\imath}_v$ maps $S^i(\frak{g}[1])$ to $S^i(\frak{g}[1])$ we
      have $\partial^i[Q,\tilde{\imath}_v]=\partial^i Q\circ \tilde{\imath}_v+ 
\tilde{\imath}_v\circ
\partial^i Q$. We need only to consider the cases $i=1,2$.  Assume
first $i=2$.  Then (using the fact that $\tilde{\imath}_v$ is also a
derivation on $S(\frak{g}[1])$, equal to $j_v$ on $\frak{g}[1]$) we
compute for $a,b\in \frak{g}$ (considered as elements of
$\frak{g}[1]$)
\begin{align*}
(\partial^2 Q\circ \tilde{\imath}_v+
\tilde{\imath}_v\circ \partial^2 Q)(a,b)&=\partial^2 Q (j_v a,b)+(-1)^{|a|-1}\partial^2 
Q(a,j_v b)+j_v\circ \partial^2 Q(a,b)\\
&=(-1)^{|a|-1} [j_va , b]-[a,j_v b]+(-1)^{|a|} j_v [a,b]\\
&=0
\end{align*}
(note that $|a|$, $|b|$ refer to the degrees of $a$, $b$ in $\frak{g}$).

Now assume $i=1$. We have
\begin{align*}
(\partial^1 Q\circ \tilde{\imath}_v+\tilde{\imath}_v \circ\partial^1 Q)(a)&=-d j_v(a)-j_v d(a)=L_v(a)=\tilde{L}_v(a)\qed
\end{align*}
\def\qed{}\end{proof}
Put 
\begin{equation}
\label{ref-7.11-97}
\frak{g}^{\frak{s}}=\{X\in \frak{g}\mid \forall v\in \frak{s}:i_v X= L_v X=0\}
\end{equation}
we call $\frak{g}^{\frak{s}}$
 the reduction of $\frak{g}$ with respect to the $\frak{s}$-action.
It is clear that $\frak{g}^{\frak{s}}$ is a DG-Lie algebra as well.
\begin{remarks} It is perhaps useful to point out that whereas the
notion of an $\frak{s}$-action only depends on the graded structure of
$\frak{g}$, the construction of $\frak{g}^{\frak{s}}$ also depends
on the differential. 
\end{remarks}
\begin{propositions} 
\label{ref-7.6.3-98} Assume that $\psi$ is an $L_\infty$-morphism
  $\frak{g}\r \frak{h}$ between DG-Lie algebras equipped with a
  $\frak{s}$-action as above. Assume that $\psi$ commutes with the $\frak{s}$ action in the sense that for all $v\in \frak{s}$
\[
\forall v: [\tilde{\imath}_v,\psi]=0
\]
(where as above $\tilde{\imath}_v$ stands for the induced coderivations on
$S(\frak{g}[1])$ and $S(\frak{h}[1])$).  Then $\psi$ descends to an
$L_\infty$-morphism $\psi^{\frak{s}}:\frak{g}^{\frak{s}}\r
\frak{h}^{\frak{s}}$.
\end{propositions}
\begin{proof} By \eqref{ref-7.11-97} we have to show that the
  restrictions of $\tilde{\imath}_v\circ
  \partial^i \psi$ and $\tilde{L}_v\circ \partial^i \psi$ to
  $S^i(\frak{g}^{\frak{s}}[1])$ are zero.  Note that since $\tilde{L}_v=[Q,
\tilde{\imath}_v]$ (by Lemma \ref{ref-7.6.1-96}, the fact
that $\tilde{\imath}_v$ commutes with $\psi$ implies that $\tilde{L}_v$ commutes with $\psi$ as well. 
\begin{align*}
\tilde{\imath}_v\circ
  \partial^i \psi&=\partial^i(\tilde{\imath}_v \circ \psi)=\partial^i(\psi \circ \tilde{\imath}_v)\\
\tilde{L}_v\circ
  \partial^i \psi&=\partial^i(\tilde{L}_v \circ \psi)=\partial^i(\psi \circ \tilde{L}_v)\
\end{align*}
Since $\tilde{\imath}_v$ and $\tilde{L}_v$ are zero on $S (\frak{g}^{\frak{s}}[1])$ this  implies the desired result.
\end{proof}
\subsection{Compatibility with twisting}
Assume that $\frak{g}$, $\frak{h}$ are topological $L_\infty$-algebras
and $\psi$ is an $L_\infty$-morphism $\frak{g}\r \frak{h}$. We make
the same assumptions as in \S\ref{ref-7.5-90} with regard to convergence
of series. Our aim to understand the behavior of $\frak{s}$-actions under
twisting.
\begin{propositions}
  \label{ref-7.7.1-99} Assume that $\frak{g}$ and $\frak{h}$ are
  equipped with a $\frak{s}$-action and assume that $\psi$ commutes
  with this action (as in Proposition \ref{ref-7.6.3-98}). Let $\omega\in
  \frak{g}_1$ be a solution to the Maurer Cartan equation. 
Since twisting does not change the Lie bracket (see \eqref{ref-7.10-94}), 
$\frak{s}$ acts on  $\frak{g}_\omega$ and $\frak{h}_\omega$ as well.

Assume
  that for $i\ge 2$ and all $v\in \frak{s}$, $\gamma\in
  S^{i-1}(\frak{g}[1])$ we have
\begin{equation}
\label{ref-7.12-100}
(\partial^i \psi)(i_v\omega \cdot \gamma)=0
\end{equation}
Then $\psi_\omega$ is compatible with the $\frak{s}$-action on  $\frak{g}_\omega$ and $\frak{h}_\omega$.
\end{propositions}
\begin{proof} $[\tilde{\imath}_v,\psi_\omega]$ is a
  ``$\psi_\omega$-coderivation'' which is the dual notion of a
  ``$\phi$-derivation'' for a map of algebras $\phi:A\r B$.  One
  verifies that in order to prove $[\tilde{\imath}_v,\psi_\omega]=0$
  it is sufficient to show that $\partial^i
  [\tilde{\imath}_v,\psi_\omega]=\tilde{\imath}_v\circ
\partial^i \psi_\omega-\partial^i \psi_\omega\circ \tilde{\imath}_v=0$.

We have for $\gamma\in S^i(\frak{g}[1])$, $i>0$
\begin{equation}
\label{ref-7.13-101}
  (\partial^i\psi_\omega\circ \tilde{\imath}_v)(\gamma)= \sum_{j\ge 0} \frac{1}{j!}
  (\partial^{i+j}\psi)( \omega^j \cdot \tilde{\imath}_v \gamma)
\end{equation}
and
\begin{equation}
\label{ref-7.14-102}
\begin{split}
  (\tilde{\imath}_v\circ \partial^i\psi_\omega)(\gamma)&= \sum_{j\ge 0} \frac{1}{j!}
  (\partial^{i+j}\psi)( \tilde{\imath}_v(\omega^j  \gamma))\\
&=\sum_{j\ge 1} \left(\frac{1}{(j-1)!}
  (\partial^{i+j}\psi)( j_v(\omega) \cdot \omega^{j-1}  \gamma)\right)+
\left( \frac{1}{j!} (\partial^{i+j}\psi)(\omega^j \cdot \tilde{\imath}_v\gamma) \right)
\end{split}
\end{equation}
(recall that $\omega$ has degree zero in $\frak{g}[1]$). The
difference between \eqref{ref-7.13-101} and \eqref{ref-7.14-102} is a linear
combination of terms of the form \eqref{ref-7.12-100} and hence it is zero.
\end{proof}
\section{Poly-differential operators and poly-vector fields revisited}
\label{ref-8-103}
In this section we remind the reader about some facts on poly-differential
operators and poly-vector fields. These notions were already introduced in
the introduction but for the convenience of the reader we repeat some definitions. From now on we assume $k=\CC$.
\label{ref-8-104}
\subsection{General definitions}
\label{ref-8.1-105}
Let $R$ be a finite adic $k$-algebra. We view $R$ as a $R^{\otimes n}$
module through the diagonal action.  We put
\[
D^{\poly,n}(R)=\Diff_{R^{\otimes n}}(R^{\otimes n},R)
\]
where ``$\Diff$'' stands for differential operators.\footnote{Note
  that since differential operators are continuous with respect to any
  adic topology, we don't have to worry about continuity.} We also
write $D^{\poly,n}(X)=D^{\poly,n}(R)$ if $X=\Spc R$.  If $X$ is a
finite adic scheme then we define $\Dscr_X^{\poly,n}$ by gluing from
the affine case.

We may view the elements of $D^{\poly,n}(R)$ as the set
of multilinear maps $R^{\otimes n}\r R$ which are differential
operators in each of their arguments.  With this interpretation it is
clear that $D^{\poly,\cdot}(R)$ is a DG-Lie subalgebra of 
$\HC^\cdot(R)$, the Hochschild complex of $R$. In particular it is
a  DG-Lie subalgebra.

We say that $R$ is formally of finite type \cite{ye4} if $R$ has
a finitely generated ideal of definition $I$ such that $R/I$ is
finitely generated. This definition is clearly independent of $I$. If $R$ is
formally of finite type then let us say that $R$ is formally smooth if
$\Omega^{1,\cont}_R$ is projective. In that case following
Yekutieli's argument in \cite{ye2} we see that the map
\[
D^{\poly,\cdot}(R)\r \HC^\cdot(R)
\]
is a quasi-isomorphism.

For $p\ge 0$ let $F^p D^{\poly,n}(R)$ be the differential operators of
degree $\le p$.  It is then easy to see that if $R$ is formally of
finite type then $F^pD^{\poly,n}(R)$ is a finite $R$-module. In that
case we will view $D^{\poly,n}(R)$ as a filtered complete linear
topological vector space.

\medskip

Similarly put
\[
T^{\text{poly},n}(R)=\Der_{R^{\otimes n}}(R^{\otimes n},R)^{S^n}
\]
The righthand side describes the set of poly-derivations which are
anti-symmetric in their arguments.\footnote{In \S\ref{ref-5.4-32}
  the notation $\Der^n(R,R)$ was used. The current notation is more
  convenient for this section}

\medskip

If $R$ is formally of finite type then $T^{\text{poly},n}(R)$ is a
finite $R$-module.  In that case we view $T^{\text{poly},n}(R)$ as a filtered
complete linear topological vector space with filtration concentrated
in degree $n$.
\begin{conventionwarnings} In this section and
  the next almost all our objects will be considered as being
  (naturally) filtered. This has serious implications for the meaning of completions
  and completed tensor products. See \S\ref{ref-5.2-28}.
\end{conventionwarnings}

If $R$ formally of finite type and formally smooth then there is an
isomorphism
\[
\wedge^n_R T^{\text{poly,1}}(R)\r T^{\text{poly},n}(R):\gamma_1\wedge\cdots\wedge \gamma_n
\mapsto \sum_{\sigma\in S_n}(-1)^\sigma\gamma_{\sigma(1)}\otimes\cdots\otimes \gamma_{\sigma(n)}
\]
where $\gamma_i\in T^{\text{poly}}(R)=\Der_k(R,R)$ and $\gamma_1\otimes\cdots \otimes \gamma_n$ acts on $R^{\otimes n}$ via
\[
(\gamma_1\otimes\cdots \otimes \gamma_n)(r_1\otimes\cdots\otimes r_n)
=\gamma_1(r_1)\cdots \gamma_n(r_n)
\]
\subsection{The formal case}
\label{ref-8.2-106}
In this section we consider $R=k[[t_1,\ldots,t_d]]$. Write $\partial_i$ for
$\partial/\partial t_i$.  In this case we can give very concrete descriptions
of $T^{\text{poly},\cdot}(R)$ and $D^{\text{poly},\cdot}(R)$. First we have
\begin{equation}
\label{ref-8.1-107}
\begin{split}
T^{\text{poly},1}(R)&=R\partial_1\oplus\cdots\oplus R\partial_d\\
D^{\text{poly},1}(R)&=R[\partial_1,\cdots,\partial_d]
\end{split}
\end{equation}
and then
\begin{equation}
\label{ref-8.2-108}
\begin{split}
T^{\text{poly},n}(R)&=\wedge^n_R T^{\text{poly},n}\\
D^{\text{poly},n}(R)&=\otimes^n_R D^{\text{poly},}
\end{split}
\end{equation}
These descriptions reflect the algebra structure on
$T^{\text{poly},\cdot}(R)$ and $D^{\text{poly},\cdot}(R)$.  It is also
easy to get the Lie algebra structure on $T^{\text{poly},\cdot}(R)$
using the fact that the product satisfies the Leibniz property with
respect to the Lie bracket. In $D^{\text{poly},\cdot}(R)$ this Leibniz
property holds only up to homotopy and therefore the situation is much
more complicated.

Kontsevich (over the reals) 
 constructs in \cite{Ko3} an $L_\infty$-quasi-isomorphism
\begin{equation}
\label{ref-8.3-109}
\Uscr:T^{\text{poly},\cdot}(R)[1]\r D^{\text{poly},\cdot}(R)[1]
\end{equation}
If we write $\Uscr_i=\partial^i\, \Uscr$ then $\Uscr_1$ is given by the
HKR formula
\begin{equation}
\label{ref-8.4-110}
\Uscr_1(\partial_{i_1}\wedge\cdots \wedge\partial_{i_p})= 
\frac{1}{p!}\sum_{\sigma\in S_p}(-1)^{\sigma}\partial_{i_{\sigma(1)}}\otimes \cdots \otimes \partial_{i_\sigma(p)}
\end{equation}
The higher $\Uscr_n$ are matrices of differential operators when expressed
in the natural $R$-bases  of $T^{\text{poly},n}(R)$ and $D^{\text{poly},n}(R)$ obtained
from \eqref{ref-8.1-107}\eqref{ref-8.2-108}.

\medskip

The quasi-isomorphism constructed by Kontsevich has two supplementary 
properties which are crucial for its extension to the non-formal case.
\begin{enumerate}
\item[(P4)]
$\Uscr_q(\gamma_1\cdots \gamma_q)=0$ for  $q\ge 2$ and $\gamma_1,\ldots,\gamma_q
\in T^{\text{poly},1}(R)$. 
\item[(P5)]
$\Uscr_q(\gamma\alpha)=0$ for $q\ge 2$ and $\gamma\in
\frak{gl}_d(k)\subset T^{\text{poly},1}(R)$.
\end{enumerate}
\begin{remarks} 
  In \cite{Tamarkin} Tamarkin constructs an
  $L_\infty$-quasi-isomorphism like \eqref{ref-8.3-109} over the
  rationals.  Halbout has informed me that the methods in
  \cite{Halbout} show that Tarmarkin's quasi-isomorphism may be defined
  in such a way that it also satisfies (P4) and (P5). Using this one
  may replace the complex numbers by the rational numbers in this
  paper.
\end{remarks}
\begin{remarks} Another property which is usually being regarded as essential
for globalization is the fact that the $\Uscr_q$ are $\Gl_d(k)$ equivariant
(condition (P3) in \cite{Ko3}). We will not explicitly use this
condition below. The explanation for this is that (P3) almost follows
from (P5). To be more precise let (P3') be the condition that $\Uscr_q$
is $\frak{gl}_d(k)$ equivariant. I.e. 
\begin{itemize}
\item[(P3')]
$[\gamma,\Uscr_q(\alpha_1\cdots\alpha_q)]=\sum_j\Uscr_q(\alpha_1\cdots
[\gamma,\alpha_j]\cdots \alpha_q)$
\end{itemize}
for $\gamma\in \frak{gl}_d(k)$, $\alpha_i\in T^{\poly,\cdot}(R)[1]$. 
\end{remarks}
Then we have $(P5)\Rightarrow (P3')$. This easily follows from
\eqref{ref-7.5-88}. In sufficiently nice situations (P3) and (P3')
are of course equivalent. 
\section{Global formality}
\label{ref-9-111}
\subsection{Lifting to coordinate spaces}
\label{ref-9.1-112}
It is easy to define relative poly-differential operators with respect
to a graded commutative base ring. Assume that $A\r B$ is a morphism
of graded commutative algebras and let $M$ be a graded $B$-module.  We
define $D^{\poly,n}_{A}(B,M)$ as the set of multilinear maps
$B\otimes_A\cdots \otimes_A B\r M$ ($n$ copies of $B$) which are
relative $B/A$ differential operators in each of their arguments and
which are finite sums of homogeneous maps.

We will use the following routine lemma to manipulate such relative
poly-differential operators.
\begin{lemmas}
  \label{ref-9.1.1-113} Assume that $A$ is a graded commutative DG-algebra and
  let $S$ be a finitely generated smooth $k$-algebra. Let $I$ be a
  finitely generated ideal in $A_0\otimes S$ ($A_0$ is the part of
  degree zero of $A$). Then the obvious map of DG-Lie algebras
\begin{equation}
\label{ref-9.1-114}
 D^{\text{poly},.}(S)\r D^{\text{poly},.}_A(A\cbtimes S)
\end{equation}
(all completions are $I$-adic completions of filtered topological
$A_0\otimes S$-modules) 
extends to an isomorphism of double complexes  of filtered complete
linear topological vector spaces
\begin{equation}
\label{ref-9.2-115}
 A\cbtimes D^{\text{poly},.}(S)\r D^{\text{poly},.}_A(A\cbtimes S)
\end{equation}
 if we define the  vertical differentials as the Hochschild
differentials (considering $A\cbtimes S$ as a graded ring) and the 
horizontal differentials  on the right as $[d_A\otimes 1,-]$ and on the
left as $d_A\otimes 1$.
\end{lemmas}
\begin{proof}
It is easy to see that the differentials are as indicated. So we only
have to show that \eqref{ref-9.2-115} is an isomorphism. 

Since differential operators are always continuous with respect to the 
$I$-adic topology we have 
\[
F^p D^{\text{poly},n}_A(A\cbtimes S)=F^p D^{\text{poly},n}_A(A\otimes S,A\cbtimes S)=
F^p D^{\text{poly},n}(S,A\cbtimes S)
\]
where as in \S\ref{ref-8.1-105} $F^\cdot$ denotes the filtration by
degree of differential operators.

From the standard theory of differential operators 
it follows that there exists a \emph{finitely generated
projective} $S^{\otimes n}$ module $J_p$ such that
\[
F^pD^{\text{poly},n}(S,-)=\Hom_{S^{\otimes n}}(J_p,-)
\]
Hence
\begin{align*}
F^pD^{\text{poly},n}(S,A\cbtimes S)&=\Hom_{S^{\otimes n}}(J_p,A\cbtimes S)\\
&=\Hom_{S^{\otimes n}}(J_p,\invlim_m (A\otimes S)/I^m)\\
&=\invlim_m\Hom_{S^{\otimes n}}(J_p, (A\otimes S)/I^m)\\
&=\invlim_m(A\otimes \Hom_{S^{\otimes n}}(J_p, S))/I^m\\
&=A\cbtimes F^pD^{\text{poly},n}(S)
\end{align*}
In the third line we have used the fact that $J_p$ is finitely generated
projective.  This allows us to replace $J_p$ by $S^{\otimes n}$.
\end{proof} 
Assume now that $R$ is smooth of dimension $d$.  Using the previous
lemma first with $A=\Omega_{R^\text{coord}}^\cdot$, $S=R$ and $I=\ker
(R^{\coord}\otimes R\r R^{\coord})$ and then with
$A=\Omega_{R^\text{coord}}^\cdot$, $S= k[t_1,\ldots,t_d]$,
$I=(t_1,\ldots,t_d)$ we have the following string of maps between of
filtered complete linear topological DG-Lie algebras.
\begin{equation}
\label{ref-9.3-116}
\begin{split}
D^{\text{poly},\cdot}(R)&\r  \Omega_{R^\text{coord}}^\cdot\cbtimes D^{\text{poly},\cdot}(R)\\
&\cong
D^{\text{poly},\cdot}_{\Omega_{R^\text{coord}}^\cdot}(\Omega^\cdot_{R^\text{coord}}\cbtimes R)\\
&\cong D^{\text{poly},\cdot}_{\Omega_{R^\text{coord}}^\cdot}(
\Omega_{R^\text{coord}}^\cdot[[t_1,\ldots,t_d]])\\
&\cong\Omega_{R^\text{coord}}^\cdot \,\hat{\boxtimes}\, D^{\text{poly},\cdot}(
k[t_1,\ldots,t_d])\\
&\cong\Omega_{R^\text{coord}}^\cdot \ctimes D^{\text{poly},\cdot}(
k[[t_1,\ldots,t_d]])
\end{split}
\end{equation}
In the third line we have used the isomorphism \eqref{ref-6.22-68}.

\medskip

In \S\ref{ref-6.5-67} we have seen that the differential on
$\Omega^{\cdot,\cont}_{R^{\coord}\cbtimes R/R}=\Omega_{R^\text{coord}}^\cdot\cbtimes R$ induces the differential
$d+\omega_{MC}$ on $\Omega_{R^\text{coord}}^\cdot[[t_1,\ldots,t_d]]$ with
\begin{equation}
\label{ref-9.4-117}
\omega_{MC}\in \Omega^1_{R^{\coord}}\ctimes \Der_k(k[[t_1,\ldots,t_d]])
=\Omega^1_{R^{\coord}}\ctimes T^{\poly,1}(k[[t_1,\ldots,t_d]])\subset
\Omega^1_{R^{\coord}}\ctimes D^{\poly,1}(k[[t_1,\ldots,t_d]])
\end{equation}

This then induces the differential
$[d,-]+[\omega_{MC},-]+d_{\text{Hoch}}$ on $
D^{\text{poly},\cdot}_{\Omega_{R^\text{coord}}^\cdot}(
\Omega_{R^\text{coord}}^\cdot[[t_1,\ldots,t_d]])$. It is easy to see
that this induces the differential $d\otimes
1+[\omega_{MC},-]+d_{\text{Hoch}}$ on the graded Lie algebra
$\Omega_{R^\text{coord}}^\cdot \ctimes
D^{\text{poly},\cdot}( k[[t_1,\ldots,t_d]])$.

\medskip

As for differential operators, we can define relative poly-vector fields
and there is an obvious analog of lemma \ref{ref-9.1.1-113}.
We obtain maps of graded vector spaces which we may employ to get the
following string of maps between filtered complete linear topological DG-Lie algebras.
\begin{equation}
\label{ref-9.5-118}
\begin{split}
T^{\text{poly},\cdot}(R)&\r \Omega_{R^\text{coord}}^\cdot\cbtimes T^{\text{poly},\cdot}(R)
\\
&\r T^{\text{poly},\cdot}_{\Omega_{R^\text{coord}}^\cdot}
(\Omega_{R^\text{coord}}^\cdot\cbtimes R)\\
&\cong  T^{\text{poly},\cdot}_{\Omega_{R^\text{coord}}^\cdot}
(\Omega_{R^\text{coord}}^\cdot[[t_1,\ldots,t_d]])\\
&\cong\Omega_{R^\text{coord}}^\cdot\, \hat{\boxtimes}\, T^{\text{poly},\cdot}(
k[t_1,\ldots,t_d])\\
&\cong \Omega_{R^\text{coord}}^\cdot\ctimes T^{\text{poly},\cdot}
(k[[t_1,\ldots,t_d]])
\end{split}
\end{equation}
On the first Lie algebra this differential is trivial and on the last
it is $d\otimes 1+[\omega_{MC},-]$.

\medskip

Now consider the following two DG-Lie algebra
\begin{equation}
\frak{t}=\Omega_{R^{\coord}}^\cdot\ctimes T^{\text{poly},\cdot}_{k[[t_1,\cdots,t_d]]}
(k[[t_1,\ldots,t_d]])
\end{equation}
\begin{equation}
\frak{d}=\Omega_{R^{\coord}}^\cdot \ctimes D^{\text{poly},\cdot}_{k[[t_1,\ldots,t_d]]}(
k[[t_1,\ldots,t_d]])
\end{equation}
where the DG-Lie algebra structures are obtained by linearly
extending the ones on the second component.   Looking at the second line
and the last line in \eqref{ref-9.3-116} and \eqref{ref-9.5-118} and the
description of the differentials we obtain by \eqref{ref-7.10-94}
\begin{equation}
\label{ref-9.8-119}
\begin{split}
 \frak{t}_{\omega_{MC}}&=\Omega_{R^\text{coord}}^\cdot\cbtimes T^{\text{poly},\cdot}(R)
\\
 \frak{d}_{\omega_{MC}}&=\Omega_{R^\text{coord}}^\cdot\cbtimes D^{\text{poly},\cdot}(R)
\end{split}
\end{equation}
We now extend the $L_\infty$ quasi-isomorphism 
\[
\Uscr :T^{\text{poly},\cdot}
(k[[t_1,\ldots,t_d]])\r  D^{\text{poly},\cdot}(
k[[t_1,\ldots,t_d]])
\]
$\Omega^\cdot_{R^{\coord}}$-linearly to an $L_\infty$-map
$\bar{\Uscr}:\frak{t}\r \frak{d}$.

Property (P4) together with \eqref{ref-9.4-117} implies that if we plug
$\omega_{MC}$ into \eqref{ref-7.9-93} we obtain $\omega'_{MC}=\omega_{MC}$ (taking
into account that the right most inclusion in \eqref{ref-9.4-117} is simply $\bar{\Uscr}_1$). 
\begin{lemmas}
\label{ref-9.1.2-120}
Assume that $\gamma\in \Omega^a_{R^{\coord}}\cbtimes
T^{\poly,b}(R)$. Then the sum in \eqref{ref-7.8-92} is finite.
\end{lemmas}
\begin{proof}
  This follows by degree considerations. Indeed $\bar{\Uscr}_i$ is
  obtained by extension of the map of degree zero
\[
\Uscr_i:S^i(T^{\poly,\cdot}(k[[t_1,\ldots,t_d]])[2])
\r D^{\poly,\cdot}(k[[t_1,\ldots,t_d]])[2]
\]
It follows that $\bar{\Uscr}_i$ considered as a map
\[
S^i(\Omega^\cdot_{R^{\coord}}\ctimes T^{\poly,\cdot}(k[[t_1,\ldots,t_d]]))
\r \Omega^\cdot_{R^{\coord}}\ctimes D^{\poly,\cdot}(k[[t_1,\ldots,t_d]])
\]
has bidegree $(0,2-2i)$. One the other hand it follows from \eqref{ref-9.4-117}
that $\omega_{MC}$ has bidegree $(1,1)$. Assume that $\gamma$ has 
bidegree $(a,b)$. Then the $j$'th term in \eqref{ref-7.8-92} has bidegree
$(a,b)+(j,j)+(0,2-2(i+j))=(a+j,b-j+2-2i)$.  Since for $j\gg 0$ we have
$b-j+2-2i<0$ it follows that the sum in \eqref{ref-7.8-92} is indeed finite
(for a given bihomogeneous $\gamma$). 
\end{proof}
Hence we obtain that
\begin{equation}
\label{ref-9.9-121}
\bar{\Uscr}_{\omega_{MC}}:\frak{t}_{\omega_{MC}}\r \frak{d}_{\omega_{MC}}
\end{equation}
is defined. From the formula \eqref{ref-7.8-92} it follows that $\bar{\Uscr}_{\omega_{MC}}$ is still $\Omega^\cdot_{R^{\coord}}$-linear. 
\subsection{Descent}
 If we let
$v\in \frak{g}=\Der_k(k[[t_1,\ldots,t_d]])$ act by $i_{\bar{v}}$ on 
$\Omega^\cdot_{R^{\text{coord},\cdot}}$ (as in Lemma \ref{ref-6.5.1-70}) 
then this defines a
$\frak{g}$-action on $\Omega^\cdot_{R^{\text{coord}}}$ in the
sense of \S\ref{ref-7.6-95}.  We define a corresponding $\frak{g}$-action on
$\frak{t}$ and $\frak{d}$ by linearly extension.
\begin{lemmas}
\label{ref-9.2.1-122}
Put $\frak{s}=\mathrm{Lie} \Gl_d=\frak{gl}_d\subset \frak{g}$. 
The $L_\infty$-morphism
\eqref{ref-9.9-121} descends to an $L_\infty$-morphism
\[
(\bar{\Uscr}_{\omega_{MC}})^{\frak{s}}:(\frak{t}_{\omega_{MC}})^{\frak{s}}\r
(\frak{d}_{\omega_{MC}})^{\frak{s}}
\]
(where $(-)^{\frak{s}}$ is defined by \eqref{ref-7.11-97}).
\end{lemmas}
\begin{proof} 
  According to Proposition \ref{ref-7.6.3-98} we need to show that the
  $\frak{s}$ action on $\frak{t}_{\omega_{MC}}$ and
  $\frak{d}_{\omega_{MC}}$ is compatible with $\bar{\Uscr}_{\omega_{MC}}$. By
    Proposition \ref{ref-7.7.1-99} it is sufficient to prove the
    following two statements
\begin{enumerate}
\item The $\frak{s}$ action on $\frak{t}$ and
  $\frak{d}$ is compatible with $\bar{\Uscr}$. 
\item For $j\ge 2$ the condition 
\[
\bar{\Uscr}_j(i_{\bar{v}}\omega\cdot \gamma)=0
\]
is satisfied. 
\end{enumerate}
Since $\bar{\Uscr}$ is a base extension of $\Uscr$,
it is easy to see that it commutes with the action of $i_{\bar{v}}$,
$v\in \frak{s}$ (in the sense of Proposition \ref{ref-7.6.3-98}). This proves (1). 

Using Lemma \ref{ref-6.5.1-70} and expanding $\gamma$ as a $\Omega^\cdot_{R^{\coord}}$-linear  combination of elements $\gamma'$ of $S(\frak{t}[1])$ it
follows that for (2) it is sufficient to prove that 
\[
\Uscr_j(v\cdot \gamma')=0
\]
But this is precisely property (P5). 
\end{proof}
\begin{lemmas}
\label{ref-9.2.2-123} The formulas \eqref{ref-9.8-119} descend to  isomorphisms of filtered complete linear
topological DG-Lie algebras
\[
(\frak{d}{\omega_{MC}})^{\frak{s}}\cong \Omega^\cdot_{R^{\aff}}\cbtimes D^{\poly,.}(R)
\]
\[
(\frak{t}_{\omega_{MC}})^{\frak{s}}\cong \Omega^\cdot_{R^{\aff}}\cbtimes T^{\poly,.}(R)
\]
where the completion is computed with respect to the ideal
$\ker(R^{\aff}\otimes R\r R^{\aff})$.
\end{lemmas}
\begin{proof}
We will concentrate ourselves on $\frak{d}$. The case
of $\frak{t}$ is similar. By \eqref{ref-9.8-119} we
have
\[
(\frak{d}_{\omega_{MC}})^{\frak{s}}=\bigl( \Omega_{R^\text{coord}}^\cdot\cbtimes D^{\text{poly},\cdot}(R) \bigr)^{\frak{s}}
\]
Thus we have to understand how the $\frak{s}$-action on
$\Omega_{R^\text{coord}}^\cdot \ctimes D^{\text{poly},\cdot}(
k[[t_1,\ldots,t_d]])$ is transported under the isomorphisms of
\eqref{ref-9.3-116} to a $\frak{s}$-action on $
\Omega_{R^\text{coord}}^\cdot\cbtimes D^{\text{poly},\cdot}(R)$.  We
will do this for the $\frak{g}$-action. Since $\frak{s}\subset \frak{g}$
this does what we want. 

\medskip

 We claim that the transported $\frak{g}$-action is just the extension of
the $\frak{g}$-action on $ \Omega_{R^\text{coord}}^\cdot$.
To prove this we observe that the
isomorphism \eqref{ref-6.12-51}
\[
\Omega_{R^{\text{coord}}}^{\cdot}\cbtimes R\r \Omega_{R^{\text{coord}}}^\cdot[[t_1,\ldots,t_d]]
: \omega\otimes f\mapsto \omega\tilde{f}
\]
commutes with the $\frak{g}$-actions on both sides (obtained from
extending the $\frak{g}$ action on $\Omega_{R^{\text{coord}}}^\cdot$).
To see this note that since $i_{\bar{v}}$ has degree $-1$,
$i_{\bar{v}}$ is $R^{\coord}\cbtimes R$-linear on the left and
$R^{\coord}[[t_1,\ldots,t_d]]$-linear on the right.

If $v\in \frak{g}$ then  $i_{\bar{v}}$ acts as a derivation on
$\Omega_{R^{\text{coord}}}^\cdot\cbtimes R$ and $\Omega_{R^{\text{coord}}}^\cdot[[t_1,\ldots,t_d]]$
which preserves $\Omega_{R^{\coord}}^\cdot$ and this implies that
$[i_{\bar{v}},-]$ acts on
\begin{equation}
\label{ref-9.10-124}
D^{\text{poly},\cdot}_{ \Omega_{R^{\text{coord}}}^\cdot} (\Omega_{R^{\text{coord}}}^\cdot\cbtimes R)
\end{equation}
and
\begin{equation}
\label{ref-9.11-125}
D^{\text{poly},\cdot}_{  \Omega_{R^{\text{coord}}}^\cdot}(\Omega_{R^{\text{coord}}}^\cdot[[t_1,\ldots,t_d]])
\end{equation}
and of course these actions are compatible with the isomorphism
between \eqref{ref-9.10-124} and \eqref{ref-9.11-125}. 

It follows that we obtain compatible $\frak{g}$-actions on
\[
\Omega_{R^{\text{coord},\cdot}}\cbtimes D^{\text{poly},\cdot}(R)
\]
and
\[
\Omega_{R^{\text{coord},\cdot}}\ctimes D^{\text{poly},\cdot}(k[[t_1,\ldots,t_d]])
\]
and it is easy to see that these are obtained from the $\frak{g}$-actions
on $\Omega_{R^{\text{coord},\cdot}}$.

\medskip

It remains to show
\[
(\Omega_{R^{\coord,\cdot}}\cbtimes D^{\text{poly},\cdot}(R))^{\frak{s}}=
\Omega_{R^{\aff,\cdot}}\cbtimes D^{\text{poly},\cdot}(R)
\]
This meant to be an isomorphism of filtered objects so we first
consider $(\Omega_{R^{\coord,\cdot}}\cbtimes F^p D^{\text{poly},\cdot}(R))^{\frak{s}}$ which can be rewritten as
\[
((\Omega_{R^{\coord,\cdot}}\cbtimes R)\otimes_R F^p D^{\text{poly},\cdot}(R))^{\frak{s}}
=(\Omega_{R^{\coord,\cdot}}\cbtimes R)^{\frak{s}}\otimes_R D^{\text{poly},\cdot}(R)
\]
since the $F^p D^{\text{poly},\cdot}(R)$ are finitely
generated projective $R$-modules.

So now we have to show
\[
(\Omega_{R^{\coord}}^\cdot\cbtimes R)^{\frak{s}}=\Omega_{R^{\aff}}^\cdot\cbtimes R
\]
We use the easily proved fact that $\Omega_{R^{\coord}}^\cdot\cbtimes
R= \Omega_{R^{\coord}\cbtimes R/R}^{\cdot,\text{cont}}$. Let $I$ be
the kernel of $R^{\coord}\otimes R\r R^{\coord}$. Using Proposition
\ref{ref-5.4.1-33} and Lemma \ref{ref-9.2.3-126} below we have
\begin{align*}
\left(\Omega_{R^{\coord}\cbtimes R/R}^{m,\text{cont}}\right)^{\frak{s}}&
=\left(\invlim_n \Omega_{(R^{\coord}\otimes R)/I^n/R}^{m}\right)^{\frak{s}}\\
&=\invlim_n \left( \Omega_{(R^{\coord}\otimes R)/I^n/R}^{m}\right)^{\frak{s}}\\
&=\invlim_n \left( \Omega_{((R^{\coord}\otimes R)/I^n)^{\Gl_d}/R}^{m}\right)\\
\end{align*}
Now let $J$ be the kernel of $R^{\aff}\otimes R\r R^{\aff}$.  From the
fact that $\Gl_d$ acts freely on $R^{\coord}$ (Proposition
\ref{ref-6.2.2-58}) and its invariants are defined as $R^{\aff}$ we
easily deduce that $((R^{\coord}\otimes
R)/I^n)^{\Gl_d}=(R^{\aff}\otimes R)/J^n$. Hence
\begin{align*}
\left(\Omega_{R^{\coord}\cbtimes R/R}^{m,\text{cont}}\right)^{\frak{s}}&=
\invlim_n \left( \Omega_{((R^{\aff}\otimes R)/J^n)/R}^{m}\right)\\
&=\Omega_{(R^{\aff}\cbtimes R)/R}^{m,\cont}
\end{align*}
where we have used Proposition \ref{ref-5.4.1-33} once again.
\end{proof}
The following lemma was used.
\begin{lemmas} 
  \label{ref-9.2.3-126} Let $S$ be a connected reductive algebraic group over
  $k$ with Lie algebra $\frak{s}$ acting rationally, $R$-linearly and
  freely on a $R$-algebra $T$. For $v\in \frak{s}$ we denote by $i_v$
  the derivation of degree $-1$ on $\Omega^\cdot_{T/R}$ which is the
  contraction with the derivation corresponding to $v$. Then
\begin{align*}
\Omega_{T^S/R}^\cdot&
=\{\omega\in \Omega^\cdot_{T/R}\mid \forall v\in \frak{s}:
i_v(\omega)=i_v(d\omega)=0\}
\\&=(\Omega_{T/R}^\cdot)^{\frak{s}}
\end{align*}\end{lemmas}
\begin{proof} We use a fragment of the Cartan model for equivariant
  cohomology.  Let $(e_j)_j$ be a basis for $\frak{s}$ and let
  $(e^\ast_j)_j$ be the corresponding dual basis.  Since $S$ acts
  freely on $T$, $T/T^S$ is smooth. We obtain an exact sequence
\[
0\r \Omega^1_{T^S/T}\otimes_{T^G} T\r \Omega^1_{T/R} 
\xrightarrow{\sum_j i_{e_j}\otimes e^\ast_j} T\otimes \frak{s}^\ast \r 0
\]
This sequence is split and hence we may transform into a Koszul type
long exact sequence
\begin{equation}
\label{ref-9.12-127}
0\r\Omega^i_{T^S/R}\otimes_{T^S} T\r
\Omega^{i}_{T/R}\xrightarrow{\delta}\Omega^{i-1}_{T/R} \otimes
\frak{s}^\ast\xrightarrow{\delta} \Omega^{i-2}_{T/R} \otimes
S^2\frak{s}^\ast\xrightarrow{\delta} \cdots
\end{equation}
where 
\[
\delta(\omega\otimes f)=\sum_j i_{e_j}(\omega)\otimes e^\ast_j f
\]
Taking invariants we obtain in particular
\[
\Omega_{T^S/R}^\cdot=\{\omega\in (\Omega_{T/R}^\cdot)^S\mid \forall v\in \frak{s}:
i_v(\omega)=0\}
\]
The differentiated $S$ action of $\Omega_T^\cdot$ is given by $L_v=di_v+i_vd$.
Hence we obtain
\[
\Omega_{T^S/R}^\cdot=\{\omega\in (\Omega_{T/R}^\cdot)\mid \forall v\in \frak{s}:
i_v(\omega)=L_v(\omega)=0\}
\]
which yields the desired result. 
\end{proof}
\subsection{Quasi-isomorphisms}
\label{ref-9.3-128}
\begin{theorems} \label{ref-9.3.1-129} The canonical maps
\begin{align*}
 T^{\text{poly,.}}(R)&\r \Omega_{R^{\text{aff}}}\cbtimes T^{\text{poly,.}}(R)\\
 D^{\text{poly,.}}(R)&\r \Omega_{R^{\text{aff}}}\cbtimes D^{\text{poly,.}}(R)\\
\end{align*}
obtained by linearly extending $R^{\aff}\r \Omega^\cdot_{R^{\aff}}$,
are filtered quasi-isomorphisms.
\end{theorems}
\begin{proof}
By acyclicity (Theorem \ref{ref-6.6.1-74}) we have a quasi-isomorphism
\[
R\r \Omega^{\cdot, \text{cont}}_{R^{\text{aff}}\cbtimes R/R}
\]
and since 
\[
\Omega^{\text{cont},\cdot}_{R^{\text{aff}}\cbtimes R/R}\cong \Omega_{R^{\text{aff}}}^\cdot\cbtimes R
\]
we obtain a quasi-isomorphism
\[
R\r  \Omega_{R^{\text{aff}}}^\cdot\cbtimes R
\]
Tensoring on the right by  the finitely generated projective $R$-modules
$F^pT^{\text{poly},\cdot}(R)$ and $F^pD^{\text{poly},\cdot}(R)$ gives what we want.
\end{proof}
\subsection{Tying it all together}
We have a commutative diagram of DG-Lie algebras and (vertical) $L_\infty$-maps
\begin{equation}
\label{ref-9.13-130}
\begin{CD}
T^{\text{poly,.}}(R)@>>> \Omega_{R^{\text{aff}}}\cbtimes T^{\text{poly,.}}(R) 
@>>>\Omega_{R^{\text{coord}}}\cbtimes T^{\text{poly,.}}(R) 
@>\cong>>\Omega_{R^{\text{coord}}}\ctimes T^{\text{poly,.}}(k[[t_1,\dots,t_d]])\\
@.                          @V \Vscr^{\frak{s}}VV @V\Vscr VV @VV\overline{\Uscr}_{\omega_{MC}}V \\
D^{\text{poly,.}}(R)@>>> \Omega_{R^{\text{aff}}}\cbtimes D^{\text{poly,.}}(R)
@>>>\Omega_{R^{\text{coord}}}\cbtimes D^{\text{poly,.}}(R) 
@>>\cong>\Omega_{R^{\text{coord}}}\ctimes D^{\text{poly,.}}(k[[t_1,\dots,t_d]])
\end{CD}
\end{equation}
where $\Vscr$ is obtained from $\overline{\Uscr}_{\omega_{MC}}$ using the
horizontal isomorphisms and
$\Vscr^{\frak{s}}$ is obtained from $(\overline{\Uscr}_{\omega_{MC}})^{\frak{s}}$
(see Lemmas \ref{ref-9.2.1-122},\ref{ref-9.2.2-123}). By Theorem \ref{ref-9.3.1-129} we know
that the left most horizontal maps are quasi-isomorphisms. 
\begin{theorems}
\label{ref-9.4.1-131}  The induced map
\[
\mu:T^{\text{poly,.}}(R)\r H^\cdot(D^{\text{poly,.}}(R))
\]
is an isomorphism. If $R$ has a system of parameters $(x_i)_i$ then 
\begin{equation}
\label{ref-9.14-132}
\mu(\partial_{i_1}\wedge \cdots \wedge \partial_{i_n})=
\frac{1}{n!}\sum_{\sigma \in S_n}
(-1)^\sigma\partial_{i_{\sigma(1)}}\otimes \cdots \otimes \partial_{i_{\sigma(n)}}
\end{equation}
where $\partial_i=\partial/\partial x_i$.
\end{theorems}
\begin{proof} 
Since everything is local on $R$ we may assume that $R$
  has a system of parameters. Denote the map defined by \eqref{ref-9.14-132} by $\mu'$.
Since $\mu'$
  is a quasi-isomorphism \cite{ye2} it is sufficient to prove that
  $\mu=\mu'$.

Let us regard the complexes occurring in the \eqref{ref-9.13-130}
as double complexes
such that the rows are obtained from the De Rham complexes.  Assume
$\gamma\in ( \Omega_{R^{\text{aff}}}\ctimes T^{\text{poly,.}}(k[[t_1,\dots,t_d]]))_{pq}$ where $p$ is
the column index. 

According to \eqref{ref-7.8-92} $\bar{\Uscr}_{\omega_{MC},1}$ is given by
\[
\bar{\Uscr}_{\omega_{MC},1}(\gamma)=\sum_{j\ge 0} \frac{1}{j!}\, \bar{\Uscr}_{j+1}
(\omega^j_{MC} \gamma)
\]
and one checks
\begin{equation}
\label{ref-9.15-133}
 \bar{\Uscr}_{j+1}
(\omega^j_{MC} \gamma)\in (\Omega_{R^{\coord}}\ctimes  D^{\text{poly,.}}(k[[t_1,\dots,t_d]]))_{p+j,q-j}
\end{equation}
(see for example the proof of Lemma \ref{ref-9.1.2-120}).
Since the horizontal maps in \eqref{ref-9.13-130} are inclusions we obtain that
$\Vscr^{\frak{s}}_1$ maps $(\Omega_{R^{\text{aff}}}^\cdot\cbtimes T^{\text{poly,.}}(R))_{pq}$ to
$\oplus_j (\Omega_{R^{\aff}}^\cdot\cbtimes D^{\poly,.}(R))_{p+j,q-j}$.

We claim that the component corresponding to $j=0$ of
$\Vscr^{\frak{s}}_1$ is equal to (the linear extension of) $\mu'$.  To
prove this we look at the component of $\bar{\Uscr}_{\omega_{MC},1}$
corresponding to $j=0$, which is equal to $\bar{\Uscr}_1$.  As
discussed in \S\ref{ref-9.1-112} the vertical arrows in \eqref{ref-9.13-130}
are linear for the action of the De Rham complexes. Hence it suffices to prove that the following diagram is
commutative
\begin{equation}
\label{ref-9.16-134}
\begin{CD}
  T^{\text{poly},n}(R) @>i>> R^{\text{coord}}\ctimes T^{\poly,n}(k[[t_1,\ldots,t_d]])\\
  @V \mu' VV @VV\bar{\Uscr}_1 V\\
  D^{\text{poly},n}(R)@>>j>
  R^{\text{coord}}\ctimes D^{\poly,n}(k[[t_1,\ldots,t_d]])
\end{CD}
\end{equation}
For convenience we have denoted the horizontal arrows by $i$ and $j$. They
are obtained from the ``expansion in local coordinates'' isomorphism \eqref{ref-6.12-51}
\[
R^{\coord}\cbtimes R\xrightarrow{\cong} R^{\coord} \ctimes k[[t_1,\ldots,t_d]]
\]
Thus a poly-differential operator or vector field on $R$ is linearly extended
to one on $R^{\coord}\cbtimes R$ and then transported to an operator
on $R^{\coord} \ctimes k[[t_1,\ldots,t_d]]$. Clearly $i$ and $j$ are compatible
with cup-product.

To avoid confusion we write $\partial_{x_i}$ for $\partial/\partial
x_i$ and $\partial_{t_j}$ for $\partial/\partial t_j$.  According to for example the proof of
Theorem \ref{ref-6.1.4-48} the matrix $(\partial x_i/\partial t_j)_{ij}$
is an invertible matrix over $R^{\coord}[[t_1,\ldots,t_d]]$. Denote
the inverse matrix by $(\partial t_j/\partial x_i)_{ij}$. Then
\begin{align*}
  i( \partial_{x_{i_1}}\wedge\cdots \wedge \partial_{x_{i_n}}
  )=\sum_{j_1,\ldots,j_n}\frac{\partial t_{j_1}}{\partial x_{i_1}}\cdots \frac{\partial
    t_{j_n}}{\partial x_{i_n}} \partial_{t_{j_1}}\wedge\cdots \wedge
  \partial_{t_{j_n}}
  \\
  j( \partial_{x_{i_1}}\otimes\cdots \otimes\partial_{x_{i_n}}
  )=\sum_{j_1,\ldots,j_n}\frac{\partial t_{j_1}}{\partial x_{i_1}}\cdots \frac{\partial
    t_{j_n}}{\partial x_{i_n}} \partial_{t_{j_1}}\otimes\cdots \otimes 
  \partial_{t_{j_n}}
\end{align*}
Comparing the formulas \eqref{ref-8.4-110} and \eqref{ref-9.14-132} (which we have
taken to define $\mu'$) we see that \eqref{ref-9.16-134} is indeed
commutative.

The first component of an $L_\infty$-map always commutes with the differential
thus we have a  map of complexes
\[
\Vscr^{\frak{s}}_1:\Omega_{R^{\aff}}^\cdot\cbtimes T^{\text{poly},\cdot}(R)
\r \Omega_{R^{\aff}}^\cdot\cbtimes D^{\text{poly},\cdot}(R)
\]
We filter the two complexes according to the column index. 
By \eqref{ref-9.15-133}
this filtration
is compatible with $\Vscr^{\frak{s}}_1$ and the graded map associated
to $\Vscr^{\frak{s}}_1$ is $(\Vscr^{\frak{s}}_1)_{j=0}$, which we have shown
to be equal to the linear extension of $\mu'$.

Denote by $H^{\text{columns}}$ the homology of the columns of a double complex.  We clearly have 
\[
H^{\text{columns}}(\Omega_{R^{\aff}}^\cdot \cbtimes T^{\text{poly},\cdot}(R))=
\Omega_{R^{\aff}}^\cdot \cbtimes T^{\text{poly},\cdot}(R)
\]
and since $D^{\text{poly},\cdot}(R)$ consists of filtered projective $R$-modules
with filtered projective homology we also have
\[
H^{\text{columns}}(\Omega_{R^{\aff}}^\cdot \cbtimes D^{\text{poly},\cdot}(R))
=\Omega_{R^{\aff}}^\cdot \cbtimes H^\cdot( D^{\text{poly},\cdot}(R))
\]
Taking homology for the rows (and using Theorem \ref{ref-9.3.1-129}) we see that $\Vscr^{\frak{s}}$ induces indeed
$\mu'$ on homology. 
\end{proof}
\subsection{The global case}
\label{ref-9.5-135}
\begin{theorems}  \label{ref-9.5.1-136} There exists a sheaf of
  DG-Lie algebras $\frak{l}[1]$ on $X$ together with $L_\infty$
  morphisms
\[
\Tscr^{\text{poly},\cdot}_X[1] \r \frak{l}[1]\l \Dscr^{\text{poly},\cdot}_X[1]
\]
Furthermore
\begin{enumerate}
\item
$\frak{l}$ as well as the given quasi-isomorphisms do not depend on any
choices.\footnote{Except for the choice of the $L_\infty$-quasi-isomorphism in the formal case
satisfying the properties (P4)(P5)}
\item If $X$ has a system of parameters $(x_i)_i$ then the resulting map 
on homology
\[
\Tscr^{\text{poly},\cdot}_X\r H^\cdot(\Dscr^{\text{poly},\cdot}_X)
\]
is given by the HKR-formula.
\[
\partial_{i_1}\wedge \cdots \wedge \partial_{i_n}\mapsto
\frac{1}{n!}\sum_{\sigma \in S_n}
(-1)^\sigma\partial_{i_{\sigma(1)}}\otimes \cdots \otimes \partial_{i_{\sigma(n)}}
\]
where $\partial_i=\partial_{x_i}$.
\end{enumerate}
\end{theorems}
\begin{proof}
Since all our constructions are
canonical we may assume that $X=\Spec R$ where $R$ is smooth of
dimension $d$.

The diagram \eqref{ref-9.13-130} in combination with Theorem
\ref{ref-9.4.1-131} furnishes us with $L_\infty$-quasi-isomorphisms
\[
T^{\text{poly},\cdot}(R)\r  \Omega_{R^{\text{aff}}}\cbtimes T^{\text{poly,.}}(R)
\xrightarrow{\Vscr^{\frak{s}}}  \Omega_{R^{\text{aff}}}\cbtimes D^{\text{poly,.}}(R)
\l D^{\text{poly},\cdot}(R)
\]
We now take $\frak{l}$ to be equal to $\Omega_{R^{\text{aff}}}\cbtimes D^{\text{poly,.}}(R)$. (2) follows directly from Theorem \ref{ref-9.4.1-131}.
\end{proof}
\begin{proof}[Proof of Theorem \ref{ref-1.1-1}]
Given Theorem \ref{ref-9.5.1-136} we only need to prove that if we have a $L_\infty$
quasi-isomorphism 
\[
\psi:\Gscr\r \Hscr
\]
between sheaves of  DG-Lie algebras then $\Gscr$ and $\Hscr$ are isomorphic
in the homotopy category of DG-Lie algebras. This is done in the standard way
using the bar-cobar construction \cite{hinich1, keller5, quillen1}.

The bar-cobar construction may be performed in any symmetric abelian
monoidal category, in particular it can be done in the categories of
presheaves and sheaves of vector spaces. Since the bar-cobar
construction involves only colimits it is compatible with
sheaffification.

Considering $\psi$ first as a morphism of presheaves there is a commutative
diagram of $L_\infty$-morphisms of presheaves of DG-Lie algebras
\[
\xymatrix{
&\Omega_{\text{pre}} \mathrm{B}_{\text{pre}} \Gscr\ar[dl]_c\ar[dr]^{\phi}&\\
\Gscr \ar[rr]_\psi&&\Hscr
}
\]
where $c,\phi'$  are morphism of presheaves of DG-Lie algebras and
$c$ is a quasi-isomorphism. 

Sheaffifying we obtain an analogous diagram of sheaves of DG-Lie algebras
\[
\xymatrix{
&\Omega  \mathrm{B} \Gscr\ar[dl]_{\tilde{c}}\ar[dr]^{\tilde{\phi}}&\\
\Gscr \ar[rr]_\psi&&\Hscr
}
\]
Since sheaffication is exact $\tilde{c}$ is still a quasi-isomorphism.
Since $\psi$ is a quasi-isomorphism by assumption we obtain that
$\tilde{\phi}$ is a quasi-isomorphism as well. We conclude that
$\Gscr$ and $\Hscr$ are isomorphic in the homotopy category of sheaves
of DG-Lie algebras.
\end{proof}
\appendix
\section{Reminder on the Thom-Sullivan normalization}
\label{ref-A-137}
The material in this section is  standard. See e.g.\
\cite{hinich}. Let $k$ be a field. Let $S$ be a small category and let
$M,N:S\r \Mod(k)$ be respectively a contravariant and a covariant
functor. We define $M\underline{\otimes}_S N$ as the subset of
$\prod_{P\in \Ob(S)} M(P)\otimes_k N(P)$ consisting of $(c_P)_P$ such
that for all $\phi:P\r Q$ we have
\[
(M(\phi)\otimes 1)(c_Q)=(1\otimes N(\phi))(c_P)
\]
inside $M(P)\otimes_k N(Q)$. We  extend the bifunctor $-\underline{\otimes}-$
in the obvious way to the case where $M,N$ take values in complexes. 

We will now consider the case where is $S$ is the simplicial
category $\Delta$.  If $N(\Delta[n])$ denotes the normalized
(combinatorial) cochain complex of $\Delta[n]$ then
\[
N(-):\Delta[n]\mapsto N(\Delta[n])
\]
is a contravariant functor from $\Delta$ to $C(k)$. 

If now $A$ is a cosimplicial $k$-vector space then we may consider
\[
N(-)\,\underline{\otimes}\,{}_{\Delta} A
\]
The following is well-known
\begin{proposition} 
  \label{ref-A.1-138} $N(-)\,\underline{\otimes}\,{}_{\Delta} A$
  is canonically isomorphic to the normalized cochain complex
  $N(A)$ (given by the common kernels of the degeneracies) of $A$.
\end{proposition}
Now fix a $k$-linear DG-operad $\Oscr$. If $A$ is a cosimplicial
$\Oscr$-algebra then $N(A)$ will in general not have the
structure of an $\Oscr$-algebra.  The Thom-Sullivan  construction
repairs this defect. The idea is to replace 
the complexes $N(\Delta[n])$ in Proposition \ref{ref-A.1-138} by quasi-isomorphic complexes which
have the structure of a commutative DG-algebra.

We now assume that $k$ has characteristic zero.  Think of 
$\Delta[n]$ as the affine space
\[
\Spec k[t_0,\ldots,t_n]/(t_0+\cdots+t_n-1)
\]
Taking the algebraic De Rham complex of $\Delta[n]$ yields a
contravariant functor $\Omega^\cdot(-)$ from $\Delta$ to commutative
DG-algebras. The Thom-Sullivan normalization of a cosimplicial
$\Oscr$-algebra $A^\cdot$ is defined as
\[
N(A)^{TS}=\Omega^\cdot(-)\,\underline{\otimes}\,{}_\Delta A
\]
From the commutativity of $\Omega^\cdot(\Delta[n])$ it easily follows that
$N(A)^{TS}$ has a canonical structure as an $\Oscr$-algebra. 
\begin{proposition}\label{ref-A.2-139} \cite{BG} There is a canonical quasi-isomorphism
  (as complexes of $k$-vector spaces) 
\[
 N(A)^{\text{TS}}\r N(A)
\]
\end{proposition}
The quasi-isomorphism is constructed using functorial
homotopy equivalences $\Omega^\cdot(\Delta[n])\r N(\Delta[n])$. The latter
are 
obtained by integrating differential forms. 

If $A$ is a complex of $k$-vector spaces then we may consider $A$ as a constant
cosimplicial object.  One easily checks that
\begin{equation}
\label{ref-A.1-140}
A\cong N(A)^{\text{TS}}\cong N(A)
\end{equation}

 \section{Derived global sections of sheaves of algebras}
\label{ref-B-141}
\subsection{Introduction}
\label{ref-B.1-142}
Let $X$ be a topological space and fix a $k$-linear DG-operad $\Oscr$
for a field of characteristic zero.  If $\Ascr$ is a sheaf of
$\Oscr$-algebras on $X$ then it is easy to see that $H^\cdot(X,\Ascr)$
has the structure of an $H^\cdot(\Oscr)$-algebra.  However one would
like to give $R\Gamma(X,\Ascr)$ the structure of an $\Oscr$-algebra as
well.

 In
\cite{hinich2} Hinich constructs a model structure on the \emph{presheaves}
of $\Oscr$-algebras on $X$ which is such that a presheaf of $\Oscr$-algebras
is weakly equivalent to its sheaffification.

It follows from Hinich's construction that $R\Gamma(X,\Ascr)$ is
quasi-isomorphic to $\Gamma(X,\Ascr')$
for an arbitrary fibrant resolution $\Ascr\r \Ascr'$. In this way
we obtain indeed an actual $\Oscr$-algebra representing  $R\Gamma(X,\Ascr)$.

Note however that the choice of $\Ascr'$ is not functorial\footnote{It
is of course functorial in a homotopy theoretic sense} and furthermore
it depend on the operad $\Oscr$. In this appendix we give an
alternative construction for the algebra structure on
$R\Gamma(X,\Ascr)$ (if $\Ascr$ has left bounded cohomology) which is
functorial and whose outcome does not depend on $\Oscr$. More
precisely: for a complex of sheaves $\Ascr$ on $X$ we construct a
complex $R\Gamma(X,\Ascr)^{\tot}$ which is functorial in $\Ascr$ and
which inherits any algebra structure present on $\Ascr$.

Our construction is a generalization of a construction originally due to Hinich
and Schechtman which first replaces $\Ascr$ by a (DG-)cosimplicial algebra
\cite{hsI,hsII} using the {\v{C}ech} construction.
The Thom-Sullivan normalization (see Appendix \ref{ref-A-137}) is
then used to transform this cosimplicial algebra into a genuine
algebra over $\Oscr$. It is clear that this procedure has the properties
mentioned in the previous paragraph.

The Hinich-Schechtman construction works well for quasi-coherent
sheaves but must be modified in more general situations. This issue is
not entirely academic as non-quasi-coherent sheaves do occur in
nature. Examples in this paper are $\frak{l}$ (Theorem
\ref{ref-9.5.1-136}) and $\Omega \mathrm{B} \Gscr$ (the proof of
Theorem \ref{ref-1.1-1}).

Our initial idea was to replace \v{C}ech cohomology by a
colimit over all hypercoverings of $X$ but as the category of
hypercoverings is only filtered in a homotopy theoretic sense
\cite{SGA4}, this creates rather unpleasant technical difficulties.
Luckily it seems we can make at least some of these difficulties go
away by replacing hypercoverings with pro-hypercoverings, which is
what we will do in this section.

Although below we will work in an arbitrary Grothendieck
topos, for simplicity we will, in this introduction, continue to use
the topological space $X$.  Let $\Alg^+(X,\Oscr)$ the category of
$\Oscr$-algebra objects in $\Sh(X)$ with left bounded cohomology and
let $\Alg(\Oscr)$ be the category of $\Oscr$-algebras. We equip both
categories with weak equivalences given by quasi-isomorphisms. We will
construct a functor (see \S\ref{ref-B.8-192})
\[
\Sigma: \Alg^+(X,\Oscr)\r \Delta\Alg(\Oscr)
\]
such that  the cochain complex associated to $\Sigma (\Ascr)$  for
$\Ascr\in \Alg^+(X,\Oscr)$ is canonically
isomorphic to the derived global sections of $\Ascr$ when viewed as 
a complex of sheaves of abelian groups. This part does not require
the presence of a basefield of characteristic zero. 

\medskip

For the benefit of the reader we indicate how $\Sigma$ is defined.  We
will construct a pro-object $F=(F_\alpha)_\alpha$ in the category of
hypercovering of $X$ which is \emph{homotopy projective} (a suitable
 lifting property, see \eqref{ref-B.12-172}) and we put
\begin{equation}
\label{ref-B.1-143}
\Sigma (\Ascr)=\dirlim_\alpha \Hom(F_\alpha,\Ascr)
\end{equation}
We show in Proposition \ref{ref-B.6.3-176} that $F$ is unique up to unique
isomorphism in a homotopy theoretic sense.  This implies that $\Sigma$
is defined up to a unique natural isomorphism when viewed as a functor
between homotopy categories.

\medskip

It follows from Proposition \ref{ref-A.2-139} that if $\Oscr$ is $k$-linear for
$k$ a field of characteristic zero and we put
\[
R\Gamma(X,-)^{\text{tot}}=N(\Sigma(-))^{\text{TS}}
\]
then we obtain a functor
\[
R\Gamma(X,-)^{\text{tot}}:\Alg^+(X,\Oscr)\r \Alg(\Oscr)
\]
such that the underlying complex of vector spaces of
$R\Gamma(X,\Ascr)^{\text{tot}}$ is isomorphic to $R\Gamma(X,\Ascr)$ in
$D(k)$.

\medskip

In the last section of this appendix we outline the connection of our construction with that of Hinich in
\cite{hinich2}. 

\subsection{Simplicial objects}
In this section we recall some standard constructions on simplicial
objects.  Let $\Pscr$ be a category with arbitrary limits and
colimits. We consider the category $\Delta^\circ\Pscr$ of simplicial
objects in $\Pscr$. If $F\in \Pscr$ then we denote by $\hat{F}$ the
associated constant simplicial object.  $F\mapsto \hat{F}$ is a left
adjoint to the functor $F\mapsto F_0$.

We may define a bifunctor
\begin{equation}
\label{ref-B.2-144}
-\times-:\Delta^\circ\Pscr\times \Delta^\circ\Set\r \Delta^\circ\Pscr:(F,S)\mapsto (F_n\times S_n)_n
\end{equation}
where $F_n\times S_n$ is the $|S_n|$-fold coproduct of $F_n$. If $F\in
\Pscr$ then we define $F\times S$ as $\hat{F}\times S$.  It is easy to
see that any object $F\in \Delta^\circ \Pscr$ is a coequalizer of the
form
\begin{equation}
\label{ref-B.3-145}
\xymatrix{\coprod_{[i]\r [j]\in \Delta} F_j\times \Delta[i]\ar@<0.5em>[r] \ar@<-0.5em>[r]&
\coprod_i F_i\times \Delta[i]\ar[r] &F}
\end{equation}
The
functor $\Pscr\times \Delta^\circ\Set\r \Delta^\circ\Pscr: (F,S)\mapsto \hat{F}\times S$ 
 has a right adjoint in its second argument given by a bifunctor.
\begin{equation}
\label{ref-B.4-146}
(\Delta^\circ\Set)^\circ\times \Delta^\circ \Pscr
\r \Pscr:(S,F)\mapsto \Hom(S,F)
\end{equation}
which is the unique functor such that  $\Hom(-,F)$ sends colimits to limits and $\Hom(\Delta[n],F)=F_n$.

The  associated derived functor
\[
(\Delta^\circ\Set)^\circ\times \Delta^\circ \Pscr
\r \Delta^\circ \Pscr:(S,F)\mapsto \underline{\Hom}(S,F)
\]
defined by 
\[
\underline{\Hom}(S,F)_n=\Hom(\Delta[n]\times S,F)
\]
is the  right adjoint in the second argument to \eqref{ref-B.2-144}.

\medskip

For $F\in \Delta^\circ\Pscr$ write $F\times I=F\times \Delta[1]$ (the \emph{cylinder} object of $F$) and $F^I=
\underline{\Hom}(\Delta[1],F)$ (the \emph{path} object of $F$).  

\medskip

The category $\Delta^\circ \Pscr$ is enriched in simplicial sets (it
is a so-called simplicial category).  Let $F,G\in \Delta^\circ \Pscr$. Then the
simplicial set $\underline{\Hom}_{\Delta^\circ \Pscr}(F,G)$ is defined by
\[
\underline{\Hom}_{\Delta^\circ \Pscr}(F,G)_n=\Hom_{\Delta^\circ\Pscr}(F\times
\Delta[n],G)
\]
Define the homotopy category $\Ho(\Delta^\circ\Pscr)$ of $\Delta^\circ\Pscr$ by
\begin{equation}
\label{ref-B.5-147}
\Hom_{\Ho(\Delta^\circ\Pscr)}(F,G)=\text{connected components of
}\underline{\Hom}_{\Delta^\circ \Pscr}(F,G)
\end{equation}
In the sequel we will use the terminology exhibited in the next definition.
\begin{definitions} 
\label{ref-B.2.1-148}
\begin{enumerate}
\item Two maps $f_{0,1}:F\r G$ in $\Delta^\circ\Pscr$ are
  \emph{strictly homotopic} if there is a map $f':F\times I\r G$ such
  that the $f_i$ is the composition of
  $F=F\times\Delta[0]\xrightarrow{\partial^i} F\times\Delta[1]=
  F\times I\xrightarrow{f'} G$. 
\item Two maps $f_{0,1}:F\r G$  in $\Delta^\circ\Pscr$ are
  \emph{combinatorially homotopic}  if they can be connected 
by a chain of strict homotopies and their inverses, or equivalently
if they represent the same maps in $\Ho(\Delta^\circ\Pscr)$. 
\item A map $f:F\r G$ in $\Delta^\circ\Pscr$ is a \emph{combinatorial
    homotopy equivalence} if there is a map $g:G\r F$ such that $fg$
  and $gf$ are combinatorially homotopy equivalent to the identity, 
or equivalently if $f$ is invertible in $\Ho(\Delta^\circ\Pscr)$. 
\end{enumerate}
\end{definitions}
\begin{lemmas} Let $F\in \Delta^\circ\Pscr$. Then the functors
\[
F\times-:\Delta^\circ \Set\r \Delta^\circ \Pscr
\]
\[
\underline{\Hom}(-,F):\Delta^\circ \Set\r \Delta^\circ \Pscr
\]
preserve strict homotopy equivalent maps (and hence also combinatorially
homotopic maps and combinatorial homotopy equivalences). 
\end{lemmas}
\begin{proof} Let us consider the second functor. Let $f':S\times I\r T$
be a homotopy between maps $f_0,f_1:S\r T$ between simplicial sets. Let $\tilde{f}'$,
$\tilde{f}_0,\tilde{f}_1$ be the maps obtained applying $\underline{\Hom}(-,F)$. Since
\[
\underline{\Hom}(S\times I,F)=\underline{\Hom}(S,F)^I
\]
we obtain that $\tilde{f}'$ is a map $\underline{\Hom}(T,F)\r \underline{\Hom}(S,F)^I$ which yields a map $\underline{\Hom}(T,F)\times I\r \underline{\Hom}(S,F)$. It
is easy to see that this is a homotopy between $\tilde{f}_0$ and $\tilde{f}_1$.
\end{proof}
\begin{corollarys}\label{ref-B.2.3-149} The ``constant path'' map $F\r F^I$ is a combinatorial
homotopy equivalence.
\end{corollarys}
\begin{proof} This follows from the fact that it is obtained from the
combinatorial homotopy equivalence $\Delta[1]\r\Delta[0]$ in $\Delta^\circ \Set$. 
\end{proof}
The following is standard. 
\begin{lemmas}\label{ref-B.2.4-150}  Assume that $\Qscr$ is abelian. For $F\in
  \Delta^\circ\Qscr$ let $C_\ast(F)$ be the usual (unnormalized) chain
  complex of $F$. If $f,g:F\r G$ are strictly homotopic maps in $\Delta^\circ
  \Pscr$ then $C_\ast(f)$ and $C_\ast(g)$ are homotopic.
\end{lemmas}
The following is standard as well.
\begin{lemmas} \label{ref-B.2.5-151} Let $f_0,f_1:F\r G$ be strictly
  homotopic maps in $\Delta^\circ \Pscr$. Then the induced maps
$\ZZ f_0,\ZZ f_1:\ZZ F\r \ZZ G$ are strictly homotopic as well. 
\end{lemmas}
\begin{proof}
  $f_0$, $f_1$ are induced from a homotopy $f':F\times I\r G$. Since
  the functor $W\mapsto \ZZ W$ is a left adjoint it commutes with
  coproduct. Hence $\ZZ(F\times I)=(\ZZ F)\times I$. Thus $f'$ yields
a homotopy in $\Delta^\circ \Pscr_\ZZ$, $\ZZ f':(\ZZ F)\times I\r \ZZ G$. It is
easy to see that $\ZZ f'$ induces $\ZZ f_0$, $\ZZ f_1$. 
\end{proof}
\begin{definitions} 
\label{ref-B.2.6-152}
Let $f:F\r H$, $g:G\r H$ be in $\Delta^\circ \Pscr$. The \emph{homotopy fiber product} $F\overset{h}{\times}_H G$ is the
limit of the following diagram.
\[
\xymatrix{
F\ar[r]^f& H\\
H^I\ar[ru]_{\partial^0}\ar[rd]^{\partial_1} &\\
G\ar[r]_g& H
}
\]
\end{definitions}
If $p_{0}:F\overset{h}{\times}_H G\r F$,  $p_{1}:F\overset{h}{\times}_H G\r G$
are the resulting projection maps then clearly $f\circ p_0$ and $g\circ p_1$ are
strictly homotopic.

\begin{definitions}
Similarly if $f,g:F\r G$ are maps in $\Delta^\circ \Pscr$ then we define
the \emph{homotopy equalizer} of $f$ and $g$ as the limit of the
following diagram
\[
\xymatrix{
F\ar[r]^{f}\ar[rd]|(0.35){g} & G\\
G^I\ar[r]_{\partial_1}\ar[ru]|(0.39){\partial_0} & G
}
\]
\end{definitions}
Let $\Delta^{\le n}$ be the simplicial category truncated in dimension
$n$ and let $(-)_{\le n}$ denote the truncation functor
$\Delta^\circ\Pscr\r \Delta^{\le n,\circ}\Pscr$. The right
adjoint to $(-)_{\le n}$ is the coskeleton functor denoted by $\cosk_n$.
Concretely
\[
(\cosk_n G)_m=\Hom(\Delta[m]_{\le n},G)
\]
The truncation functor also has a left adjoint which is denoted by
$\sk_n$. If $F\in \Delta^{\le n,\circ}\Pscr$ then using the truncated
version of  \eqref{ref-B.3-145} we see that $\sk_n F$ is the coequalizer in $\Delta^\circ \Pscr$ of
\[
\xymatrix{\coprod_{[i]\r [j]\in \Delta^{\le n}} F_j\times \Delta[i]\ar@<0.5em>[r] \ar@<-0.5em>[r]&
\coprod_{i\le n} F_i\times \Delta[i]}
\]
As is customary we will also use the the notations $\sk_n$, $\cosk_n$ for the
compositions $\sk_n\circ (-)_{\le n}$,  $\cosk_n\circ (-)_{\le n}$.

\subsection{Grothendieck topoi}
\label{ref-B.3-153}
From here on
$\Pscr$ is a Grothendieck topos \cite{SGA4}.  This means that $\Pscr$
has properties very reminiscent of those of the category of sets. By
Giraud's theorem \cite{SGA4} $\Pscr$ may be realized as the category of
sheaves on a small site $\Cscr$. Therefore we sometimes we refer to
the objects of $\Pscr$ as ``sheaves''.  Recall that a \emph{site} is a
category $\Cscr$ equipped with a so-called \emph{Grothendieck
  topology}. I.e. for every $A\in \Cscr$ there is a collection
$\Tscr_A$ of subfunctors of $\Cscr(-,A)$ (called coverings) satisfying
the axioms of \cite[Def I.1.1]{SGA4}.

We recall the following standard result (see e.g.\ \cite[Prop.\ 6.20]{lowenvdb1}).
\begin{lemmas} 
\label{ref-B.3.1-154} Let $\Cscr$ be a small site and let $\Pre(\Cscr)$ and
  $\Sh(\Cscr)$ be respectively the categories of presheaves and
  sheaves on $\Cscr$. Let $a:\Pre(\Cscr)\r \Sh(\Cscr)$ be the
  sheaffification functor.  For $F\in \Pre(\Cscr)$ define
  $|F|=\sum_{C\in \Cscr} |F(C)|$. Let $|\Cscr|$ be the sum of the
  cardinalities of the $\Hom$'sets in $\Cscr$.   Then we have the following
bound
\begin{equation}
\label{ref-B.6-155}
|aF|\le |\Cscr| (2|F|)^{|\Cscr|}
\end{equation}
If $\mathfrak{a}=2^{\mathfrak{b}}$ where $\frak{b}\ge
\max(|\Cscr|,|\NN|)$ then $|F|\le \mathfrak{a}$ implies $|aF|\le
\mathfrak{a}$.
\end{lemmas}
\begin{proof} To prove \eqref{ref-B.6-155} we may assume that $F$ is separated.
Indeed if we identify sections in $F$ which are locally identical then
we only reduce $|F|$.

So assume that $F$ is separated. For any $P\in \Cscr$ we have
\[
(aF)(P)=\dirlim \Hom_{R\in \Tscr_P}(R,F)
\]
Thus
\[
|aF|\le \sum_{P\in \Cscr} |\Hom_{R\in \Tscr_P}(R,F)|
\]
Since the existence of identities implies $|\Ob(\Cscr)|\le |\Cscr|$ we
deduce from this
\[
|aF|\le |\Cscr| 2^{|\Cscr|} |F|^{|\Cscr|}
\]
which yields \eqref{ref-B.6-155}.

Thus if $|F|\le \mathfrak{a}$ with $\mathfrak{a}$ as in the statement of the
lemma.
\begin{equation}
\label{ref-B.7-156}
|aF|\le \mathfrak{a}^{|\Cscr|}=2^{\mathfrak{b}|\Cscr|}=2^{\mathfrak{b}}=\frak{a}\qed
\end{equation}
\def\qed{}\end{proof}
\begin{lemmas}
\label{ref-B.3.2-157}
Let $\Pscr$ be the category of sheaves over a small site $\Cscr$. Put
\[
\Pscr_{\mathfrak{a}}=\{F\in \Pscr\mid |F|\le \mathfrak{a}\}
\]
where $\mathfrak{a}$ is as in Lemma \ref{ref-B.3.1-154}. Then
$\Pscr_{\mathfrak{a}}$ is closed under finite limits, finite colimits, epimorphisms
and monomorphisms.

Furthermore $\Pscr_{\mathfrak{a}}$ satisfies the following cofinality
property.  For any epimorphism $f:F\r G_0$ with $G_0\in
\Pscr_{\mathfrak{a}}$ there exists a map $F_0\r F$ such that $F_0\in
\Pscr_{\mathfrak{a}}$ and the composed map $F_0\r G_0$ is an
epimorphism
\end{lemmas}
\begin{proof} This follows easily from Lemma
  \ref{ref-B.3.1-154} and the corresponding results for presheaves.
\end{proof}

Let $\Pscr_{\ZZ}$ be the category of abelian group objects in $\Pscr$. 
\begin{lemmas} \label{ref-B.3.3-158} 
\begin{enumerate}
\item $\Pscr_{\ZZ}$ is a Grothendieck category. 
\item
The forgetful functor $\Pscr_{\ZZ}\r \Pscr$ has a left adjoint.
\item If $F,G\in \Pscr_\ZZ$ then the functor of bilinear maps $\Bilin(F\times G,-)$ is representable by an object $F\otimes G$. In this way $\Pscr_\ZZ$
becomes a symmetric monoidal category. 
\end{enumerate}
\end{lemmas}
\begin{proof} 
  These facts may be proved by realizing $\Pscr$ as the category
  of sheaves on a small site $\Cscr$.   Then $\Pscr_\ZZ$
  is precisely the category of sheaves of abelian groups and the statements
are standard. 
\end{proof}
We will denote the left adjoint to $\Cscr_{\ZZ}\r \Cscr$ by $F\mapsto \ZZ F$.
If $e$ is the final object of $\Pscr$ then we write $\underline{\ZZ}$ for
$\ZZ e$.

\subsection{Hypercoverings}
We recall briefly some results about hypercoverings.
 An object $F$ in $\Delta^\circ\Pscr$ is a
 \emph{hypercovering} if for all $m$ the
 the canonical morphism 
\begin{equation}
\label{ref-B.8-159}
 F_{m+1}\r (\cosk_m F)_{m+1}
\end{equation}
is an epimorphism (see e.g.\ \cite[\S1.1]{fredneu}) and if the map of $F_0$ to 
the final object $e$ or $\Pscr$ is an epimorphism as well. Similarly $F\in
\Delta^{\le n\circ }\Pscr$ is a truncated hypercovering if
\eqref{ref-B.8-159} holds for $m\le n-1$.

\medskip

There are many equivalent characterizations for the notion of a
hypercovering.  Put $\partial\Delta[n]=\sk_{n-1} \Delta[n]$. The
following is a direct translation of the 
definition. $F$ is a hypercovering if and only if for all $n$ the morphism
\begin{equation}
\label{ref-B.9-160}
\Hom(\Delta[n],F)\r \Hom(\partial\Delta[n],F)
\end{equation}
is an epimorphism and if $F_0\r e$ is an epimorphism. This is called the \emph{local lifting property}.  From
this way of writing the definition we see that if $\Pscr$ has enough
points \cite{SGA4} then $F$ is a hypercovering if and only for every
point $p$ the simplicial set $(p^\ast F_n)_n$ is non-empty, acyclic and Kan.

\begin{remarks} The definition of hypercovering we use is in fact
a slight modification of the one used by Verdier (which depends on a site 
representing $\Pscr$).  For the original definition see \S\ref{ref-B.10-198} below.
\end{remarks}

\medskip

The following result follows from \cite[Lemma V.7.2.1]{SGA4}. 
\begin{propositions} \label{ref-B.4.2-161}  Let $F\in \Delta^\circ\Pscr$ be a hypercovering.
  Then the chain complex $C_\ast(\ZZ F)$ associated to $\ZZ F$ is
  acyclic in degrees $>0$ and its cohomology is equal to
  $\underline{\ZZ}$ in degree zero.
\end{propositions}
Note that this result is clear if $\Pscr$ has enough points since in that
case we may check it on stalks (see \cite{AM2}).

\medskip

We will frequently use the following results which are proved in the
same way as for acyclic Kan simplicial sets. In case $\Pscr$ has
enough points, they can also be checked on stalks.
\begin{propositions} 
\label{ref-B.4.3-162}
\begin{enumerate}
\item Let $F$ be a hypercovering and $S$ a finite simplicial set
  (i.e.\ $S$ has only a finite number of non-degenerate simplices).
  Then $\underline{\Hom}(S,F)$ is a hypercovering. In particular the
  path object of $F$ is a hypercovering.
\item Homotopy fiber products and homotopy equalizers of hypercoverings are
hypercoverings.
\end{enumerate}
\end{propositions}

We quote some results from \cite{AM2}.
\begin{propositions} \cite{AM2} \label{ref-B.4.4-163} Let $G$ be a hypercovering and
  let $\psi_0:F_0\r G_{\le n}$ be a morphism of hypercoverings
  trunctated in degree $n$. Then there is a hypercovering $F$ and a
  morphism of hypercoverings $\psi:F\r G$ such that $\psi_{\le n}$
is equal to $\psi_0$.
\end{propositions}

If $F\in \Delta^\circ \Pscr$ then $D_n(F)=\cup_{\sigma:[n]\r [m]\text{
    surj},m<n} \sigma F_m \subset F_n$. We call $D_n(F)$ the
\emph{degenerate part} of $F_n$. We say that $F$ is \emph{split} in degree
$n$ if $D_n(F)$ has a (necessarily unique) complement $N_n(F)$ in
$\Pscr$. $N_n(F)$ (if existing) is the \emph{non-degenerate} part of $F_n$.

If $F$ is split up to degree $n$ then one may write
\[
F_n=\coprod_{\sigma:[n]\r [m]\text{
    surj}} \sigma N_m(F)
\]
The following proposition shows that we may restrict ourselves to
split hypercoverings, if necessary.
\begin{propositions}  \cite{AM2}
\label{ref-B.4.5-164}
Assume that $G$ is a hypercovering in $\Pscr$ split up to
  degree $n$. Then there exists a map $\psi:F\r G$ where $F$ is a
  split hypercovering in $\Pscr$ and $\psi_{\le n}$ is the identity.
\end{propositions}
The next propositions shows that we may arbitrarily refine the
non-degenerate part of a split hypercovering.
\begin{propositions} \cite{AM2}
\label{ref-B.4.6-165}
Let $G$ be a split hypercovering in $\Pscr$ and let $\phi:N'\r
N_n(G)$ be an epimorphism. Then there exists a map $\psi:F\r G$ of
split hypercoverings in $\Pscr$ such that $\psi_{\le n-1}$ is the
identity and furthermore $N_n(F)=N'$ in such a way that $\psi_n$
restricts to the map $\phi$.
\end{propositions}

Throughout we fix a full small subcategory $\Pscr_0$ of $\Pscr$
which is closed under finite limits, finite  colimits, monomorphisms
and epimorphisms and which satisfies the cofinality condition of
Lemma \ref{ref-B.3.2-157}. Such a $\Pscr_0$ may be constructed
by taking $\Pscr_0$ to be a skeletal subcategory of some $\Pscr_{\mathfrak{a}}$ 
where $\Pscr_{\mathfrak{a}}$ is as in Lemma \ref{ref-B.3.2-157}.

 $\Hscr(\Pscr)$ is the category
of hypercoverings in $\Pscr$ and $\Hscr(\Pscr_0)$ is the full subcategory of
hypercoverings $F$ such that $F_n\in \Pscr_0$ for all $n$.
\begin{lemmas} \label{ref-B.4.7-166} If $G\in \Hscr(\Pscr)$ then there
  exists $F\in \Hscr(\Pscr_0)$ together with a morphism $F\r G$.
\end{lemmas}
\begin{proof} We construct $F$ step by step.  Our first step is to
  select a map $F'_0\r G$ such that the composition $F'_0\r G_0\r e$
  is an epimorphism using Lemma \ref{ref-B.3.2-157}. We then extend
  $F'_0$ to a map of  hypercoverings $F'\r G$ using Proposition
  \ref{ref-B.4.4-163}. Using Proposition \ref{ref-B.4.5-164} we may assume that
$F'$ is split. 

Assume now that we have 
constructed a map of hypercoverings $F'\r G$ such that $F'_i\in \Pscr_0$ for $i\le n$. Assume in addition that $F'$ is split. 

Consider the epimorphism $F'_{n+1}\r (\cosk_n F')_{n+1}$.  We have
$(\cosk_n F')_{n+1}\in \Pscr_0$ since the construction of the coskeleton
involves only finite limits. Let $N$ be the image of $N_{n+1}(F')$ in
$\cosk_n F'$ and choose $N_0\in \Pscr_0$ together with a map $N_0\r N_{n+1}(F')$
such that the composition $N_0\r N_{n+1}(F')\r N$ is an epimorphism. 
Put
\[
F''_{n+1}=N_0\coprod \coprod_{\sigma:[n+1]\r [m]\text{
    surj}, m\le n} \sigma N_m(F)
\]
and extend the truncated hyperovering $F''_{n+1},F_n,\ldots,F_0$ to a
hypercovering $F''$ mapping to $F'$ using Proposition \ref{ref-B.4.4-163}.
Then $F''$ coincides with $F'$ in degrees $\le n$ and is in $\Pscr_0$
in degrees $\le n+1$. Using Proposition \ref{ref-B.4.5-164} we may
assume that $F''$ is split. Repeating this procedure we ultimately
construct the desired $F$.
\end{proof}

\subsection{Pro-objects}
\label{ref-B.5-167}
Recall that if $\Dscr$ is any category then $\Pro \Dscr$ is the
category with objects denoted by formal symbols
$''\invlim''_{\alpha\in I} A_i$ where $I^\circ$ is a (small) filtered
category and $A$ is a functor $I\r \Dscr$. The $\Hom$-sets are given by
\[
\Hom_{\Pro \Dscr}(''\invlim\nolimits''_{\alpha\in I} A_\alpha,''\invlim\nolimits''_{\beta\in I} B_\beta)
=
\invlim_\beta \dirlim_\alpha \Hom_\Dscr(A_\alpha,B_\beta)
\]
Below we will omit the quotes around $\invlim$.

By \cite[\S A.4.4]{AM2} $\Pro \Dscr$ is closed under filtered inverse
limits. By \cite[Cor. 3.3]{AM2} any finite diagram $D\r \Pro \Dscr$ where
$D$ is directed (``contains no loops'') is the image 
of an object in $\Pro \Fun(D,\Dscr)$.  Informally we say that
the diagram can be constructed ``levelwise''. \cite[Prop. A.4.1]{AM2}
states that limits and colimits of finite levelwise defined limits
of pro-object can be computed levelwise as well. 

Let $L(\Dscr)$ be the category of left exact covariant functors $\Dscr\r \Set$. Then there is fully faithful embedding
\begin{equation}
\label{ref-B.10-168}
\Pro \Dscr\r L(\Dscr)^\circ:(A_\alpha)_\alpha\mapsto \dirlim_\alpha \Hom_{\Dscr}(A_\alpha,-)
\end{equation}
The construction of filtered inverse limits in $\Pro \Dscr$ in \cite[Prop.\ A.4.4]{AM2} shows that the functor \eqref{ref-B.10-168} commutes with filtered
inverse limits.  In particular the
objects in $\Dscr$ are ``cofinitely presented''. Let $(F_i)_{i\in I}$ be
a filtered inverse system of objects in $\Pro \Dscr$ and $F\in \Dscr$. Then
\begin{equation}
\label{ref-B.11-169}
\Hom_{\Pro \Dscr}(\invlim_i F_i,F)=\dirlim_i \Hom_{\Pro \Dscr}(F_i,F)
\end{equation}
Below we will work in (full) subcategories of $\Pro \Delta^\circ
\Pscr$.  It is clear that $\Pro \Delta^\circ \Pscr$ is a simplicial
category (it may be enriched in simplicial sets). 
The functors $-\times S$ and
$\underline{\Hom}(S,-)$ for $S\in \Delta^\circ \Set$ may be
extended to $\Pro \Delta^\circ \Pscr$ and they remain adjoints. 
In particular cylinder and path objects exist in $\Pro \Delta^\circ
\Pscr$ and we may define homotopy fiber products and equalizers in 
$\Pro \Delta^\circ
\Pscr$.

It also clear that the Definition  \ref{ref-B.2.1-148}
make sense in this context and furthermore we can define $\Ho\Pro \Delta^\circ\Pscr$ using the
formula \eqref{ref-B.5-147}.
\subsection{Pro-hypercoverings}
\label{ref-B.6-170}
We consider the full subcategory $\Pro \Hscr(\Pscr)$ of $\Pro
\Delta^\circ \Pscr$. We refer to the objects in $\Pro \Hscr(\Pscr)$ as 
pro-hypercoverings. 

We note the following generalization of Proposition \ref{ref-B.4.3-162}.
\begin{propositions} 
\label{ref-B.6.1-171}
\begin{enumerate}
\item Let $F$ be a pro-hypercovering and $S$ a finite simplicial set
  (i.e.\ $S$ has only a finite number of non-degenerate simplices).
  Then $\underline{\Hom}(S,F)$ is a pro-hypercovering. In particular the
  path object of $F$ is a pro-hypercovering.
\item Homotopy fiber products and homotopy equalizers of pro-hypercoverings are
pro-hypercoverings.
\end{enumerate}
\end{propositions}
\begin{proof} (1) follows directly from Proposition \ref{ref-B.4.3-162}(1)
and (2) follows from Proposition \ref{ref-B.4.3-162}(1) and the fact that the diagrams for computing homotopy
fiber products and equalizers can be constructed levelwise (see \S\ref{ref-B.5-167}). 
\end{proof}
We say that $F\in \Pro
\Hscr(\Pscr)$ is \emph{homotopy projective} (with respect to $\Hscr(\Pscr)$) if
every diagram of solid arrows
\begin{equation}
\label{ref-B.12-172}
\xymatrix{
& F\ar@{.>}[dl]\ar[d]\\
C'\ar[r]&C 
}
\end{equation}
in $ \Pro
\Hscr(\Pscr)$ with $C,C'\in \Hscr(\Pscr)$ can be completed with 
a dotted arrow in $\Pro \Hscr(\Pscr)$ such that the resulting diagram is commutative in $\Ho\Pro \Hscr(\Pscr)$. Let us denote the category of 
projective pro-hypercoverings by $\Proj \Hscr(\Pscr)$.

The following is our main technical result.
\begin{propositions} \label{ref-B.6.2-173} For every pro-hypercovering $F$ there exists
a map of pro-hypercoverings $F^\ast\r F$ such that $F^\ast$ is homotopy projective. 
\end{propositions}
\begin{proof} The proof is adapted from the proof of \cite[Thm 2.7]{Barr}.
Let $\Pscr_0\subset \Pscr$ be as in \S\ref{ref-B.3-153} but choose $\Pscr_0$
large enough such that $F\in \Pro \Hscr(\Pscr_0)$ (up to isomorphism). This
may be done by choosing the cardinal $\mathfrak{a}$ in Lemma \ref{ref-B.3.2-157}
large enough.

We will temporarily work in $\Pro \Hscr(\Pscr_0)$. 
We start by well-ordering the diagrams in $\Pro \Hscr(\Pscr_0)$
\[
\xymatrix{
& F\ar[d]\\
C'\ar[r] & C
}
\]
with $C$, $C'\in \Hscr(\Pscr_0)$.

We construct an ordinal sequence 
\begin{equation}
\label{ref-B.13-174}
\cdots \r F_\omega\r\cdots \r F_1\r F_0
\end{equation}
in $\Pro\Hscr(\Pscr_0)$ as follows: $F_0=F$; at a limit ordinal $\lambda$
let $F_\lambda=\invlim_{\mu< \lambda} F_{\mu}$. To define $F_\lambda$ for
a successor cardinal $\lambda=\mu+1$ let $C'\r C\l F$ be the least
diagram (if existing) for the well ordering such that the diagram of solid arrows
\[
\xymatrix{
F_\mu\ar[r]\ar@{.>}[d]& F\ar[d]\\
C'\ar[r] & C
}
\]
cannot be completed with the dotted arrow.

Put $F_\lambda=C'\overset{h}{\times}_{C}
F_\mu$. 
Then the resulting diagram 
\[
\xymatrix{
F_{\lambda}\ar[r]\ar[dr]&F_\mu\ar[r]& F\ar[d]\\
&C'\ar[r] & C
}
\]
is commutative in $\Pro \Hscr(\Pscr_0)$ up to a strict homotopy. 

Since $\Hscr(\Pscr_0)$ is small it follows that this procedure
has to stop for some ordinal $\sigma$. Put $F^\sharp=F_\sigma$ 
Then it follows that any diagram of solid arrows
\[
\xymatrix{
F^\sharp\ar[r]\ar@{.>}[d]& F\ar[d]\\
C'\ar[r] & C
}
\]
can be completed with the dotted arrow up to a strict homotopy.

\medskip

Now define a sequence 
\[
\cdots \r F'_1\r F'_0=F^\sharp
\]
where $F'_{n+}=(F'_n)^\sharp$ and put $F^\ast=\invlim_n F'_n$. We claim that any 
diagram of solid arrows 
\begin{equation}
\label{ref-B.14-175}
\xymatrix{
& F^\ast\ar@{.>}[dl]\ar[d]\\
C'\ar[r]&C
}
\end{equation}
with $C$, $C'\in \Hscr(\Pscr_0)$ can be completed with the dotted arrow up to a strict homotopy. Indeed by \eqref{ref-B.11-169}
the vertical map is obtained from some map $F'_n\r C$. But then by construction
we may factor $F'_{n+1}$ through $C'$. 

\medskip

We now claim that $F^\ast$ is homotopy projective. So we consider a
solid diagram as in \eqref{ref-B.14-175} but now we only require $C$,
$C'\in \Hscr(\Pscr)$.  We have $F^\ast=\invlim
(F^\ast_\alpha)_{\alpha\in A}$ with $F^\ast_\alpha\in \Hscr(\Pscr_0)$.
So the vertical map in \eqref{ref-B.14-175} is obtained from some map 
$F^\ast_\alpha\r C$.  Put $D=F^\ast_\alpha$. By Lemma \ref{ref-B.4.7-166}
we can construct a map of hypercoverings $D'\r C'\overset{h}{\times}_C D$
with $D'\in \Hscr(\Pscr_0)$. We may then construct a diagram in $\Pro \Hscr(\Pscr)$ 
\[
\xymatrix{
& F^\ast\ar[dl]\ar[d]\\
D'\ar[r]\ar[d]&D\ar[d]\\
C'\ar[r]&C
}
\]
which is commutative in $\Ho\Pro \Hscr(\Pscr)$. This finishes the proof.
\end{proof}
Let $\Wscr$ be the set of maps in $\Ho\Pro \Hscr(\Pscr)$ between homotopy
projectives. We have the following result:
\begin{propositions} 
\label{ref-B.6.3-176}
For any
$F,G\in \Pro \Hscr(\Pscr)$ with $F$ homotopy projective there
is precisely one map $F\r G$ in $\Wscr^{-1}\Ho\Pro \Hscr(\Pscr)$.
\end{propositions}

\begin{proof}
\begin{step} If $f,g:F\r G$ are maps in $\Pro \Hscr(\Pscr)$ and 
$F$ is homotopy projective then the images of $f$ and $g$ are the same in
  $\Wscr^{-1}\Ho \Pro \Hscr(\Pscr)$.

\medskip

To see this let $K'\r F$ be the homotopy equalizer of $f$ and $g$ and let
$K\r K'$ be a homotopy projective object mapping to $K$ (constructed
using \ref{ref-B.6.2-173}). If $k:K\r F$ is the composed map then $fk$ and
$gk$ are the same in $\Ho \Pro \Hscr(\Pscr)$ . Since $k$ is 
in $ \Wscr$  this implies that $f$ and $g$ are the same in $\Wscr^{-1}\Ho \Pro \Hscr(\Pscr)$
\end{step}
\begin{step} Any map $f:F \r G$ in $\Wscr^{-1} \Ho\Pro \Hscr(\Pscr)$ with $F$
  homotopy projective can be written as $vu^{-1}$ where $u,v$
fit in a diagram in $\Pro \Hscr(\Pscr)$ 
\[
\xymatrix{
& U \ar[dr]^v\ar[dl]_u &\\
F & & G
}
\]
with $U$ homotopy projective. 

\medskip

It is easy to see that it is sufficient to prove that if $f$ is of the indicated form then so is $w^{-1}f$ with $w:H\r G$ in $\Wscr$. To see
this we make the following commutative diagram in $\Ho \Pro \Hscr(\Pscr)$.
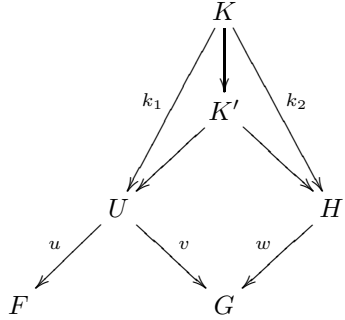
\[
\xymatrix{
&&K\ar[d]\ar[ddl]_{k_1}\ar[ddr]^{k_2}&\\
&& K'\ar[dl]\ar[dr]&\\
& U \ar[dr]^v\ar[dl]_u &&H\ar[dl]_w \\
F & & G
}
\] 
where $K'=U\overset{h}{\times}_G H$ and $K$ is homotopy projective. 
Then in $\Ho \Pro \Hscr(\Pscr)$ we have $wk_2=vk_1$ and thus in $\Wscr^{-1} \Ho\Pro \Hscr(\Pscr)$: $w^{-1}v=k_2 k_1^{-1} $. Hence $w^{-1}f=w^{-1} v u^{-1}=
k_2(uk_1)^{-1}$. 
\end{step}
\begin{step} 
If $F, G\in  \Pro \Hscr(\Pscr)$ with $F$ homotopy projective
then there is at most one map $F\r G$ in $\Wscr^{-1}\Ho\Pro \Hscr(\Pscr)$.

\medskip

Assume that there are two maps $vu^{-1}$, $v' u^{\prime-1}$ with
``middle objects'' $U$ and $U'$ as in Step 2. Let $U''$ be a homotopy
projective mapping to $U\times U'$ (using Proposition
\ref{ref-B.6.2-173}) . Using Step 1 we have a commutative diagram in
$\Wscr^{-1} \Ho\Pro \Hscr(\Pscr)$
 
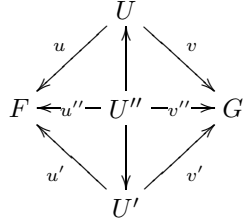
\[
\xymatrix{
& U\ar[dr]^{v}\ar[dl]_{u} &\\
F & U''\ar[r]|{v''}\ar[l]|{u''}\ar[d]|{}\ar[u]|{}& G\\
& U'\ar[ur]_{v'}\ar[ul]^{u'} &
}
\]
from which we obtain $vu^{-1}=v'' u^{\prime\prime -1}=v'u^{\prime -1}$. 
\end{step}
\begin{step}If $F, G\in  \Pro \Hscr(\Pscr)$ with $F$ homotopy projective
then there is precisely one map $F\r G$ in $\Wscr^{-1}\Ho\Pro \Hscr(\Pscr)$.
\medskip

By Step 3 we only have to show that there is a map $F\r G$. Let $K$ be
a homotopy projective mapping to $F\times G$. Denote the maps of $K$
to $F$ and $G$ by $u$ and $v$. Then $vu^{-1}$ is the required map. \qed
\end{step}
\def\qed{}\end{proof}
Let $\Proj \Ho\Pro \Hscr(\Pscr)$ be the full subcategory of $\Ho\Pro \Hscr(\Pscr)$
consisting of homotopy projective objects. The same proof as the
previous proposition, replacing $ \Ho\Pro \Hscr(\Pscr)$  by $\Proj \Ho\Pro \Hscr(\Pscr)$ yields the following result.
\begin{corollarys}
\label{ref-B.6.4-177}
The category $\Wscr^{-1}\Proj \Ho \Pro\Hscr(\Pscr)$ is
equivalent to the singleton category. I.e.\ the category with one object
and one (identity) arrow. 
\end{corollarys}
\subsection{Complexes of sheaves of abelian
groups}
\label{ref-B.7-178}
For an abelian category $\Ascr$ let  $C(\Ascr)$ be the category
of cochain complexes over $\Ascr$.

Let us say that contravariant functor $H:\Pscr\r \Ab$ is \emph{weakly
effaceable} if for every $G\in \Pscr$ and for every $h\in H(G)$ there
exists an epimorphism $\phi:F\r G$ in $\Pscr$ such that $H(\phi)(h)=0$.

 We
say that a contravariant functor $H:\Hscr(\Pscr)\r \Ab$ is
\emph{weakly effaceable} if for every $G\in \Hscr(\Pscr)$ and for every
$h\in H(G)$ there exists a map of hypercoverings $\phi:F\r G$  such
that $H(\phi)(h)=0$. 

We will need the following result. 
\begin{lemmas} Let $H:\Pscr\r \Ab$ be a weakly effaceable functor which sends
finite coproducts to products and let
 $G$ be a hypercovering. Let $m\in \ZZ$ and let $a\in H(G_m)$.
Then there exists a map of hypercoverings $\psi:F\r G$ such that $H(\psi_m)(a)=0$ in $H(F_m)$.
\end{lemmas}
\begin{proof} By Proposition \ref{ref-B.4.5-164} we may assume that $G$
  is split. I.e.
\[
G_m=\coprod_{\sigma:[m]\r [p]\text{
    surj}} \sigma N_p(G)
\]
and hence $a=\sum_{\sigma:[m]\r [p]\text{
    surj}} \sigma a_\sigma$ where $a_\sigma\in H^p(N_p(G))$.

Let $N'\r N_p(G)$ is an epimorphism which effaces $a_\sigma$. Using Proposition
\ref{ref-B.4.6-165} we may refine $G$ to a split hypercovering $G'$
whose non-degenerate part is $N'$ in degree $p$ and which is unchanged
in lower degrees.

Starting with the maximal $p$ such that $a_p\neq 0$ and work our way down
we eventually find a hypercovering in which the image of all $a_p$ is zero.
\end{proof}
\begin{corollarys} 
\label{ref-B.7.2-179} If $H:\Pscr\r \Ab$ is weakly effaceable and sends finite coproducts to products then for
  all $m$ the functor $\Hscr(\Pscr)\r \Ab: G\mapsto H(G_m)$ is
  effaceable as well.
\end{corollarys}

\begin{lemmas} 
\label{ref-B.7.3-180}
For all acyclic
complexes $\Lscr\in C(\Pscr_\ZZ)$  the functor
$H^0(\Hom(-,\Lscr)):\Pscr\r \Ab$ is weakly effaceable.
\end{lemmas}
\begin{proof}
  Let $a\in H^0(\Hom(G,\Lscr))$ with $G\in\Pscr$. Thus $a$ is
  represented by a map $G\r \ker(\Lscr^0\r \Lscr^{1})= \im
  (\Lscr^{-1}\r \Lscr^0)$. Let $F$ be the pullback of the diagram
\[
\xymatrix{
 & \Lscr^{-1}\ar[d]\\
G\ar[r]&\Lscr^{0}
}
\]
Then $F\r G$ is an epimorphism and the image of $a$ in
$H^0(\Hom(F,\Lscr))$ is zero. 
\end{proof}
If $A$ is a cosimplicial abelian group then as usual we denote by
$C^\ast(A)$ the (unnormalized) cochain complex associated to $A$. If
$A$ is a cosimplicial object in the category of  complexes of abelian groups then by $C^\ast(A)$ we will denote
the total (product) complex of the double complex with rows
$C^\ast(A^n)$.

If $F\in \Delta^\circ \Pscr$ and $\Lscr$ is a complex in $\Pscr_{\ZZ}$ then by $\Hom(F,\Lscr)$ we denote the
cosimplicial object in the category of complexes of abelian groups defined by
\[
\Hom(F,\Lscr)^n=\Hom(F_n,\Lscr)
\]
 The following formula is clear
\begin{equation}
\label{ref-B.15-181}
\underline{\Hom}(C_\ast(\ZZ F),\Lscr)=C^\ast(\Hom(F,\Lscr))
\end{equation}
where the left hand side is the usual differentially graded Hom of complexes.

For $\Lscr\in \Pscr_\ZZ$ we define $H_\Lscr:\Hscr(\Pscr)\r \Ab$ as the
functor $H^0(C^\ast(\Hom(-,\Lscr)))$.
\begin{lemmas} 
\label{ref-B.7.4-182} If $\Lscr\in \Pscr_\ZZ$ is acyclic then the functor 
$H_\Lscr$ is weakly effaceable.
\end{lemmas}
\begin{proof} Assume $a\in H_\Lscr(G)$ is represented by a morphism
\[
a:C_\ast(\ZZ G)\r \Lscr
\]
Let $N_\ast(\ZZ G)$ be the normalized chain complex of $G$ (i.e.\ the quotient
of $C_\ast(\ZZ G)$ by the images of the degeneracies). It follows
from the proof of \cite[Thm 8.3.8]{Weibel} that the canonical map $C_\ast(\ZZ G)
\r N_\ast(\ZZ G)$ is a homotopy equivalence. Indeed the proof shows that the
kernel $D$ of this map is of the form $\bigcup_p D_p$ such that $D_{p+1}/D_p$ is
contractible. Therefore $D$ is itself contractible which is sufficient.

Hence up to homotopy we may view $a$ as a map
\[
b:N_\ast(\ZZ G)\r \Lscr
\]
Without loss of generality we may assume that $G$ is split. Then
$N_\ast(\ZZ G)_n=\ZZ N_n(G)$.

We must construct a map of hypercoverings $\phi:F\r G$ and a homotopy
$h:N_\ast(\ZZ F)\r \Lscr[-1]$ such that $b\circ \phi=dh+hd$.

We will construct $F$ and $h$ step by step.  Suppose we have constructed
a morphism of split hypercoverings $\phi':F'\r G$ and maps $h'_i:\ZZ N_i(F')
\r \Lscr_{i+1}$ for $i< n$ such that $b'_i=dh'_i+h'_{i-1}d$ for
$i=0,\ldots,n-1$ where $b'=b\circ \phi'$ and $h'_{-1}=0$.

Put $c=b'_n-h'_{n-1}d$. Then $dc=0$. 
Thus $c$ defines an element $\bar{c}$ of 
$H^0(\Hom(N_n(F'),\Lscr[-n]))$. By Lemma \ref{ref-B.7.3-180} there exists an epimorphism
$f:N'\r N_n(F')$ in $\Pscr$ which effaces $\bar{c}$. 

By Proposition \ref{ref-B.4.6-165} there is a map of split hypercoverings
$\psi:F''\r F'$ such that $\psi_{\le n-1}$ is the identity and furthermore 
$N_n(F'')=N'$ in such a way that $\psi_n$
restricts to the map $f$. 

Put $h''_i=h'_i\circ \psi$, $b''_i=b'_i\circ \psi$.  Then still 
$b''_i=dh''_i+h''_{i-1}d$ for $i<n$ but now $b''_n-h''_{n-1}d$ is of the form $dh''_n$
for some map $h''_n:\ZZ N_n(F'')\r \Lscr_{n+1}$. Repeating this procedure we ultimately
construct the desired $F$ and $h$. 
\end{proof}
Below a contravariant functor $H:\Hscr(\Pscr)\r \Ab$ will be extended implicitly to
a
contravariant functor $\Pro \Hscr(\Pscr)\r \Ab$ by putting
\[
H(\invlim_{\alpha\in I} F_\alpha)=\dirlim_{\alpha\in I} H(F_\alpha)
\]
Let us say that a contravariant functor $H:\Hscr(\Pscr)\r \Ab$ is
\emph{homotopy insensitive} if it factors through $\Ho \Hscr(\Pscr)$.
This is equivalent with demanding that $H$ inverts constant path
maps. Since this condition lifts to pro-objects it
follows in particular that $H$ extends to a functor $\Ho \Pro \Hscr(\Pscr)\r
\Ab$. 
\begin{lemmas} \label{ref-B.7.5-183} $H_\Lscr(-)$ is homotopy insensitive for any $\Lscr\in
   C(\Pscr_\ZZ)$.
\end{lemmas}
\begin{proof} If $F\in \Hscr(\Pscr)$ then according to Corollary
  \ref{ref-B.2.3-149}, the constant path map $F\r F^I$ is a combinatorial
  homotopy equivalence. It follows from Lemma \ref{ref-B.2.5-151} 
that $\ZZ F\r \ZZ (F^I)$ is a combinatorial homotopy equivalence in $\Delta^\circ
\Pscr_\ZZ$. 

Hence by Lemma \ref{ref-B.2.4-150} the induced map $C_\ast(\ZZ F)\r C_\ast(\ZZ F^I)$ is a homotopy
equivalence. It follows that $\underline{\Hom}(
C_\ast(\ZZ F^I),\Lscr)\r \underline{\Hom}(
C_\ast(\ZZ F),\Lscr)$. is a homotopy equivalence. Then formula \eqref{ref-B.15-181} finishes the proof. 
\end{proof}
\begin{lemmas} \label{ref-B.7.6-184} Let $H:\Hscr(\Pscr)\r \Ab$ be a homotopy insensitive
  weakly effaceable functor. Then for any homotopy projective $F\in
  \Ho\Pro \Hscr(\Pscr)$ we have $H(F)=0$.
\end{lemmas}
\begin{proof} We have $F=\invlim_{\alpha\in A} F_\alpha$ with
  $F_\alpha \in \Hscr(\Pscr)$ and $H(F)=\dirlim_\alpha
  H(F_\alpha)$. Let $h\in H(F)$.
  Then $h$ is represented by some $h_\alpha \in
  H(F_\alpha)$. Since $H$ is weakly effaceable
   there exists a map of hypercoverings $F'\r
  F_\alpha$ such that the  image
  of $h_\alpha$ in $H(F')$ is zero. Since $F$ is homotopy projective
the map $F\r F_\alpha$ factors through $F'$ in $\Ho\Pro\Hscr(\Pscr)$. This
implies that $h$ is zero.
\end{proof}
%

If $F$ is the pro-object $(F_\alpha)_\alpha$ then  we define
\[
C^\ast(\Hom(F,\Lscr))_{\text{pro}}=\dirlim_\alpha C^\ast(\Hom(F_\alpha,\Lscr))
\]
and 
\[
  C^\ast(\Hom(F,\Lscr))=C^\ast(\dirlim_\alpha
  \Hom(F_\alpha,\Lscr))
\]
We show below that $C^\ast(\Hom(F,\Lscr))_{\text{pro}}$ is well-behaved
in its first argument and $C^\ast(\Hom(F,\Lscr))$ is well behaved in
its second argument. Furthermore there is an obvious map
\[
C^\ast(\Hom(F,\Lscr))_{\text{pro}}\r C^\ast(\Hom(F,\Lscr))
\]
which is an isomorphism if $\Lscr$ has
left bounded cohomology (see Lemma \ref{ref-B.7.9-189} below).

\begin{propositions} \label{ref-B.7.7-185} If $f:F\r G$ is a map between homotopy projective 
pro-hypercoverings then $H_\Lscr(f)$ is invertible for any $\Lscr\in
C(\Pscr_\ZZ)$.
\end{propositions}
\begin{proof}
Note first that
\begin{equation}
\label{ref-B.16-186}
H_\Lscr(F)= H^0( C^\ast(\Hom(F,\Lscr)_{\pro})
\end{equation}
As $\Pscr_\ZZ$ is a Grothendieck category
there is a quasi-isomorphism $q:\Lscr\r\Escr$ where $\Escr$ is homotopy
injective \cite{TLS}. I.e.\ $\underline{\Hom}_{\Pscr_\ZZ}(\Qscr,\Escr)$
is acyclic for every acyclic $\Qscr$. 

Using Proposition \ref{ref-B.4.2-161} and 
formula \eqref{ref-B.15-181} we 
obtain that for any hypercovering $E$ the canonical map
\[
C^\ast(\Hom(E,\Escr))\r \Hom_{\Pscr_\ZZ}(\underline{\ZZ},\Escr)
\]
is a quasi-isomorphism.  Taking direct limits we obtain that for any
pro-hypercovering $E$ we have a canonical quasi-isomorphism
\begin{equation}
\label{ref-B.17-187}
C^\ast(\Hom(E,\Escr))_{\pro}\r \Hom_{\Pscr_\ZZ}(\underline{\ZZ},\Escr)
\end{equation}
 Let $\Cscr$ be the cone of $q$. Then
$\Cscr$ is acyclic. 
We obtain a morphism of distinguished triangles in $K(\Ab)$
(the homotopy category of $\Ab$):
\[
\xymatrix{
C^\ast(\Hom(F,\Lscr))_{\pro}\ar[r] &C^\ast(\Hom(F,\Escr))_{\pro}\ar[r]& C^\ast(\Hom(F,\Cscr))_{\pro}
\ar[r]&\\
C^\ast(\Hom(G,\Lscr))_{\pro}\ar[r]\ar[u] &C^\ast(\Hom(G,\Escr))_{\pro}\ar[r]\ar[u]& C^\ast(\Hom(G,\Cscr))_{\pro}
\ar[r]\ar[u]&
}
\]
By Lemmas \ref{ref-B.7.4-182},\ref{ref-B.7.5-183},\ref{ref-B.7.6-184} and
\eqref{ref-B.16-186} $C^\ast(\Hom(F,\Cscr))_{\pro}$ and
$C^\ast(\Hom(G,\Cscr))_{\pro}$ are acyclic.  By \eqref{ref-B.17-187} the
middle vertical map is a quasi-isomorphism. Hence it follows that the
left most vertical map is a quasi-isomorphism as well.
\end{proof}
\begin{lemmas} 
  \label{ref-B.7.8-188}  
  Assume that $F\in \Pro \Hscr(\Pscr_\ZZ))$ is homotopy projective. If
  $\Lscr\in C(\Pscr_\ZZ)$ is acyclic then we have that $C^\ast(\Hom(F,\Lscr))$ is
  acyclic. 
\end{lemmas}
\begin{proof}
Let the index of $\Lscr$ be denoted by $q\in \ZZ$ and let  the
index of an object in $\Hscr(\Pscr)$ be denoted by $p\in \NN$. Let
$H$ be the functor $\Hscr(\Pscr)\r \Ab$ given by
\[
H(G)=\bigoplus_{p,q} H^p(H^q(\Hom(G,\Lscr)))
\]
By Lemma \ref{ref-B.7.3-180} and Corollary \ref{ref-B.7.2-179} we see that
$U^q=H^q(\Hom(-,\Lscr))$ is weakly effaceable. Thus $U=\oplus_q U^q$
is weakly effaceable as well. By an argument similar to Lemma
\ref{ref-B.7.5-183} we deduce that $H^p(U)$ is homotopy insensitive.  Thus we conclude that $H$ is both weakly effaceable
and homotopy insensitive. By Lemma \ref{ref-B.7.6-184} we conclude $H(F)=0$
for $F\in \Hscr(\Pscr)$ and hence for $F\in \Pro\Hscr(\Pscr)$.

Hence the $E_2$ term of the spectral sequence computing the cohomology
of $C^\ast(\Hom(F,\Lscr))$ vanishes. Either by invoking the correct convergence
criterion or by a direct diagram chase (which the author did) this implies
that $C^\ast(\Hom(F,\Lscr))$ is acyclic. 
\end{proof}
Let $C^+(\Pscr_\ZZ)$ be the full subcategory of $C(\Pscr_\ZZ)$ consisting 
of complexes with left bounded cohomology. 
\begin{lemmas} 
\label{ref-B.7.9-189}
Assume that $\Lscr\in C^+(\Pscr_\ZZ)$. Then the canonical map
\[
C^\ast(\Hom(F,\Lscr))_{\pro}\r C^\ast(\Hom(F,\Lscr))
\]
is a quasi-isomorphism.
\end{lemmas}
\begin{proof} 
 Choose a quasi-isomorphism $\Lscr\r \Escr$ to a left bounded complex
  of injectives $\Escr$. Let $\Cscr$ be the cone. Then we have the
  following morphism of distinguished triangles in $K(\Ab)$
\[
\xymatrix{
C^\ast(\Hom(F,\Lscr))_{\pro}\ar[r]\ar[d]&C^\ast(\Hom(F,\Escr))_{\pro}\ar[r]\ar[d]&C^\ast(\Hom(F,\Cscr))_{\pro}\ar[r]\ar[d]&\\
C^\ast(\Hom(F,\Lscr))\ar[r]&C^\ast(\Hom(F,\Escr))\ar[r]&C^\ast(\Hom(F,\Cscr))\ar[r]&
}
\]
By Lemmas \ref{ref-B.7.5-183},\ref{ref-B.7.4-182},\ref{ref-B.7.6-184} and Lemma
\ref{ref-B.7.8-188} $C^\ast(\Hom(F,\Cscr))_{\pro}$ and  $C^\ast(\Hom(F,\Cscr))$
are acyclic. Furthermore since $\Escr$ is left bounded it is easy to
see that the middle map is an isomorphism. Hence the left most map 
is a quasi-isomorphism.
\end{proof}
 For a homotopy projective $F$ in $\Pro \Hscr(\Pscr)$ let
$\Pi_F$ be the functor
\[
\Pi_F:C(\Pscr_\ZZ)\r \Delta C(\Ab):\Lscr\mapsto \Hom(F,\Lscr)
\]
\begin{lemmas} 
  \label{ref-B.7.10-190} The functor $C^\ast\circ \Pi_F$
 sends weak
  equivalences to quasi-isomorphisms.
\end{lemmas}
\begin{proof} By considering the cones of quasi-isomorphisms, it is sufficient
to prove that for any acyclic $\Lscr\in C(\Pscr_\ZZ)$ we have that
$C^\ast(\Hom(F,\Lscr))$ is acyclic. This is precisely Lemma \ref{ref-B.7.8-188}.
\end{proof}
The following proposition is the raison d'\^etre for the functor $\Pi_F$.
\begin{lemmas}\label{ref-B.7.11-191} When restricted to $C^+(\Pscr_\ZZ)$ the composition $C^\ast\circ \Pi_F$ is canonically isomorphic
  to $\RHom_{\Pscr_\ZZ}(\underline{\ZZ},-)$.
\end{lemmas}
\begin{proof} If $\Escr$ is a left bounded complex of injectives then
\[
C^\ast(\Hom(F,\Escr))=C^\ast(\Hom(F,\Escr))_{\pro}
\]
and the latter is equal to $\Hom_{\Pscr_\ZZ}(\underline{\ZZ},\Escr)$ by 
\eqref{ref-B.17-187}
\end{proof}
\subsection{Sheaves of algebras}
\label{ref-B.8-192}
In addition to the above notations, in this section $\Oscr(n)_n$ will
be a fixed DG-operad of abelian groups. We write $\Alg(\Oscr)$ for the
category of $\Oscr$-algebras.

Since by Lemma \ref{ref-B.3.3-158} $\Pscr_\ZZ$
is a symmetric monoidal category we may speak of $\Oscr$-algebra
objects on $\Pscr$. We define  $\Alg(\Pscr,\Oscr)$ as the
category of $\Oscr$-algebras in $\Pscr$.

We make the following definitions.
\begin{enumerate}
\item A weak equivalence in $\Alg(\Pscr,\Oscr)$ is a quasi-isomorphism.
\item A map $A\r B$ in $\Delta \Alg(\Oscr)$ is a weak equivalence
if $C^\ast(A)\r C^\ast(B)$ is a quasi-isomorphism. 
\end{enumerate}
Following custom the classes of weak equivalences will be denoted by
$\Wscr$. Note that if $F\in \Pscr$ and $\Ascr\in \Alg(\Pscr,\Oscr)$
then by construction $\Hom(F,\Ascr)\in \Delta\Alg(\Oscr)$.  Let $\Pi_\ast$
be the bifunctor
\[
\Pi_\ast:\Pro \Hscr(\Pscr)\times \Alg(\Pscr,\Oscr)\r \Delta
\Alg(\Oscr):\Ascr\mapsto \Hom(F,\Ascr)
\]
Let $\Wscr^{\text{cp}}$ be the constant path maps in $\Pro \Hscr(\Pscr)$.
  Let $\Alg^+(\Pscr,\Oscr)$ be the full subcategory of $\Alg(\Pscr,\Oscr)$
whose objects have left bounded cohomology. 

According to Lemma \ref{ref-B.7.5-183} and Lemma \ref{ref-B.7.10-190} we obtain 
a bifunctor
\[
\Pi_\ast:\Wscr^{\text{cp},-1} \Pro \Hscr(\Pscr)\times \Alg(\Pscr,\Oscr)\r \Wscr^{-1}\Delta
\Alg(\Oscr):\Ascr\mapsto \Hom(F,\Ascr)
\]
and hence a bifunctor
\[
\Pi_\ast:\Ho \Pro \Hscr(\Pscr)\times \Alg(\Pscr,\Oscr)\r \Wscr^{-1}\Delta
\Alg(\Oscr):\Ascr\mapsto \Hom(F,\Ascr)
\]
$\Pi_\ast$ restricts
to a bifunctor
\[
\Pi_\ast:\Proj \Ho \Pro \Hscr(\Pscr)\times \Alg(\Pscr,\Oscr)\r \Wscr^{-1}\Delta
\Alg(\Oscr):\Ascr\mapsto \Hom(F,\Ascr)
\]
Using Proposition \ref{ref-B.7.7-185} we obtain a bifunctor
\[
\Pi_\ast:\Wscr^{-1}\Proj \Ho \Pro \Hscr(\Pscr)\times \Alg(\Pscr,\Oscr)\r \Wscr^{-1}\Delta
\Alg(\Oscr):\Ascr\mapsto \Hom(F,\Ascr)
\]
By Corollary \ref{ref-B.6.4-177} the first argument of $\Pi_\ast$ is now a singleton
category. 

Below we define $\Sigma=\Pi_\ast(F,-)$  for an arbitrary pro-hypercovering
$F$. It follows from the above discussion that
 $\Sigma$ is well defined up
to a unique natural isomorphism. It follows from Lemma \ref{ref-B.7.11-191} that
the following diagram is commutative
\[
\begin{CD} 
\Alg^+(\Pscr,\Oscr) @>\Sigma>> \Delta \Alg(\Oscr)\\
@VVV @VV C^\ast V\\
C(\Pscr_\ZZ)@>>\RHom_{\Pscr_\ZZ}(\underline{\ZZ},-)> D(\Ab)
\end{CD}
\]
where the left arrow is the forgetful functor.

\medskip

Let $F\in \Pro\Hscr(\Pscr)$ and $\Ascr\in \Alg^+(\Oscr)$. Choose an
arbitrary projective pro-hypercovering $P$. According to Proposition
\ref{ref-B.6.3-176} there is a unique map $P\r F$ in $\Wscr^{-1}\Pro
\Hscr(\Pscr)$. In this way we obtain a canonical map
\begin{equation}
\label{ref-B.18-193}
\Hom(F,\Ascr)\r \Hom(P,\Ascr)\cong\Sigma \Ascr
\end{equation}
in $\Wscr^{-1}\Delta \Alg(\Oscr)$.
\begin{propositions} 
  \label{ref-B.8.1-194} Assume $\Ascr\in \Alg^+(\Pscr,\Oscr)$ and $F\in
  \Hscr(\Pscr)$ is such that $\Ext^i_{\Pscr_\ZZ}(\ZZ F_m, \Ascr^n)=0$
  for all $m\ge 0$, $n\in \ZZ$ and $i>0$. Then
\eqref{ref-B.18-193} is an isomorphism 
in $\Wscr^{-1} \Delta \Alg(\Oscr)$. 
\end{propositions}
\begin{proof} 
 We have to show that
\[
C^\ast(\Hom(F,\Ascr))\r C^\ast(\Hom(P,\Ascr))
\]
is a  quasi-isomorphism. By formula \eqref{ref-B.15-181} and our
vanishing hypotheses we obtain $C^\ast(\Hom(F,\Ascr))\cong
\RHom_{\Pscr_\ZZ}(\underline{\ZZ},\Ascr)$. It follows from Lemma
\ref{ref-B.7.11-191} that $C^\ast(\Hom(P,\Ascr))\cong \RHom_{\Pscr_\ZZ}(\underline{\ZZ},\Ascr)$ as well. 
\end{proof}
\subsection{{\v{C}ech} cohomology}
\label{ref-B.9-195}
In this section we discuss the important special case of {\v{C}ech} cohomology. 
Let $X$ be a topological space and let $\Uscr=\{U_i\mid i\in I\}$ be 
an  open covering of $X$. As usual we identify $U\in \Open(X)$ with the
representable sheaf $\Hom_{\Open(X)}(-,U)$. Then the unordered {\v{C}ech} covering of
$X$
is the simplicial sheaf on $X$ which in degree $m$ is equal to
\[
C(\Uscr)_m=\coprod_{i_0,\ldots,i_m} U_{i_0}\cap\cdots \cap U_{i_m}
\]
It is well-known and easy to see that this a hypercovering.

If given an ordering on $I$ we may also define 
\[
C^o(\Uscr)_m=\coprod_{i_0\le \ldots\le i_m} U_{i_0}\cap\cdots \cap U_{i_m}
\]

Note that the inclusion map
\[
C^o(\Uscr)\r C(\Uscr)
\]
is a map of simplicial sheaves. 

Let $\Ascr\in \Alg^+(\Pscr,\Oscr)$. The unordered and ordered {\v{C}ech}
complexes of $\Ascr$ are respectively defined as the cosimplicial
complexes of $\Oscr$-algebras
\begin{equation}
\label{ref-B.19-196}
\begin{split}
\Ch\nolimits(\Uscr,\Ascr)&=\Hom(C(\Uscr),\Ascr)\\
\Ch\nolimits^o(\Uscr,\Ascr)&=\Hom(C^o(\Uscr),\Ascr)
\end{split}
\end{equation}
\begin{lemmas} \label{ref-B.9.1-197} Assume that for all $m\ge 0$,
  $\{i_0,\ldots,i_m\}\subset I$, $j>0$ and $n\in \ZZ$ we have
\[
H^i(U_{i_0}\cap\cdots \cap U_{i_m},\Ascr^n)=0
\]
Then 
\[
\Ch(\Uscr,\Ascr)\cong \Ch\nolimits^o(\Uscr,\Ascr)\cong \Sigma(\Ascr)
\]
in $\Wscr^{-1} C^+(\Pscr_\ZZ)$.
\end{lemmas}
\begin{proof} Since $\Ch(\Uscr,\Ascr)$ is a hypercovering the
  isomorphism $\Ch(\Uscr,\Ascr)\cong \Sigma(\Ascr)$ follows from
  Proposition \ref{ref-B.8.1-194}. 

  The ordered {\v{C}ech} covering is not a hypercovering but nevertheless,
  by looking at stalks, it is easy to see that $C^\ast(\ZZ
  C^o(\Uscr))$ is a resolution of the constant sheaf $\ZZ_X$. From this we deduce that source and target of the map 
\[
C^\ast(\Ch\nolimits^o(\Uscr,\Ascr))\r \Ch(\Uscr,\Ascr)
\]
compute $\RHom(\ZZ_X,\Ascr)$. Hence it is a quasi-isomorphism. 
\end{proof}
\subsection{Relation to Hinich's construction}
\label{ref-B.10-198}
We now assume that $\Pscr=\Sh(\Cscr)$ for a small site $\Cscr$. 
A presheaf over $\Cscr$ is said to be \emph{semi-representable} if it
is a coproduct of representable presheaves. We will say that 
a simplicial presheaf $F$  is a presheaf-hypercovering if the associated simplicial sheaf $aF$ is
a hypercovering in the above sense. We will say that $F$
is a Verdier hypercovering if each $F_n$ 
is semi-representable.  We denote the corresponding categories by
$\Hscr^{\Pre}(\Cscr)$ and $\Hscr^{\text{V}}(\Cscr)$.

If $U\in \Cscr$
then a simplicial presheaf $F$ with an augmentation $F\r U$ will be
called a Verdier-hypercovering of $U$ if $F$ is a Verdier-hypercovering of $U$ in
the site $\Cscr/U$. 

Following \cite{hinich2} we say that a complex of presheaves is fibrant if
if for any  $U\in \Cscr$ and for any Verdier-hypercovering $F\r U$ we have that $M(U)\r
C^\ast(\Hom(F,M))$ is a quasi-isomorphism.

Hinich proves under some hypotheses on $\Oscr$ (which hold if $\Oscr$
is $k$-linear over a field of characterstic zero) that for any
presheaf of $\Oscr$-algebras $\Ascr$ there is a map of presheaves of
$\Oscr$-algebras $\Ascr\r \Ascr'$ with $\Ascr'$ fibrant which is a
quasi-isomorphism after sheaffication. The derived global sections of
$\Ascr$ are then given by $\Hom(F,\Ascr')^{TS}$ for a Verdier
hypercovering $F$ of $e$. If $e$ itself is in $\Cscr$ then we
may consider it as its own hypercovering and in this case we may
dispense with the Thom-Sullivan normalization. I.e.\ we may define the derived
global sections of $\Ascr$ as $\Ascr'(e)$. 

We will show that in case $\Ascr$ is a sheaf of $\Oscr$-algebras with
left bounded grading this yields the same result as our construction.
Mimicking the proof of Proposition \ref{ref-B.6.2-173} we may produce
a pro-object $P=(P_\alpha)_\alpha$ in $\Hscr^{V}(\Cscr)$ mapping to
the hypercovering $F$ such that any diagram of solid arrows
\[
\xymatrix{
&P\ar@{.>}[dl]\ar[d]\\
H_1\ar[r] & H_2
}
\]
with $H_1$, $H_2$ in $\Hscr^{\Pre}(\Cscr)$ can be factored like the
dotted arrow (up to homotopy).  It follows in particular that $aP$ is
a homotopy projective object in $\Pro \Hscr(\Pscr)$.  Hence we need to
prove that $\Hom(P,\Ascr)$ is weaky equivalent to $\Hom(F,\Ascr')$.

We have now maps
\[
\Hom(P,\Ascr)\r \Hom(P,\Ascr')\l \Hom(F,\Ascr')
\]
and it it is sufficient prove that these are weak equivalences.  By a
suitable analogue of Lemma \ref{ref-B.7.10-190} the first maps is a
weak equivalence.  By an analogue of Lemma \ref{ref-B.7.9-189} we
have that $\Hom(P,\Ascr')$ is weakly equivalent to $\injlim_\alpha
\Hom(P_\alpha,\Ascr')$. Hence it is sufficient to show that
\[
\Hom(P_\alpha,\Ascr')\l \Hom(F,\Ascr')
\]
is a weak equivalence. This follows from the fact that
$\Hom(F,\Ascr')$ is, up to weak equivalence, independent of the
Verdier hypercovering $F$. See \cite[\S1.4.3]{hinich2}.

\def\cprime{$'$} \def\cprime{$'$}
\ifx\undefined\bysame
\newcommand{\bysame}{\leavevmode\hbox to3em{\hrulefill}\,}
\fi

\end{document}